\documentclass[10pt]{amsart}
\usepackage{amssymb}
\usepackage[all]{xy}

\newtheorem{prop}{Proposition}[section]
\newtheorem{prop:def}{Proposition-Definition}[section]

\newtheorem{lemma}{Lemma}[section]

\newtheorem{thm}{Theorem}[section]

\theoremstyle{remark}

\newtheorem{remark}{Remark}

\begin{document}

\newcommand{\nc}{\newcommand} \nc{\on}{\operatorname}

\nc{\pa}{\partial}

\nc{\cA}{{\mathcal A}} \nc{\cB}{{\cal B}}\nc{\cC}{{\mathcal C}} 
\nc{\cD}{{\mathcal D}} 
\nc{\cE}{{\mathcal E}} \nc{\cG}{{\mathcal G}}\nc{\cH}{{\cal H}} 
\nc{\cI}{{\cal I}} \nc{\cJ}{{\cal J}}\nc{\cK}{{\cal K}} 
\nc{\cL}{{\mathcal L}} \nc{\cR}{{\cal R}} \nc{\cS}{{\mathcal S}}   
\nc{\cV}{{\cal V}} \nc{\cX}{{\cal X}}

\nc{\¦}{{|}}

\nc{\sh}{\on{sh}}\nc{\Id}{\on{Id}}\nc{\Diff}{\on{Diff}}
\nc{\Perm}{\on{Perm}}\nc{\conc}{\on{conc}}\nc{\Alt}{\on{Alt}}
\nc{\ad}{\on{ad}}\nc{\Der}{\on{Der}}\nc{\End}{\on{End}}
\nc{\no}{\on{no\ }} \nc{\res}{\on{res}}\nc{\ddiv}{\on{div}}
\nc{\Sh}{\on{Sh}} \nc{\card}{\on{card}}\nc{\dimm}{\on{dim}}
\nc{\Sym}{\on{Sym}} \nc{\Jac}{\on{Jac}}\nc{\Ker}{\on{Ker}}
\nc{\Spec}{\on{Spec}}\nc{\Cl}{\on{Cl}}
\nc{\Imm}{\on{Im}}\nc{\limm}{\lim}\nc{\Ad}{\on{Ad}}
\nc{\ev}{\on{ev}} \nc{\Hol}{\on{Hol}}\nc{\Det}{\on{Det}}
\nc{\Bun}{\on{Bun}}\nc{\diag}{\on{diag}}\nc{\pr}{\on{pr}} 
\nc{\Span}{\on{Span}}\nc{\Comp}{\on{Comp}}\nc{\Part}{\on{Part}}
\nc{\tensor}{\on{tensor}}\nc{\ind}{\on{ind}}\nc{\id}{\on{id}}
\nc{\Hom}{\on{Hom}}\nc{\Quant}{\on{Quant}} \nc{\Dequant}{\on{Dequant}}\nc{\Def}{\on{Def}}
\nc{\AutLBA}{\on{AutLBA}}\nc{\AutQUE}{\on{AutQUE}}
\nc{\LBA}{{\on{LBA}}}\nc{\Aut}{\on{Aut}}\nc{\QUE}{{\on{QUE}}}
\nc{\Lyn}{\on{Lyn}}\nc{\Cof}{\on{Cof}}\nc{\LCA}{\underline{\on{LCA}}}\nc{\FLBA}{\on{FLBA}}
\nc{\LA}{\on{LA}}\nc{\FLA}{\on{FLA}}\nc{\EK}{\on{EK}}
\nc{\class}{\on{class}}\nc{\br}{\on{br}}\nc{\co}{\on{co}}
\nc{\Prim}{\on{Prim}}\nc{\ren}{\on{ren}}\nc{\lbr}{\on{lbr}}
\nc{\SC}{\on{SC}}\nc{\ol}{\overline}\nc{\FL}{\on{FL}}
\nc{\FA}{\on{FA}}\nc{\alg}{{\on{alg}}}\nc{\KZ}{\on{KZ}}
\nc{\op}{{\on{op}}}\nc{\cop}{{\on{cop}}}
\nc{\Inn}{\on{Inn}}\nc{\OutDer}{\on{OutDer}}
\nc{\inv}{{\on{inv}}}\nc{\gr}{{\on{gr}}}\nc{\Lie}{{\on{Lie}}}
\nc{\Out}{{\on{Out}}}\nc{\univ}{{\on{univ}}}
\nc{\GT}{{\on{GT}}}\nc{\GRT}{{\on{GRT}}}
\nc{\GS}{{\on{GS}}} 
\nc{\restr}{{\on{restr}}}\nc{\lin}{{\on{lin}}}
\nc{\mult}{{\on{mult}}}\nc{\qt}{{\on{qt}}}
\nc{\compl}{{\on{compl}}} \nc{\Trees}{{\on{Trees}}}
\nc{\Ord}{{\on{Ord}}} \nc{\norm}{{\on{norm\ ord}}}
\nc{\fd}{{\on{fd}}}\nc{\Maps}{{\on{Maps}}}
\nc{\CYBE}{{\on{CYBE}}} \nc{\CY}{{\on{C}}} \nc{\can}{{\on{can}}}
\nc{\Alg}{{\underline{\on{Alg}}}} \nc{\Poisson}{{\underline{\on{Poisson}}}}
\nc{\Coalg}{{\underline{\on{Coalg}}}} \nc{\Vect}{{\underline{\on{Vect}}}}
\nc{\UE}{{\on{UE}}}\nc{\Bialg}{{\underline{\on{Bialg}}}}
\nc{\QTA}{{\underline{\on{QTA}}}}\nc{\QTQUE}{{\underline{\on{QTQUE}}}}
\nc{\QTLBA}{{\underline{\on{QTLBA}}}}\nc{\TA}{{\underline{\on{TA}}}}
\nc{\MT}{{\underline{\on{MT}}}}\nc{\QYBE}{{\underline{\on{QYBE}}}}
\nc{\QTBialg}{{\underline{\on{QTBialg}}}}
\nc{\qcocomm}{{\on{quasi-cocomm}}}
\nc{\qcomm}{{\on{quasi-comm}}}\nc{\cocomm}{{\on{cocomm}}}
\nc{\dual}{{\on{dual}}} \nc{\sym}{{\on{sym}}}
\nc{\cycl}{{\on{cycl}}} \nc{\Prop}{{\underline{\on{Prop}}}}
\nc{\uV}{{\underline{V}}} \nc{\uR}{{\underline{R}}}
\nc{\Mor}{{\on{Mor}}}\nc{\Ob}{{\on{Ob}}}
\nc{\Bil}{{\on{Bil}}} \nc{\invt}{{\on{invt}}}

\nc{\al}{\alpha}\nc{\de}{\delta}
\nc{\eps}{\epsilon}\nc{\la}{{\lambda}}
\nc{\si}{\sigma}\nc{\z}{\zeta}

\nc{\La}{\Lambda}

\nc{\ve}{\varepsilon} \nc{\vp}{\varphi} 

\nc{\AAA}{{\mathbb A}}\nc{\BB}{{\mathbb B}}
\nc{\CC}{{\mathbb C}}\nc{\ZZ}{{\mathbb Z}} 
\nc{\QQ}{{\mathbb Q}} \nc{\NN}{{\mathbb N}}\nc{\VV}{{\mathbb V}} 
\nc{\KK}{{\mathbb K}} 

\nc{\ff}{{\mathbf f}}\nc{\bg}{{\mathbf g}}
\nc{\ii}{{\mathbf i}}\nc{\kk}{{\mathbf k}}
\nc{\bl}{{\mathbf l}}\nc{\zz}{{\mathbf z}} 
\nc{\pp}{{\mathbf p}}\nc{\qq}{{\mathbf q}} 

\nc{\cF}{{\cal F}}\nc{\cM}{{\cal M}}\nc{\cO}{{\cal O}}
\nc{\cT}{{\cal T}}\nc{\cW}{{\cal W}}\nc{\cP}{{\mathcal P}}

\nc{\Assoc}{{\mathbf Assoc}}
\nc{\ul}{\underline}
\def\Sha{{\mathop{\scriptstyle\amalg\!\hspace{-1.8pt}\amalg}}}

\nc{\ub}{{\underline{b}}}
\nc{\uk}{{\underline{k}}} 
\nc{\un}{{\underline{n}}} \nc{\um}{{\underline{m}}}
\nc{\up}{{\underline{p}}}\nc{\uq}{{\underline{q}}}
\nc{\ur}{{\underline{r}}}
\nc{\us}{{\underline{s}}}\nc{\ut}{{\underline{t}}}
\nc{\uw}{{\underline{w}}}
\nc{\uz}{{\underline{z}}}
\nc{\ual}{{\underline{\alpha}}}\nc{\ualpha}{{\underline{\alpha}}}
\nc{\ubeta}{{\underline{\beta}}}\nc{\ugamma}{{\underline{\gamma}}}
\nc{\ueps}{{\underline{\epsilon}}}\nc{\ueta}{{\underline{\eta}}}
\nc{\uzeta}{{\underline{\zeta}}}\nc{\ula}{{\underline{\lambda}}}
\nc{\umu}{{\underline{\mu}}}\nc{\unu}{{\underline{\nu}}}
\nc{\usigma}{{\underline{\sigma}}}\nc{\utau}{{\underline{\tau}}}
\nc{\uI}{{\underline{I}}}\nc{\uJ}{{\underline{J}}}
\nc{\uK}{{\underline{K}}}\nc{\uM}{{\underline{M}}}
\nc{\uN}{{\underline{N}}}

\nc{\A}{{\mathfrak a}}
\renewcommand{\t}{{\mathfrak t}}\nc{\G}{{\mathfrak g}}
\nc{\x}{{\mathfrak x}}
\nc{\B}{{\mathfrak b}} \nc{\C}{{\mathfrak c}} 
\nc{\D}{{\mathfrak d}} \nc{\HH}{{\mathfrak h}}
\nc{\iii}{{\mathfrak i}}\nc{\mm}{{\mathfrak m}}\nc{\N}{{\mathfrak n}} 
\nc{\ttt}{{\mathfrak{t}}}\nc{\U}{{\mathfrak u}}\nc{\V}{{\mathfrak v}}
\nc{\grt}{{\mathfrak{grt}}}\nc{\gt}{{\mathfrak gt}}
\nc{\SL}{{\mathfrak{sl}}}\nc{\out}{{\mathfrak{out}}}

\nc{\SG}{{\mathfrak S}}

\nc{\wt}{\widetilde} \nc{\wh}{\widehat}

\nc{\bn}{\begin{equation}}\nc{\en}{\end{equation}} \nc{\td}{\tilde}

\title{Quantization of quasi-Lie bialgebras}

\begin{abstract}
We construct quantization functors of quasi-Lie bialgebras. We establish a bijection 
between this set of quantization functors, modulo equivalence and twist equivalence, 
and the set of quantization functors of Lie bialgebras, modulo equivalence. This is 
based on the acyclicity of the kernel of the natural morphism from  
the universal deformation complex of quasi-Lie bialgebras to that of Lie bialgebras. 
The proof of this acyclicity consists in several steps, ending up in the acyclicity 
of a complex related to free Lie algebras,  namely, the universal version of the 
Lie algebra cohomology complex of a Lie algebra in its enveloping algebra, viewed as
the left regular module. Using the same arguments, we also prove the compatibility of 
quantization functors of quasi-Lie bialgebras with twists, which allows us to recover
our earlier results on compatibility of quantization functors 
with twists in the case of Lie bialgebras. 
\end{abstract}

\author{Benjamin Enriquez}
\address{B.E.: IRMA (CNRS) et Universit\'e Louis Pasteur, 
7, rue Ren\'e Descartes, F-67084 Strasbourg, France}
\email{enriquez@math.u-strasbg.fr}

\author{Gilles Halbout}
\address{G.H.: D\'ept. Maths., Universit\'e Montpellier 2,  Place E. Bataillon, 
F-34095 Montpellier, France}
\email{ghalbout@math.univ-montp2.fr}

\maketitle

Let $\kk$ be a field of characteristic $0$. Unless specified otherwise, 
``Lie algebra'', ``vector space'', etc.,  means ``Lie algebra over $\kk$", etc. 

\section{Introduction and main results}

The main result of this paper is the construction of quantization functors for 
quasi-Lie bialgebras. This problem was posed in \cite{Dr:unsolved}, Section 5; 
let us recall its formulation (\cite{Dr:QH}). A quasi-bialgebra is a 6-uple
$(A,m,\Delta,\Phi,\eps,\eta)$, where $(A,m,\eta)$ is an algebra with unit, 
$\Delta:A\to A^{\otimes 2}$ is an algebra morphism with counit $\eps$, 
and $\Phi\in A^{\otimes 3}$
satisfies the identity $(\on{id}\otimes\Delta)(\Delta(a)) = 
\Phi (\Delta\otimes\on{id})(\Delta(a))\Phi^{-1}$ and the pentagon identity. 
Examples of quasi-bialgebras are $A_{0}=U(\A)$, where $\A$ is a 
Lie algebra, equipped with its bialgebra structure and $\Phi=1$. A deformation of 
$A_{0}$ in the category of topologically free $\kk[[\hbar]]$-modules (i.e., 
a QUE quasi-bialgebra) gives rise to 
a quasi-Lie bialgebra $(\A,\mu,\delta,\varphi)$, the classical limit of $A$. 
A quantization functor of quasi-Lie bialgebras is a section of the classical limit
functor $\{$QUE quasi-bialgebras$\}\to\{$quasi-Lie bialgebras$\}$. The paper
\cite{Dr:QH} contains examples of quasi-Lie bialgebras, whose quantization is 
not explicitly known (e.g., the quasi-Lie bialgebras from p. 1437, in the 
non-split case).

As in the case of Lie bialgebras, the problem of construction of quantization 
functors of quasi-Lie bialgebras can be rephrased in the language of props (\cite{McL}). 
In Section \ref{sect:qlba}, we introduce a complete prop ${\bf QLBA}$ of quasi-Lie 
bialgebras and define the notion of a quantized symmetric quasi-bialgebra (QSQB) 
in ${\bf QLBA}$.  A QSQB in ${\bf QLBA}$ (which we will also call a 
universal quantization functor of quasi-Lie bialgebras)
then gives rise to a quantization functor as above (it also allows to construct
quantizations of quasi-Lie bialgebras in symmetric monoidal categories, more general than 
that of vector spaces, like super-vector spaces, d.g.vector spaces, etc.). We introduce
the notions of equivalence and twist equivalence on the set of QSQB's
in ${\bf QLBA}$. We similarly introduce the notion of a QSB (quantized symmetric 
bialgebra) in the prop ${\bf LBA}$ of Lie bialgebras, which is the same as the universal 
quantization functors from \cite{EK2}; there is a notion of equivalence for these QSB's. 
Our main result is (Theorem \ref{thm:main}): 

\begin{thm} \label{thm:main0}
The map $\{$universal quantization functors of quasi-Lie bialgebras$\}/($equivalence, 
twist equivalence$)
\to \{$universal quantization functors of Lie bialgebras$\}/($equivalence$)$ is a bijection. 
\end{thm}

Together with the results of \cite{EK1,EK2}, where is constructed a map 
$\{$associators over $\kk\}\to \{$universal quantization functors of Lie bialgebras$\}$, 
and of \cite{Dr:Gal} on the existence of associators over $\kk$, 
this result implies the existence of universal quantization functors for quasi-Lie bialgebras. 

Let us explain the idea of the proof of Theorem \ref{thm:main0}. 
According to deformation theory, there are complexes
$C_{\on{QLBA}}$ and $C_{\on{LBA}}$, equipped with gradings, whose second and 
first cohomology groups respectively contain the obstruction to lifting a quantization from 
degree $n$ to degree $n+1$, and parametrize such liftings. This viewpoint in not
used in the quantization of Lie bialgebras, since the groups $H^{i}_{\on{LBA}}$
are not known. (In the same way, it is not known how to construct associators using
deformation theory, see Remark 2, p. 854 in \cite{Dr:Gal}.) However, this viewpoint 
can be used in our context. Namely, we will prove: 

\begin{thm} \label{thm:isom0}
The canonical map $H^{i}_{\on{QLBA}}
\to H^{i}_{\on{LBA}}$ is an isomorphism for any $i\geq 0$.
\end{thm} 

This immediately implies our main result (see Section \ref{pf:main}). 

Let us give some details of the proof of Theorem \ref{thm:isom0}. Our aim is to 
prove that the relative complex $\on{Ker}(C_{\on{QLBA}} \to C_{\on{LBA}})$
is acyclic. We introduce a filtration of the prop $\on{QLBA}$ by the powers of
an ideal $\langle\varphi\rangle$. Our first main result result is Theorem 
\ref{thm:QLBA}, 
which gives an isomorphism $\on{gr}\on{QLBA}\simeq \on{LBA}_{\alpha}$
of the associated graded prop with an explicitly presented prop. For this, one constructs a 
morphism $\on{LBA}_{\alpha}\to \on{gr}\on{QLBA}$, which is clearly surjective
(here $\on{LBA}_\alpha$ is an explicitly presented prop); 
to prove its injectivity, we use the existence of ``many'' quasi-Lie bialgebras, namely, 
the classical twists of Lie bialgebras of the form $F(\C)$ (where $\C$ is a Lie coalgebra)
by an element $r\in \wedge^{2}(\C)$ (in the same way, the existence of the Lie bialgebras
$F(\C)$ is the argument underlying the structure theorem for the prop $\on{LBA}$, 
see \cite{Enr:univ,Pos}). 

The associated graded (for the filtration of $\on{QLBA}$) of the relative complex
is then the positive degree part of the complex $C_{\on{LBA}_{\alpha}}$. To 
prove that it is acyclic, we prove that its lines are. These lines are of the form
$0\to  \on{LBA}_{\kappa}({\bf 1},\wedge^{q})
\to \on{LBA}_{\kappa}({\bf id},\wedge^{q}) \to \on{LBA}_{\kappa}(\wedge^{2},
\wedge^{q})\to ...$, where $\on{LBA}_{\kappa}(X,Y) = \on{Coker}(\on{LBA}(
D\otimes X,Y)\to \on{LBA}(C\otimes X,Y))$ and the map corresponds to 
$\kappa\in\on{LCA}(C,D)$. Here $C,D$ are sums of Schur functors of positive degree, 
and the differential is the universal version of the differential of Lie algebra
cohomology. 

The proof of the acyclicity of this complex (Theorem \ref{thm:isom}, proof in 
Section \ref{sec:comparison}) involves several reductions. 
We first show that in this complex, the spaces of cochains may be modified as
follows: $\on{LBA}_{\kappa}(\wedge^{p},\wedge^{q})$ is replaced by 
$\on{LBA}_{Z}(\wedge^{p},\wedge^{q}) = \on{LBA}(Z\otimes \wedge^{p},
\wedge^{q})$, where $Z$ 
is an irreducible Schur functor, and the space of cochains is reduced to the sum of 
its components, where the ``intermediate
Schur functor between $Z$ and $\wedge^{q}$'' is $Z$ (this notion is based 
on the structure theorem of $\on{LBA}$, see Proposition \ref{prop:PBW}; we say that 
the intermediate Schur functor between $X_{i}$ and $Y_{j}$ in the summand appearing 
in the r.h.s. of (\ref{str:LBA}) is $Z_{ij}$). 
We next introduce a filtration on the complex, viewing $\wedge^{p}$ as a subobject of 
${\bf id}^{\otimes p}$ and counting the number of intermediate Schur functor between 
the $p$ factors
${\bf id}$ and $\wedge^{q}$ which equal ${\bf id}$. We identify the associated graded
complex with a subcomplex of $0\to \on{LBA}(Z\otimes{\bf 1}\otimes 
{\bf id}^{\otimes p''},\wedge^{q})\to  ...\to  \on{LBA}(Z\otimes 
\wedge^{p'} \otimes 
{\bf id}^{\otimes p''},\wedge^{q})\to...$, where the differential involves 
Lie brackets between the components of ${\bf id}^{\otimes p'}\supset \wedge^{p'}$ 
and of these components with $Z$, formed by the sums of components, 
where the intermediate Schur functor between a component ${\bf id}$ of
${\bf id}^{\otimes p'}$ (resp., of ${\bf id}^{\otimes p''}$)
and $\wedge^{q}$ is ${\bf id}$ (resp., has degree $>1$) and antisymmetric w.r.t.
$\SG_{p''}$. 
This subcomplex decomposes according to the intermediate Schur functors
between the factors of ${\bf id}^{\otimes p''}$ and $\wedge^{q}$, 
and these subcomplexes are obtained from the complexes 
$C^{\bullet}_{\ul Z''} = (0\to \on{LA}(Z\otimes {\bf 1}\otimes 
(\otimes_{i=1}^{p''}Z''_{i}),\wedge^{q})\to...\to \on{LA}(Z\otimes \wedge^{p'}\otimes 
(\otimes_{i=1}^{p''}Z''_{i}),\wedge^{q})\to...)$ with the same differentials
(the $Z''_{i}$ are irreducible Schur functors of degree $>1$) by taking tensor products
with vector spaces $\on{LCA}({\bf id},Z''_{i})$
and taking antiinvariants under $\SG_{p''}$. We therefore have to show
the acyclicity of the complexes $C^{\bullet}_{\ul Z''}$. 

For this, we show that $\wedge^{q}$ may be replaced by ${\bf id}^{\otimes q}$,  
$Z$ by ${\bf id}^{\otimes z}$, and $(\otimes_{i=1}^{p''}Z''_{i})$ 
by ${\bf id}^{\otimes N}$, and express the corresponding complex as 
a sum of tensor products of complexes, which reduces the problem to a complex
$0\to \on{LA}({\bf id}^{\otimes z}\otimes {\bf 1}\otimes{\bf id}^{\otimes N},
{\bf id})\to...\to \on{LA}({\bf id}^{\otimes z}\otimes \wedge^{p'}
\otimes{\bf id}^{\otimes N}, {\bf id})\to...$. The spaces of chains are now spaces of
multilinear Lie polynomials. Using Dynkin's correspondence between free Lie and 
free associative polynomials, we identify the complex with a complex 
$\cA_{z,N,1}^{\bullet}$, 
defined in terms of associative polynomials, which we decompose as a direct sum 
$\oplus_{\sigma}\cA^{\bullet}_{\sigma}$ of subcomplexes, indexed by permutations. 
We next identify each summand $\cA^{\bullet}_{\sigma}$ with a tensor product of
``elementary'' complexes. These complexes $\cE_{\eps,\eps'}^{\bullet}$ ($\eps,\eps'
\in \{0,1\}$)
are 1-dimensional in each degree, and are universal versions of the complexes computing 
the Lie algebra cohomology of a Lie algebra $\G$ in $U(\G)$, 
equipped with one of its trivial, adjoint, 
left or right $\G$-module structures. We show that two of these complexes are 
acyclic, using the PBW filtration of free associative algebras (when 
$\G$ is a finite dimensional Lie algebra, the corresponding complexes have
1-dimensional cohomology, concentrated in degree $\on{dim}\G$). 
As $\cE_{0,1}^{\bullet}$ necessarily enters the tensor product
decomposition of each subcomplex $\cA^{\bullet}_{\sigma}$, each of the 
$\cA^{\bullet}_{\sigma}$ is acyclic, which implies that $\cA_{z,N,1}^{\bullet}$ is
acyclic.  

In the final section of the paper, we apply Theorem \ref{thm:isom} for proving that 
quantization functors of quasi-Lie bialgebras are compatible with twists (Proposition 
\ref{prop:twists}). This allows us to generalize our earlier results (\cite{EH}) 
on compatibility of quantization functors of Lie bialgebras with twists, 
see Proposition \ref{prop:5:2} (in \cite{EH}, this result was established 
for Etingof-Kazhdan 
quantization functors, while Proposition \ref{prop:5:2} applies to any quantization 
functor of Lie bialgebras).

\section{Quantization of (quasi)Lie bialgebras} \label{sect:qlba}

In this section, we recall the general formalism of props and its relation 
with the quantization problems of (quasi)Lie bialgebras.  In particular, 
we show that this formalism allows to recover biquantization results of  
\cite{KT}. We also formulate our main result (Theorem \ref{thm:main}) 
and explain the strategy of its proof.

\subsection{Props}

Recall the definitions of the Schur categories $\on{Sch}$ and ${\bf Sch}$ 
(\cite{EH}). These are braided symmetric tensor categories, defined as 
follows. The objects of $\on{Sch}$ are Schur functors, i.e., finitely supported families
$X = (X_{\rho})_{\rho}$ of finite dimensional vector spaces, where 
$\rho\in\sqcup_{n\geq 0}\wh\SG_{n}$ ($\rho$ is therefore a pair $(n,\pi_{\rho})$, 
where $n\geq 0$ and $\pi_{\rho}$ is an irreducible representation of $\SG_{n}$; $n$ 
is called the degree of $\rho$; by
convention, $\SG_{0}$ is the trivial group). The set of morphisms from $X$
to $Y$ is\footnote{For $U,V$ finite dimensional vector spaces, 
$\on{Vect}(U,V)$ is the set of morphisms $U\to V$.} $\on{Sch}(X,Y):= \oplus_{\rho} 
\on{Vect}(X_{\rho},Y_{\rho})$. The direct sum of objects is 
$X\oplus Y = (X_{\rho}\oplus Y_{\rho})_{\rho\in\sqcup_{n\geq 0}
\wh\SG_{n}}$. Their tensor product is $X\otimes Y = (\oplus_{\rho',\rho''}
M^{\rho}_{\rho'\rho''} \otimes X_{\rho'}\otimes Y_{\rho''})_{\rho}$, where
for $\rho\in \wh\SG_{n}$, $\rho'\in\wh\SG_{n'}$, $\rho''\in\wh\SG_{n''}$, 
then $M^{\rho}_{\rho'\rho''} = [\pi_{\rho}:\on{Ind}^{\SG_{n}}_{
\SG_{n'}\times\SG_{n''}}(\pi_{\rho'}\otimes\pi_{\rho''})]$ if 
$n=n'+n''$ and $0$ otherwise. ${\bf Sch}$ is defined similarly, dropping the 
condition that $X$ is finitely supported. 

An object $X$ of $\on{Sch}$ or ${\bf Sch}$ is called homogeneous 
of degree $n$ iff $X_{\rho}=0$ if the degree of $\rho$ is $\neq n$. 
If $X$ is homogeneous, we denote by $|X|$ its degree. 

We have a bijection $\on{Irr(Sch)}\simeq \sqcup_{n\geq 0}\wh\SG_{n}$, 
where $\on{Irr(Sch)}$ is the set of equivalence classes of irreductible objects in 
$\on{Sch}$. The unit object of $\on{Sch}$ is ${\bf 1}$, which corresponds to the 
element of $\wh\SG_{0}$. We also define ${\bf id},S^{p},\wedge^{p}$
as the objects corresponding to: the element of $\wh\SG_{1}$, the trivial 
and the signature character of $\SG_{p}$. We set $T_{p}:= {\bf id}^{\otimes p}$
and $S:= \oplus_{p\geq 0}S^{p}\in \on{Ob}({\bf Sch})$. 

Recall that a prop (resp., a ${\bf Sch}$-prop) is an additive symmetric monoidal 
category $\cC$, equipped with a tensor functor $\on{Sch}\to \cC$
(resp., ${\bf Sch}\to\cC$), which is the identity on objects. 

A prop morphism $\cC\to \cD$ is a tensor functor, such that the functors
$\on{Sch}\to\cC\to\cD$ and $\on{Sch}\to\cD$ coincide. An ideal $I$ 
of the prop $\cC$ is an assignment $(X,Y)\mapsto I(X,Y)$, such that 
$(X,Y)\mapsto \cC/I(X,Y)$ is a monoidal category. $\cC/I$ is then the 
quotient prop. If $f_{\alpha}\in \cC(X_{\alpha},Y_{\alpha})$ are 
morphisms in $\cC$, then the ideal $\langle f_{\alpha}\rangle$ is the 
smallest ideal $I$ in $\cC$ such that $f_{\alpha}\in I(X_{\alpha},Y_{\alpha})$. 

Props may be presented by generators and relations. 
If $\ul V = (V_{n,m})_{n,m\geq 0}$ is a collection of vector spaces, there is a 
unique prop $\on{Free}(\ul V)$, such that for any prop $\cC$, we have a 
bijection $\prod_{n,m}\on{Vect}(V_{n,m},\cC(T_{n},T_{m})) \simeq 
\on{Prop}(\on{Free}(\ul V),\cC)$ (where $\on{Prop}$ denotes the set of 
prop morphisms). If $\ul R = (R_{n,m})_{n,m\geq 0}$ is a collection of 
subspaces of the $\on{Free}(\ul V)(T_{n},T_{m})$, then the ideal generated
by $\ul V$ with relations $\ul R$ is the quotient prop 
$\on{Free}(\ul V)/\langle\ul R\rangle$. 

We say that $\cC$ is a topological prop if for any $X,Y\in \on{Ob(Sch)}$
we have a filtration $\cC(X,Y) = \cC^{\geq 0}(X,Y)\supset \cC^{\geq 1}(X,Y)
\supset...$,  complete and separated, and compatible with the 
prop operations; and if $\cC(X,Y) = \cC^{\geq v(|X|,|Y|)}(X,Y)$
for any homogeneous $X,Y$, where $v(x,y)\to \infty$ when 
$x\to\infty$, $y$ being fixed, or $y\to \infty$, $x$ being
fixed.  Such a $\cC$ gives rise to a ${\bf Sch}$-prop $\hat\cC$, 
given by $\hat\cC(X,Y) = \hat\oplus_{i,j}\cC(X_{i},Y_{j})$, where $X = \oplus_{i}X_{i}$, 
$Y = \oplus_{j}Y_{j}$ are the homogeneous decompositions of $X,Y\in
\on{Ob}({\bf Sch})$ ($\hat\oplus$ is the direct product); $\hat\cC$ is then equipped 
with a complete and separated filtration, compatible with the prop operations.

\subsection{Quantization functors of Lie bialgebras}

If $\cC$ is a topological prop, the associated graded prop is $\on{gr}(\cC)$ 
with $\on{gr}(\cC)(X,Y) := \oplus_{i\geq 0}\cC^{\geq i}/\cC^{\geq i+1}(X,Y)$ for 
$X,Y\in\on{Ob}(\on{Sch})$. 
We denote by $\hat{\on{gr}}(\cC)$ its degree completion. 

A {\it quantized symmetric bialgebra} (QSB) in $\hat\cC$ is a bialgebra structure on $S$ 
in $\hat\cC$, i.e., morphisms 
$m_{\cC}\in \hat\cC(S^{\otimes 2},S)$, $\Delta_{\cC}\in \hat\cC(S,S^{\otimes 2})$, 
$\eps_{\cC}\in \hat\cC(S,{\bf 1})$, 
$\eta_{\cC}\in \hat\cC({\bf 1},S)$, satisfying the bialgebra relations, and whose reductions
mod $\hat\cC^{\geq 1}$ coincide with the standard (commutative, cocommutative)
bialgebra structure on $S$, induced by the morphism ${\bf Sch}\to \hat\cC$.

A QSB in $\hat\cC$ gives rise to a Lie bialgebra 
structure\footnote{If ${\mathcal X}$ is a symmetric monoidal category and $X$ 
is an object of ${\mathcal X}$, recall that a Lie bialgebra 
(resp., algebra, Lie algebra, bialgebra...) structure on $X$ is a pair 
of morphisms $\mu_X\in {\mathcal X}(X^{\otimes 2},X)$ and $\delta_X\in 
{\mathcal X}(X,X^{\otimes 2})$ (resp., $m_X\in {\mathcal X}(X^{\otimes 2},X)$, 
etc.), 
satisfying the axioms of Lie bialgebra (resp., algebra, etc.).} 
on ${\bf id}$ in $\on{gr}(\cC)$, 
with morphisms of degree $1$ (its classical limit), as follows: 
one shows that there exist unique morphisms $\mu_{\cC}\in \on{gr}^{1}(\cC)(\wedge^{2},
{\bf id})$ and $\delta_{\cC}\in \on{gr}^{1}(\cC)({\bf id},\wedge^{2})$, such that 
$(m_{\cC} \circ \on{Alt}_{2}\circ inj_{1}^{\otimes 2}$ mod $\hat\cC^{\geq 2})
= inj_{1}\circ \mu_{\cC}$
and $(\on{Alt}_{2}\circ \Delta_{\cC}\circ inj_{1}$ mod $\hat\cC^{\geq 2})
= inj_{1}^{\otimes 2}\circ\delta_{\cC}$ 
(here $\on{Alt}_{2}:S^{\otimes 2}\to S^{\otimes 2}$ is the antisymmetrization 
(without factor 1/2), we identify 
$\cC(\wedge^{p},\wedge^{q})$ with a subspace of $\cC({\bf id}^{\otimes p},
{\bf id}^{\otimes q})$ and $inj_{1}:{\bf id}\to S$ is the canonical morphism); 
$(\mu_{\cC},\delta_{\cC})$ then obey the Lie bialgebra relations. 
We say that the QSB $(m_{\cC},\Delta_{\cC},\eps_{\cC},\eta_{\cC})$ 
is a quantization of $(\mu_{\cC},\delta_{\cC})$. 

Let $\hat\cC(S,S)_{1}$ be the preimage of $\on{id}_{S}$ under $\hat\cC(S,S)\to 
\hat\cC/\hat\cC^{\geq 1}(S,S)$; this is a group under composition. 
This group acts on the set of QSB's by $i_{\cC}*(m_{\cC},\Delta_{\cC},
\eps_{\cC},\eta_{\cC}):= 
(i_{\cC}\circ m_{\cC}\circ (i_{\cC}^{\otimes 2})^{-1},i_{\cC}^{\otimes 2}
\circ\Delta_{\cC} \circ i_{\cC}^{-1},\eps_{\cC}\circ i_{\cC}^{-1},
i_{\cC}\circ\eta_{\cC})$. Two QSB's related 
by this group action are called equivalent. Equivalent QSB's have the same classical limit. 

Recall that $\on{LBA}$ is the prop with generators $\mu\in\on{LBA}(\wedge^{2},
{\bf id})$, $\delta\in\on{LBA}({\bf id},\wedge^{2})$ and relations 
$$
\mu\circ (\mu\otimes\on{id}_{{\bf id}}) \circ \on{Alt}_{3}=0, \quad 
\on{Alt}_{3}\circ (\delta\otimes\on{id}_{{\bf id}})\circ\delta=0, \quad
\delta\circ\mu = \on{Alt}_{2} \circ (\mu\otimes\on{id}_{{\bf id}})
\circ (\on{id}_{{\bf id}}\otimes\delta) \circ \on{Alt}_{2} 
$$
(recall that $\mu,\delta$ are identified with morphisms in $\on{LBA}
({\bf id}^{\otimes 2},{\bf id})$ and $\on{LBA}({\bf id},{\bf id}^{\otimes 2})$). 
Then ${\bf id}$ is a Lie bialgebra in $\on{LBA}$, and it is an initial object in 
the category of props equipped with a Lie bialgebra structure on ${\bf id}$. 

$\on{LBA}$ is graded by $\NN^{2}$, with $\mu,\delta$ of degrees $(1,0)$, $(0,1)$; 
we denote by $(\on{deg}_{\mu},\on{deg}_{\delta})$ this grading; 
$\on{LBA}$ is then $\NN$-graded by the total degree 
$\on{deg}_{\mu}+\on{deg}_{\delta}$.  
If $x\in \on{LBA}(X,Y)$ and $X,Y,x$ are homogeneous, then 
$\on{deg}_{\mu}(x)-\on{deg}_{\delta}(x) = |X|-|Y|$, so 
$\on{LBA}(X,Y) = \on{LBA}^{\geq ||X|-|Y||}(X,Y)$. So the total degree 
completion of $\on{LBA}$ is a topological prop, with associated 
graded $\on{LBA}$. We denote by ${\bf LBA}$ the corresponding ${\bf Sch}$-prop. 

We then define a {\it quantization functor of Lie bialgebras} as a QSB in ${\bf LBA}$
quantizing $(\mu,\delta)$. Two quantization functors are equivalent if they are equivalent
as QSB's. Quantization functors of Lie bialgebras were constructed in \cite{EK1,EK2}. 

If now $\cC$ is a topological prop and 
$(\mu_{\cC},\delta_{\cC})$ is a Lie bialgebra structure on ${\bf id}$ 
in $\cC$, where the structure morphisms have positive valuation, then 
a quantization functor gives rise to a QSB in $\hat\cC$, 
quantizing $(\mu_{\cC}$ mod $\cC^{\geq 2},
\delta_{\cC}$ mod $\cC^{\geq 2})$; the structure morphisms
of this QSB are the images of the morphisms of the QSB in ${\bf LBA}$
under the prop morphism ${\bf LBA}\to \hat\cC$ taking $\mu,\delta$
to $\mu_{\cC},\delta_{\cC}$. 

\begin{remark} (QUE-QFSH equivalence)
If $V$ is a vector space and $X,Y\in\on{Ob(Sch)}$, we set 
$\cC_{0}^{V}(X,Y):= \on{Vect}(X(V),Y(V))$; this defines a prop. 
By extension of scalars, it gives rise to props $\cC_{0}^{V}[[\hbar]]$, 
$\cC_{0}^{V}((\hbar))$, with  
$\cC_{0}^{V}((\hbar))(X,Y):= \cC_{0}^{V}(X,Y)((\hbar))$, etc. 
We define two gradings on $\cC_{0}^{V}((\hbar))$: for $X,Y$ homogeneous, 
$\on{deg}_{QUE}(\cC^{V}_{0}(X,Y)\hbar^{i}) = |X|-|Y|+i$ and 
$\on{deg}_{QFSH}(\cC^{V}_{0}(X,Y)\hbar^{i}) = |Y|-|X|+i$. 
We then define $\cC^{V}_{QUE}$, $\cC^{V}_{QFSH}$ by 
$\cC^{V}_{QUE}(X,Y):=$ the part of $\cC^{V}_{0}(X,Y)((\hbar))$ with QUE-degree
$\geq ||X|-|Y||$ and $\cC^{V}_{QFSH}(X,Y):=$ the part of $\cC_{0}^{V}(X,Y)((\hbar))$ 
of QFSH-degree $\geq ||X|-|Y||$, for $X,Y$ homogeneous; these are topological props. 
Explicitly, 
$\cC_{QUE}^{V}(X,Y) = \hbar^{max(|Y|-|X|,0)}\cC_{0}^{V}(X,Y)[[\hbar]]$ and 
$\cC_{QFSH}^{V}(X,Y) = \hbar^{max(|X|-|Y|,0)}\cC_{0}^{V}(X,Y)[[\hbar]]$. 
Then a QSB in $\cC_{QUE}^{V}$ (resp., $\cC_{QFSH}^{V}$) gives rise to a QUE
deforming $U(V)$ (resp., a QFSH deforming $\hat S(V)$), where in both cases 
$V$ is equipped with a Lie bialgebra structure. The prop automorphism of 
$\cC_{0}^{V}((\hbar))$ given by $\cC_{0}^{V}(X,Y)((\hbar))\to 
\cC_{0}^{V}(X,Y)((\hbar))$, $x\mapsto \hbar^{|X|-|Y|}x$
for $X,Y$ homogeneous, restricts to a prop isomorphism  
$\cC_{QUE}^{V}\simeq \cC^{V}_{QFSH}$, which is   
compatible with the correspondence between QUE and QFSH algebras 
(\cite{Dr:QG,Gav}); namely, we have a commuting diagram 
$$\begin{matrix}
\{\on{QSB's\ in\ }\cC^{V}_{QUE}\} & \to & \{\on{QUE\ algebras}\} \\
{\scriptstyle\simeq} \downarrow & & \downarrow{\scriptstyle\simeq}\\
\{\on{QSB's\ in\ }\cC^{V}_{QFSH} \}  &\to & \{\on{QFSH\ algebras}\}
\end{matrix}$$
\end{remark}

\begin{remark} (Biquantization)
 Let $\cA_{0} := \kk[[uv]]\otimes_{\kk[uv]}\kk[u,v]\subset \kk[[u,v]]$; this is  
the subring of $\kk[[u,v]]$ of series $a(u,v) = \sum_{i,j\geq 0}
a_{ij}u^{i}v^{j}$ such that for some $N_{a}$, the support of $(a_{ij})$ is contained in 
$\{(i,j)||i-j|\leq N_{a}\}$. We have $\cA_{0}\simeq \kk[[uv]][u]\oplus 
v\kk[[uv]][v]$ as vector spaces. We have ring morphisms $\cA_{0}\to \kk[[u]]$
and $\cA_{0}\to \kk[[v]]$ obtained by setting $v=1$, resp., $u=1$. 

$\cC_{0}^{V}$ gives rise to the prop $\cC_{0}^{V}[[u,v]]$, where 
$\tilde{\cC}_{0}^{V}[[u,v]](X,Y):= 
\cC_{0}^{V}(X,Y)[[u,v]]$ and to the subprop $\cC^{V}_{\cA_{0}}(X,Y):=
\cC_{0}^{V}(X,Y)[[uv]]\otimes_{\kk[uv]}\kk[u,v] =  
\{\sum_{i,j}a_{ij}u^{i}v^{j}|$ for some $N_{a}, \on{supp}(a_{ij})\subset 
\{(i,j)||i-j|\leq N_{a}\}\}$ (this is $\cC^{V}_{0}(X,Y)\otimes \cA_{0}$ if
$V$ is finite dimensional). $\cC^{V}_{\cA_{0}}$ is graded by: 
$\on{deg}_{KT} := \on{deg}_{v} - \on{deg}_{u}+|X|-|Y|$. 
We denote by $\cC^{V}_{KT}$ the degree zero part of $\cC^{V}_{\cA_{0}}$
for this degree; this is a prop. Explicitly, we have 
$\cC^{V}_{KT}(X,Y) = u^{|X|-|Y|}\cC_{0}^{V}(X,Y)[[uv]]$ if $|X|\geq |Y|$
and $\cC^{V}_{KT}(X,Y) = v^{|Y|-|X|}\cC_{0}^{V}(X,Y)[[uv]]$ if $|X|<|Y|$. 
$\tilde\cC^{V}_{0}[[u,v]]$ is also graded by the total degree $\on{deg}_{u}
+\on{deg}_{v}$; this induces a filtration on $\cC^{V}_{KT}$, which is 
topological as $\cC^{V}_{KT}(X,Y) = (\cC^{V}_{KT})^{\geq ||X|-|Y||}(X,Y)$. 

A QSB in the prop $\cC^{V}_{KT}$ then gives rise to a biquantization of 
a Lie bialgebra structure on $V$, in the sense of \cite{KT}. 

Setting $v=1$ (resp., $u=1$), we get a prop morphisms $\cC^{V}_{\cA_{0}}\to 
\cC^{V}_{0}[[u]]$ (resp., $\cC^{V}_{\cA_{0}}\to \cC^{V}_{0}[[v]]$), which 
restrict to prop morphisms $\cC^{V}_{KT}\to \cC^{V}_{QFSH}$ 
(resp., $\cC^{V}_{KT}\to \cC^{V}_{QUE}$), 
where in the target props $\hbar$ is replaced by $u$ (resp., by $v$). The diagram of props 
$$
\begin{matrix} &\nearrow & \cC^{V}_{QUE}\\
\cC^{V}_{KT} & & \downarrow\scriptstyle{\simeq}\\
 & \searrow & \cC^{V}_{QFSH} \end{matrix}$$
commutes. It follows that the corresponding diagram between sets of QSB's in 
all three props commutes, so a biquantization of $V$ arising from a QSB 
in $\cC^{V}_{KT}$
gives rise to QUE and QFSH algebras, which correspond to each other 
under the category equivalence between QUE and QFSH algebras.
 
If now $(V,\mu_{V},\delta_{V})$ is a Lie bialgebra, 
$(u\mu_{V},v\delta_{V})$ defines a Lie bialgebra structure on 
${\bf id}$ in $\hat{\on{gr}}(\cC^{V}_{KT}) = \cC^{V}_{KT}$, where the
morphisms have positive valuation (in fact, total degree $1$). A quantization 
functor of Lie bialgebras then gives rise to a QSB in $\cC^{V}_{KT}$, i.e., 
a biquantization of $(V,\mu_{V},\delta_{V})$ in the sense of \cite{KT}. 
We recover in this way  Theorem 2.3 of \cite{KT}
(where a biquantization was constructed using the Etingof-Kazhdan construction). 
\end{remark}

\subsection{Quantization functors of quasi-Lie bialgebras}

Let again $\cC$ be a topological prop. A {\it quantized symmetric quasi-bialgebra}
(QSQB) in $\hat\cC$ is a quasi-bialgebra structure on $S$ in $\hat\cC$, whose reduction 
mod $\hat\cC^{\geq 1}$ coincides with the standard bialgebra structure on $S$, 
induced by the morphism ${\bf Sch}\to S$. Explicitly, this is the data of
morphisms $(m_{\cC},\Delta_{\cC},\Phi_{\cC},\eps_{\cC},\eta_{\cC})$
(we often drop $\eps_{\cC},\eta_{\cC}$)
in the same spaces as above, with $\Phi_{\cC}\in \hat\cC({\bf 1},S^{\otimes 3})$, 
satisfying the quasi-bialgebra axioms 
$$
m_{\cC}\circ (m_{\cC}\otimes \on{id}_{S}) 
= m_{\cC}\circ (\on{id}_{S}\otimes m_{\cC}), \quad 
\Delta_{\cC}\circ m_{\cC} = m_{\cC}^{\otimes 2} \circ \beta_{2} \circ 
\Delta_{\cC}^{\otimes 2}, $$
$$
[(\on{id}_{\cC}\otimes \Delta_{\cC}) \circ \Delta_{\cC}] * \Phi_{\cC}
= \Phi_{\cC} * [(\Delta_{\cC}\otimes \on{id}_{\cC}) \circ \Delta_{\cC}], 
$$
$$ (\eta_{\cC}\otimes \Phi_{\cC}) * [(\on{id}_{S}\otimes\Delta_{\cC}
\otimes \on{id}_{S})\circ\Phi_{\cC}] * (\Phi_{\cC}\otimes \eta_{\cC}) = 
[(\on{id}_{S}^{\otimes 2}\otimes\Delta_{\cC})\circ\Phi_{\cC}]*[(\Delta_{\cC}
\otimes\on{id}_{S}^{\otimes 2})\circ\Phi_{\cC}],$$
$$
(\on{id}_{S}\otimes\eps_{\cC}\otimes\on{id}_{S})\circ\Phi_{\cC}
=\eta_{\cC}^{\otimes 2}, \quad 
(\eps_{\cC}\otimes\on{id}_{S})\circ \Delta_{\cC}
= (\on{id}_{S} \otimes \eps_{\cC})\circ \Delta_{\cC}=
m_{\cC} \circ (\eta_{\cC}\otimes\on{id}_{S}) = 
m_{\cC}\circ (\on{id}_{S}\otimes \eta_{\cS})=\on{id}_{S}, 
$$
(we use the associativity of $m_{\cC}$ to define an associative operation  
$* : \hat\cC(X,S^{\otimes k}) \otimes \hat\cC(Y,S^{\otimes k})
\to \hat\cC(X\otimes Y,S^{\otimes k})$ by 
$x*y = m_{\cC}^{\otimes k}\circ \beta_{k}\circ (x\otimes y)$, 
where $\beta_{k}$ is the categorical version of $x_{1}\otimes x_{2}\otimes ...
\otimes y_{k}\mapsto x_{1}\otimes y_{1}\otimes...\otimes x_{k}\otimes y_{k}$; 
the isomorphisms ${\bf 1}\otimes X \simeq X \simeq X \otimes {\bf 1}$ for 
$X={\bf 1},{\bf id}$ are implied); and the same reduction conditions as 
above, together with $(\Phi_{\cC}$ mod $\hat\cC^{\geq 1}) = inj_{0}^{\otimes 3}$, 
where $inj_{0} : {\bf 1}\to S$ is the canonical morphism (note that this condition 
implies that $\Phi_{\cC}$ is invertible for the product $*$ on 
$\cC({\bf 1},S^{\otimes 3})$). 

Let $\hat\cC({\bf 1},S^{\otimes 2})_{1}\subset \hat\cC({\bf 1},S^{\otimes 2})$ 
be the set of $F_{\cC}$ such that $(F_{\cC}$ mod $\hat\cC^{\geq 1}) = 
inj_{0}^{\otimes 2}$ and $(\eps_{\cC}\otimes \on{id}_{S})\circ F_{\cC}=
(\on{id}_{S}\otimes \eps_{\cC})\circ F_{\cC}=\eta_{\cC}$. 
Let $(m_{\cC},...)$ be a QSQB in $\hat\cC$; when equipped with the product $*$
induced by $m_{\cC}$, $\hat\cC({\bf 1},S^{\otimes 2})_{1}$ is a group. 
This group acts on the set of all QSQB's in $\hat\cC$ with fixed $m_{\cC}$ 
by $F_{\cC}\star(m_{\cC},\Delta_{\cC},
\eps_{\cC},\eta_{\cC},\Phi_{\cC}) = (m_{\cC},F_{\cC}*\Delta_{\cC}
*F_{\cC}^{-1},\eps_{\cC},\eta_{\cC}, [(\eta_{\cC}\otimes F_{\cC}) * 
((\on{id}_{S}\otimes\Delta_{\cC}) \circ F_{\cC})]*\Phi_{\cC}*
[(F_{\cC}\otimes \eta_{\cC}) * 
((\Delta_{\cC}\otimes \on{id}_{S}) \circ F_{\cC})]^{-1})$. Two 
QSQB's related in this way are called twist equivalent. 

The group $\hat\cC(S,S)_{1}$ acts on the set of QSQB's in $\hat\cC$ by 
$i_{\cC}*(m_{\cC},...,\Phi_{\cC}):= (i_{\cC}\circ m_{\cC}\circ 
(i_{\cC}^{-1})^{\otimes 2},
...,i_{\cC}^{\otimes 3}\circ\Phi_{\cC})$. Two QSQB's related by this group action 
are called equivalent. 

As in \cite{Dr:QH}, one proves: 

\begin{prop}
1) For any QSQB $(m_{\cC},\Delta_{\cC},\Phi_{\cC})$ in $\hat\cC$, there exists 
$(\mu_{\cC},\delta_{\cC},\varphi_{\cC})$, where $\mu_{\cC}\in \on{gr}^{1}(
\wedge^{2},{\bf id})$, $\delta_{\cC}\in \on{gr}^{1}(\cC)({\bf id},\wedge^{2})$, 
$\varphi_{\cC}\in\on{gr}^{1}(\cC)({\bf 1},\wedge^{3})$, such that: 
$(m_{\cC}\circ \on{Alt}_{2}\circ inj_{1}^{\otimes 2}$ mod $\hat\cC^{\geq 2})
= inj_{1}\circ \mu_{\cC}$, $(\on{Alt}_{2}\circ \Delta_{\cC}\circ inj_{1}$ mod
$\hat\cC^{\geq 2}) = inj_{1}^{\otimes 2}\circ\delta_{\cC}$, $(\on{Alt}_{3}
\circ \Phi_{\cC}$ mod $\hat\cC^{\geq 2}) = inj_{1}^{\otimes 3}\circ 
\varphi_{\cC}$. The triple $(\mu_{\cC},\delta_{\cC},\varphi_{\cC})$ 
equips the object ${\bf id}$ with a quasi-Lie bialgebra structure in 
$\on{gr}(\cC)$, where the morphisms have degree $1$. We call it the classical limit of 
$(m_{\cC},\Delta_{\cC},\Phi_{\cC})$. 

2) Equivalent and twist equivalent QSQB's have the same classical limit. 
\end{prop}

{\em Proof.} Let $(m_{\cC},\Delta_{\cC},\Phi_{\cC})$ in 
$\hat\cC$. Let $\tilde\mu_{\cC}:= (m_{\cC}\circ \on{Alt}_{2}\circ 
inj_{1}^{\otimes 2}$ mod $\hat\cC^{\geq 2})$; we have $\tilde\mu_{\cC}
\in \on{gr}^{1}(\cC)(\wedge^{2},S)$. Composing $\Delta_{\cC} \circ 
(m_{\cC}\circ \on{Alt}_{2}) = [(m_{\cC}\circ \on{Alt}_{2})\otimes m_{\cC}
+ (m_{\cC}\circ \beta) \otimes (m_{\cC}\circ \on{Alt}_{2})] \circ \beta_{23}
\circ (\Delta_{\cC}\otimes \Delta_{\cC})$ with $inj_{1}^{\otimes 2}$, and 
$\Delta_{\cC}\circ inj_{1} = inj_{1}\otimes inj_{0} + inj_{0}\otimes inj_{1}$
mod $\hat\cC^{\geq 2}$, we get $\Delta_{S} \circ \tilde\mu_{\cC} = 
\tilde\mu_{\cC}\otimes inj_{0} + inj_{0}\otimes\tilde\mu_{\cC}$, which 
implies the existence of $\mu_{\cC}$ (here $\Delta_{S}:S\to S^{\otimes 2}$
is the image in $\on{gr}^{0}(\cC)$ of the coproduct morphism of $S\in 
\on{Ob}({\bf Sch})$). 

We next prove that if $x\in \hat\cC^{\geq k}(X,S^{\otimes n})$ and $\Psi\in \cC({\bf 1},
S^{\otimes 2})$ is such that $(\Psi$ mod $\hat\cC^{\geq 1}) = inj_{0}^{\otimes n}$,
then $\Psi*x*\Psi^{-1} = x\quad  \on{mod}\quad  \hat\cC^{\geq k+2}$. 
Indeed, $\Psi*x*\Psi^{-1} = [(\Psi-\eta_{\cC}^{\otimes n})*x - x*(\Psi-
\eta_{\cC}^{\otimes n})] * \Psi^{-1}$; the result then follows from 
$\Psi-\eta_{\cC}^{\otimes n}\in\hat\cC^{\geq 1}(X,S^{\otimes n})$
and $m_{\cC} - m_{\cC}\circ\beta\in \hat\cC^{\geq 1}(S^{\otimes 2},S)$. 

It follows that $(\on{id}_{S}\otimes\Delta_{\cC})\circ\Delta_{\cC} = 
(\Delta_{\cC}\otimes \on{id}_{S})\circ\Delta_{\cC}$ mod $\hat\cC^{\geq 2}$. 
Let $\beta_{321};S^{\otimes 3}\to S^{\otimes 3}$ be the analogue of 
$x_{1}\otimes x_{2}\otimes x_{3}\to x_{3}\otimes x_{2}\otimes x_{1}$. 
As $(\Delta_{\cC}\otimes \on{id}_{S})\circ\Delta_{\cC} - \beta_{321}\circ
(\on{is}_{S}\otimes\Delta_{\cC})\circ\Delta_{\cC}= 
[(\Delta_{\cC} - \beta\circ\Delta_{\cC})\otimes\on{id}_{S}]\circ\Delta_{\cC}
+[(\beta\circ\Delta_{\cC})\otimes \on{id}_{S}]\circ [\Delta_{\cC} - 
\beta\circ \Delta_{\cC}]$, we get 
$[(\Delta_{\cC}-\beta\circ\Delta_{\cC})\otimes\on{id}_{S} 
- \on{id}_{S}\otimes (\Delta_{\cC}-\beta\circ\Delta_{\cC})] \circ\Delta_{\cC}= 
[\on{id}_{S}\otimes (\beta\circ\Delta_{\cC}) - (\beta\circ\Delta_{\cC})
\otimes\on{id}_{S}]\circ [\Delta_{\cC} - \beta\circ\Delta_{\cC}]$ mod
$\hat\cC^{\geq 2}$. Let $\tilde\delta_{\cC}:= (\on{Alt}_{2}\circ\Delta_{\cC}
\circ inj_{1}$ mod $\hat\cC^{\geq 2})$; we have $\tilde\delta_{\cC}\in 
\on{gr}^{1}({\bf id},S^{\otimes 2})$. Composing this with $inj_{1}$, we get 
$\tilde\delta_{\cC}\otimes inj_{0} - inj_{0}\otimes\tilde\delta_{\cC}
+ (\Delta_{0}\otimes \on{id}_S -\on{id}_S
\otimes\Delta_{0})\circ\tilde\delta_{\cC}=0$.
The degree $n$ cohomology group 
of the co-Hochschild (or Cartier)
complex $...\to\on{gr}^{1}(\cC)(X,S^{\otimes n})\to \on{gr}^{1}(\cC)(X,
S^{\otimes n+1})\to ...$ is $\on{gr}^{1}(\cC)(X,\wedge^{n})$, and the map 
taking a coboundary to its cohomology class is $(1/n!)\on{Alt}_{n}$. As
$\tilde\delta_{\cC}$ is a coboundary and is antisymmetric, it coincides with its
cohomology class, so there exists the announced $\delta_{\cC}$. 

The pentagon equation implies that $(\Phi_{\cC}-\eta_{\cC}^{\otimes 3}$ 
mod $\hat\cC^{\geq 2})\in \on{gr}^{1}({\bf 1},S^{\otimes 3})$ is a co-Hochschild
coboundary. It follows that its image by $\on{Alt}_{3}$ coincides with its cohomology 
class, which implies the existence of $\varphi_{\cC}$. 

Before showing that the classical limit $(\mu_{\cC},\delta_{\cC},\varphi_{\cC})$
of $(m_{\cC},\Delta_{\cC},\Phi_{\cC})$ satisfies the quasi-Lie bialgebra relations, 
let us show that it is invariant under twist equivalence. Let $F_{\cC}\in 
\hat\cC(S,S)_{1}$, let $(\tilde m_{\cC},\tilde\Delta_{\cC},\tilde\Phi_{\cC}):= 
F_{\cC}\star (m_{\cC},\Delta_{\cC},\Phi_{\cC})$ and $(\tilde\mu_{\cC},
\tilde\delta_{\cC},\tilde\varphi_{\cC})$ be its classical limit. As $\tilde m_{\cC}
= m_{\cC}$, $\tilde\mu_{\cC}=\mu_{\cC}$.  
We have  
$\tilde\Delta_{\cC} - \beta \circ \tilde\Delta_{\cC} = 
F_{\cC} * \Delta_{\cC}*F_{\cC}^{-1} - (\beta\circ F_{\cC}) * 
(\beta\circ \Delta_{\cC}) * (\beta\circ F_{\cC}^{-1}) = 
F_{\cC}*(\Delta_{\cC}-\beta\circ\Delta_{\cC})*F_{\cC}^{-1}
+ F_{\cC} * [(\beta\circ\Delta_{\cC}) * (F_{\cC}^{-1}*(\beta\circ F_{\cC}))
- (F_{\cC}^{-1}*(\beta\circ F_{\cC})) * (\beta\circ\Delta_{\cC}) ] * 
(\beta\circ F_{\cC}^{-1})$; it follows that $\tilde\Delta_{\cC} - 
\beta\circ\tilde\Delta_{\cC} = \Delta_{\cC} - \beta\circ\Delta_{\cC}$
mod $\hat\cC^{\geq 3}$, so $\tilde\delta_{\cC}=\delta_{\cC}$. 
Finally, $(\tilde\Phi_{\cC} - \eta_{\cC}^{\otimes 3}$ mod $\hat\cC^{\geq 2})
= (\Phi_{\cC} - \eta_{\cC}^{\otimes 3}$ mod $\hat\cC^{\geq 2}) + d
(F_{\cC}-\eta_{\cC}^{\otimes 2}$ mod $\hat\cC^{\geq 2})$, where $d$ is the 
co-Hochschild differential, which implies that $\tilde\varphi_{\cC}=\varphi_{\cC}$.  

We have $ (m_{\cC} \circ \on{Alt}_{2}) \circ 
[(m_{\cC}\circ \on{Alt}_{2})\otimes \on{id}_{S}] \circ \on{Alt}_{3}= 0$ 
(equality in $\hat\cC^{\geq 2}$); composing with $inj_{1}$ and taking the class modulo 
$\hat\cC^{\geq 3}$, we get $\mu_{\cC}\circ (\mu_{\cC}\otimes \on{id}_{{\bf id}})
\circ \on{Alt}_{3} = 0$. 

We have $[\on{Alt}_{2} \circ \Delta_{\cC}] \circ [m_{\cC} \circ 
\on{Alt}_{2}] = [m^{(2)}_{\cC} - m^{(2)}_{\cC} \circ \beta_{S^{\otimes 2}}] \circ 
[(\on{Alt}_{2}\circ\Delta_{\cC})\otimes\Delta_{\cC} + (\beta\circ\Delta_{\cC})\otimes
(\on{Alt}_{2}\circ\Delta_{\cC})]$, where $m^{(2)}_{\cC}:S^{\otimes 2}\otimes S^{\otimes 2}
\to S^{\otimes 2}$ is $(m_{\cC}\otimes m_{\cC})\circ \beta_{2}$ and 
$\beta_{S^{\otimes 2}} : S^{\otimes 2}\otimes S^{\otimes 2} \to 
S^{\otimes 2}\otimes S^{\otimes 2}$ is the braiding; this equality takes place in 
$\hat\cC^{\geq 2}$. Composing with $inj_{1}^{\otimes 2}$ and taking the class modulo 
$\hat\cC^{\geq 3}$, we get  $\delta_{\cC}\circ\mu_{\cC} = \on{Alt}_{2} \circ 
(\mu_{\cC}\otimes\on{id}_{{\bf id}}) \circ (\on{id}_{{\bf id}}\otimes
\delta_{\cC}) \circ \on{Alt}_{2}$. 

We have $(\on{id}_{S}\otimes\Delta_{\cC} - \Delta_{\cC}\otimes \on{id}_{S})
\circ \Delta_{\cC} = [m_{\cC}^{(3)} 
- m_{\cC}^{(3)} \circ \beta_{S^{\otimes 3}}] \circ 
\{(\Phi_{\cC}-\eta_{\cC}^{\otimes 3}) \otimes  [(\Delta_{\cC}
\otimes \on{id}_{S}) \circ \Delta_{\cC}] \}$; 
also $\on{Alt}_{3}\circ (\on{id}_{S}\otimes\Delta_{\cC}-\Delta_{\cC}
\otimes\on{id}_{S})\circ\Delta_{\cC} = -{1\over 2}\on{Alt}_{3} \circ 
[(\on{Alt}_{2}\circ \Delta_{\cC})\otimes \on{id}_{S}]\circ 
[\on{Alt}_{2}\circ\Delta_{\cC}]$; so  (equality in $\hat\cC^{\geq 2}$; 
as $m_{\cC}^{(3)}$, $\beta_{S^{\otimes 3}}$ commute with the action of $\SG_{3}$, 
we get $-{1\over 2}\on{Alt}_{3} \circ [(\on{Alt}_{2}\circ\Delta_{\cC})
\otimes \on{id}_{S}] \circ [\on{Alt}_{2}\circ\Delta_{\cC}] = 
[m_{\cC}^{(3)} - m_{\cC}^{(3)}\circ \beta_{S^{\otimes 3}}] \circ 
\{[\on{Alt}_{3}\circ (\Phi_{\cC}-\eta_{\cC}^{\otimes 3})]\otimes 
[(\Delta_{\cC}\otimes \on{id}_{S})\circ\Delta_{\cC}]\}$; this equality holds
in $\hat\cC^{\geq 2}$; taking its class modulo 
$\hat\cC^{\geq 3}$ and composing with $inj_{1}$, we get $\on{Alt}_{3}\circ
(\delta_{\cC}\otimes \on{id}_{{\bf id}}) \circ \delta_{\cC} = 
\on{Alt}_{3} \circ (\mu\otimes\on{id}_{{\bf id}^{\otimes 3}}) \circ 
(\on{id}_{{\bf id}} \otimes \varphi_{\cC})$. 

Using the co-Hochschild complex, one may find a twist $F_{\cC}$ with 
$F_{\cC}\star (m_{\cC},\Delta_{\cC},\Phi_{\cC}) = (m_{\cC},\tilde\Delta_{\cC},
\tilde\Phi_{\cC})$, where $(\tilde\Phi_{\cC}-\eta_{\cC}^{\otimes 3}$ mod
$\hat\cC^{\geq 2}) = inj_{1}^{\otimes 3}\circ \varphi_{\cC}$. 
Let $\tilde\varphi_{\cC}:= \Phi_{\cC}-\eta_{\cC}^{\otimes 2}$. Set 
$\tilde\varphi_{\cC}^{2,3,4}:= \eta_{\cC}\otimes\tilde\varphi_{\cC}$, 
$\tilde\varphi_{\cC}^{1,2,34}:= (\on{id}_{S}^{\otimes 2} \otimes 
\tilde\Delta_{\cC})
\circ \tilde\varphi_{\cC}$, etc. Then the pentagon equation yields
$\tilde\varphi_{\cC}^{1,2,34}+\tilde\varphi_{\cC}^{12,3,4} 
- \tilde\varphi_{\cC}^{2,3,4} - \tilde\varphi_{\cC}^{1,23,4}
-\tilde\varphi_{\cC}^{1,2,3} 
=  \tilde\varphi_{\cC}^{2,3,4} * \tilde\varphi_{\cC}^{1,23,4}
+\tilde\varphi_{\cC}^{2,3,4} * \tilde\varphi_{\cC}^{1,2,3} + \tilde\varphi_{\cC}^{1,23,4}
* \tilde\varphi_{\cC}^{1,2,3} - \tilde\varphi_{\cC}^{1,2,34} * \tilde\varphi_{\cC}^{12,3,4} 
+ \tilde\varphi_{\cC}^{1,2,3} * \tilde\varphi_{\cC}^{1,23,4} * \tilde\varphi_{\cC}^{1,2,3}$. 
As $\varphi_{\cC}\in\hat\cC^{\geq 1}({\bf 1},\wedge^{3})$, we have 
$\tilde\varphi_{\cC}^{12,3,4} = \tilde\varphi_{\cC}^{1,3,4} + \tilde
\varphi_{\cC}^{2,3,4}$ mod $\hat\cC^{\geq 2}$. It follows that 
$\tilde\varphi_{\cC}^{1,2,34} - \tilde\varphi_{\cC}^{1,23,4}+\tilde\varphi_{\cC}^{12,3,4}
- \tilde\varphi_{\cC}^{1,2,3}-\tilde\varphi_{\cC}^{2,3,4} = \tilde\varphi_{\cC}^{2,3,4}
*\tilde\varphi_{\cC}^{1,3,4} -  \tilde\varphi_{\cC}^{1,2,4}
*\tilde\varphi_{\cC}^{1,3,4} +  \tilde\varphi_{\cC}^{1,2,4}
*\tilde\varphi_{\cC}^{1,2,3}$ mod $\hat\cC^{\geq 3}$. 
The r.h.s. of this equality belongs to $\hat\cC^{\geq 2}$, and because 
of the form of $(\tilde\varphi_{\cC}$ mod $\hat\cC^{\geq 2})$, so does its l.h.s. 
Compose both sides of this equality with $\on{Alt}_{4}$. 
As $\tilde\varphi_{\cC}^{1,2,3} = -\tilde\varphi_{\cC}^{1,3,2}$ mod $\hat\cC^{\geq 2}$, 
$\tilde\varphi_{\cC}^{2,3,4}*\tilde\varphi_{\cC}^{1,3,4}$ is symmetric in indices 3,4
modulo $\hat\cC^{\geq 3}$; the image in in $\on{gr}^{2}(\cC)$ of the composition 
of this term with $\on{Alt}_{4}$ is then zero; in the same way, 
so it the image of the composition of the r.h.s. with $\on{Alt}_{4}$. It
follows that the image in $\on{gr}^{2}(\cC)$ of $\on{Alt}_{4} \circ 
(\tilde\varphi_{\cC}^{12,3,4} - \tilde\varphi_{\cC}^{1,23,4}+\tilde\varphi_{\cC}^{1,2,34}
-\tilde\varphi_{\cC}^{2,3,4}-\tilde\varphi_{\cC}^{1,2,3})=0$. 
As $\tilde\delta_{\cC}=\delta_{\cC}$,  
$(\on{Alt}_{2}\otimes \on{id}_{S}^{\otimes 2}) \circ 
\tilde\varphi_{\cC}^{12,3,4} = (\delta_{\cC}\circ 
\on{id}_{{\bf id}}^{\otimes 2}) \circ \varphi_{\cC}$ in 
$\on{gr}^{2}(\cC)$. We then get $\on{Alt}_{4}\circ (\delta_{\cC}
\otimes\on{id}_{{\bf id}}^{\otimes 2})\circ \varphi_{\cC}=0$. 

If now $i_{\cC}\in \hat\cC(S,S)_{1}$, we have $i_{\cC}\circ inj_{1} = inj_{1}$
mod $\hat\cC^{\geq 2}$, which implies that the classical limit is invariant under 
equivalence. \hfill \qed\medskip

Let $\on{QLBA}$ be the prop with generators $\mu\in\on{QLBA}(\wedge^{2},{\bf id})$,
$\delta\in\on{QLBA}({\bf id},\wedge^{2})$, 
$\varphi\in \on{QLBA}({\bf 1},\wedge^{3})$ and relations 
$$
\mu\circ (\mu\otimes \on{id}_{{\bf id}}) \circ \on{Alt}_{3} = 0, \quad
\delta \circ \mu = \on{Alt}_{2} \circ (\mu\otimes\on{id}_{{\bf id}})
\circ (\on{id}_{{\bf id}}\otimes \delta) \circ \on{Alt}_{2}, $$
$$
\on{Alt}_{3} \circ (\delta\otimes\on{id}_{{\bf id}}) \circ\delta = 
\on{Alt}_{3}\circ (\mu\otimes\on{id}_{{\bf id}}^{\otimes 3})
\circ (\on{id}_{{\bf id}}\otimes \varphi), \quad
\on{Alt}_{4}\circ (\delta\otimes\on{id}_{{\bf id}}^{\otimes 2})
\circ\varphi=0. 
$$
This prop is graded by $\{(u,v)\in\ZZ^{2}|v\geq 0,2u+v\geq 0\}$, 
with $\mu,\delta,\varphi$ of degrees
$(1,0)$, $(0,1)$, $(-1,2)$; we denote this degree $(\on{deg}_{\mu},\on{deg}_{\delta})$.
QLBA is then $\NN$-graded by the total degree $\on{deg}_{\mu}+\on{deg}_{\delta}$; 
the generators have then degree $1$. 
If $x\in \on{QLBA}(X,Y)$ and $x,X,Y$ are homogeneous, then $\on{deg}_{\mu}(x)
-\on{deg}_{\delta}(x)=|X|-|Y|$, so $\on{QLBA}(X,Y) = 
\on{QLBA}^{\geq {1\over 3}||X|-|Y||}(X,Y)$, which implies that the total degree
completion of QLBA is a topological prop, with associated graded QLBA. 
We denote by ${\bf QLBA}$ the corresponding ${\bf Sch}$-prop. In the prop 
$\on{QLBA}$, the object ${\bf id}$ is equipped with a quasi-Lie bialgebra structure, 
and $\on{QLBA}$ is an initial object in the category of props equipped with a quasi-Lie
bialgebra structure on ${\bf id}$. We have an obvious prop morphism $\on{QLBA}
\to \on{LBA}$, given by $\mu,\delta,\varphi\mapsto \mu,\delta,0$. 

A {\it quantization functor of quasi-Lie bialgebras} is then a QSQB in 
${\bf QLBA}$ admitting $(\mu,\delta,\varphi)$ as its classical limit. 
Two quantization functors are (twist) equivalent if they are as QSB's. 

We have  a map $\{$quantization functors of quasi-Lie bialgebras$\}\to\{$quantization 
functors of Lie bialgebras$\}$, induced by the above morphism $\on{QLBA}\to
\on{LBA}$; this map takes equivalent functors to equivalent functors, and it
takes two twist equivalent functors to the same image. 

The main result of our paper is: 

\begin{thm} \label{thm:main}
The map $\{$quantization functors of quasi-Lie bialgebras$\}/$(equivalence, 
twist equivalence)$\to\{$quantization functors of Lie bialgebras$\}/$(equivalence)
is a bijection. 
\end{thm}

\begin{remark} If $V$ is a vector space, then a QSQB in $\cC^{V}_{QUE}$
gives rise to a QHQUE algebra deforming $U(V)$, in the sense of \cite{Dr:QH}. 
A quasi-Lie bialgebra structure $(\mu_{V},\delta_{V},\varphi_{V})$ on
$V$ then gives rise to a quasi-Lie bialgebra structure on ${\bf id}$ in
$\cC^{V}_{QUE}$, namely $(\mu_{V},\hbar\delta_{V},\hbar^{2}\varphi_{V})$. 

As above, if $\cC$ is a topological prop and 
$(\mu_{\cC},\delta_{\cC},\varphi_{\cC})$ is a quasi-Lie bialgebra
structure on the object ${\bf id}$ in $\cC$, with morphisms of valuation 
$\geq 1$, then a quantization functor for quasi-Lie bialgebras yield a 
QSQB in $\hat\cC$ quantizing $(\mu_{\cC},\delta_{\cC},\varphi_{\cC})$
mod $\hat\cC^{\geq 2}$. 

In particular, for each quasi-Lie bialgebra structure $(\mu_{V},\delta_{V},\varphi_{V})$
on a vector space $V$, a quantization functor for quasi-Lie bialgebras gives rise to a 
QHQUE algebra quantizing $V$, in the sense of \cite{Dr:QH}. 
\end{remark}

 \subsection{Deformation complexes}

Let $\cD$ be a prop. Set $C_{\cD}^{p,q}:= \cD(\wedge^{p+1},\wedge^{q+1})$
for $p,q\in\ZZ$ and $C_{\cD}:= \oplus_{p,q\in\ZZ}C_{\cD}^{p,q}$; this is a 
$\ZZ^{2}$-graded vector space (hence $\ZZ$-graded by the total degree). We define
a $\ZZ^{2}$-graded Lie superalgebra structure on $C_{\cD}$ as follows. For 
$a\in C^{p,q}_{\cD}$, 
$a'\in C^{p',q'}_{\cD}$, we set $a'\diamond a:= \on{Alt}_{q+q'+1} \circ (a'\otimes 
\on{id}_{{\bf id}^{\otimes q}})\circ (\on{id}_{{\bf id}^{\otimes p'}}\otimes a)
\circ \on{Alt}_{p+p'+1}$, and define their Schouten bracket 
$[a,a']:= a\diamond a' - (-1)^{(p+q)(p'+q')}a'\diamond a$. 
The condition that $(\mu_{\cD},\delta_{\cD})$ defines a Lie bialgebra structure on 
${\bf id}$ in $\cD$ is then $[\mu_{\cD}\oplus\delta_{\cD},\mu_{\cD}\oplus
\delta_{\cD}]=0$; the condition that $(\mu_{\cD},\delta_{\cD},\varphi_{\cD})$
defines a quasi-Lie bialgebra structure on ${\bf id}$ is 
$[\mu_{\cD}\oplus\delta_{\cD}\oplus\varphi_{\cD},
\mu_{\cD}\oplus\delta_{\cD}\oplus\varphi_{\cD}]=0$. 
The bracket with $\mu_{\cD}\oplus\delta_{\cD}$ (resp., $\mu_{\cD}\oplus
\delta_{\cD}\oplus\varphi_{\cD}$) then defines a 
complex structure on $C_{\cD}$ (graded by the total degree); we denote 
by $H^{i}_{\cD}$ the corresponding cohomology groups. 

If $\cC$ is a topological prop, we define a quotient prop $\cC_{\leq n}$ by 
$\cC_{\leq n}(X,Y)=\cC/\cC^{\geq n+1}(X,Y)$; this is also a topological prop; 
we have a projective system $...\to \cC_{\leq n}\to...\to \cC_{\leq 0}$ and 
$\cC(X,Y) = \on{lim}_{\leftarrow}\cC_{\leq n}(X,Y)$. 

If now $(m_{\cC}^{\leq n},\Delta_{\cC}^{\leq n})$ is a QSB in 
$\cC_{\leq n}$ with classical 
limit $(\mu_{\cC},\delta_{\cC})$, then the obstruction to extend it to a QSB in 
$\cC_{\leq n+1}$ belongs to $H^{2}_{\on{gr}(\cC)}[n+2]$, and the set of such 
extensions modulo equivalence is an affine space over $H^{1}_{\on{gr}(\cC)}[n+1]$
(here $[n]$ means the degree $n$ part for the grading of $\on{gr}(\cC)$); see
\cite{Dr:QG,GS,SS}. 

In the same way, if $(m_{\cC}^{\leq n},\Delta_{\cC}^{\leq n},
\Phi_{\cC}^{\leq n})$ is a QSQB in $\cC_{\leq n}$ with classical limit
$(\mu_{\cC},\delta_{\cC},\varphi_{\cC})$, then the obstruction to 
extend it to degree $n+1$ belongs to $H^{2}_{\on{gr}(\cC)}[n+2]$, and the set of
extensions, modulo equivalence and twist equivalence, is an affine space 
over $H^{1}_{\on{gr}(\cC)}[n+1]$ (see \cite{SS}). In both cases, cohomologies 
are relative to the differentials defined by the classical limits. 

Note that when $\cC = \on{LBA}$ or $\on{QLBA}$, $\on{gr}(\cC) = \cC$. 
The prop morphism $\on{QLBA}\to \on{LBA}$ induces morphisms 
of complexes $C_{\on{QLBA}}\to C_{\on{LBA}}$ and of the 
cohomologies $H^{i}_{\on{QLBA}}\to H^{i}_{\on{LBA}}$.

\subsection{Proof of the main theorem} \label{pf:main}

To show Theorem \ref{thm:main}, we will prove: 

\begin{thm} \label{thm:isom} The maps $H^{i}_{\on{QLBA}}\to 
H^{i}_{\on{LBA}}$ are isomorphisms for any $i\geq 0$. 
\end{thm}

We will prove this in Section \ref{sec:comparison}. 

Let us explain why this implies Theorem \ref{thm:main}. 
We will prove inductively over $n$ that the map 
$\on{red}_{n}:\{$QSQB's in ${\bf QLBA}_{\leq n}$ quantizing $(\mu,\delta,\varphi)\}
/($equivalence, twist equivalence$)\to \{$QSB's in ${\bf LBA}_{\leq n}$
quantizing $(\mu,\delta)\}$ is a bijection. We denote by $pr^{\cC}_{n+1,n}$
the reduction map $\{$QSB in $\hat\cC_{\leq n+1}\} \to \{$QSB 
in $\hat\cC_{\leq n}\}$. Then $\on{red}_{n}\circ pr_{n+1,n}^{\on{QLBA}} = 
pr_{n+1,n}^{\on{LBA}}\circ \on{red}_{n+1}$. 

Assume that $\on{red}_{n}$ is bijective and let us show that $\on{red}_{n+1}$
is bijective. We first show that it is injective. If $[(m,\Delta,\Phi)]$ and 
$[(m',\Delta',\Phi')]$ are two classes of QSQB's in ${\bf QLBA}_{\leq n+1}$
with the same image by $\on{red}_{n+1}$, then the injectivity of 
$\on{red}_{n}$ implies that their images by $pr_{n+1,n}$ coincide. So 
$[(m,\Delta,\Phi)]$ and $[(m',\Delta',\Phi')]$ differ by an element  
$\omega\in H^{1}_{\on{QLBA}}[n+1]$. Their images by $\on{red}_{n+1}$ are 
classes of QSB's in ${\bf LBA}_{\leq n+1}$, whose reductions in 
${\bf LBA}_{\leq n}$ are equivalent; these classes differ by the image of $\omega$
under $H^{1}_{\on{QLBA}}[n+1]\to H^{1}_{\on{LBA}}[n+1]$. As this map 
is injective, $\omega=0$ so $[(m,\Delta,\Phi)] = [(m',\Delta',\Phi')]$, which proves
the injectivity of $\on{red}_{n+1}$. 

Let us now show that $\on{red}_{n+1}$ is surjective. Let $[(m_{\leq n+1},
\Delta_{\leq n+1})]$ be a 
class of QSB in ${\bf LBA}_{\leq n+1}$. Set $[(m_{\leq n},\Delta_{\leq n})]
:= p_{n+1,n}^{\on{LBA}} ([(m_{\leq n+1},\Delta_{\leq n+1})])$; this 
is the class of a QSB in ${\bf LBA}_{\leq n}$. Let 
$[(\tilde m_{\leq n},\tilde \Delta_{\leq n},\Phi_{\leq n})]$ 
be the preimage of $[(m_{\leq n},\Delta_{\leq n})]$ by $\on{red}_{n}$. 
The obstruction to extending it to a QSQB in ${\bf QLBA}_{\leq n+1}$
is a cohomology class in $H^{2}_{\on{QLBA}}[n+2]$. The image of this 
class by $H^{2}_{\on{QLBA}}[n+2]\to H^{2}_{\on{LBA}}[n+2]$ is the 
obstruction to extending $[(m_{\leq n},\Delta_{\leq n})]$ to a QSB in 
${\bf LBA}_{\leq n+1}$. The existence of $[(m_{\leq n+1},\Delta_{\leq n+1})]$
implies that this class in $H^{2}_{\on{LBA}}[n+2]$ is zero. As the map 
$H^{2}_{\on{QLBA}}[n+2]\to H^{2}_{\on{LBA}}[n+2]$ is injective, the 
obstruction class in $H^{2}_{\on{QLBA}}[n+2]$ is zero, and 
$[(\tilde m_{\leq n},\tilde \Delta_{\leq n},\Phi_{\leq n})]$ may be extended to
a QSQB in ${\bf QLBA}_{\leq n+1}$. Let 
$[(\tilde m_{\leq n+1},\tilde \Delta_{\leq n+1},\Phi_{\leq n+1})]$ be 
such an extension. The difference between $\on{red}_{n+1}
([(\tilde m_{\leq n+1},\tilde \Delta_{\leq n+1},\Phi_{\leq n+1})])$
and $[(m_{\leq n+1},\Delta_{\leq n+1})]$ is a cohomology class in 
$H^{1}_{\on{LBA}}[n+1]$. As the map $H^{1}_{\on{QLBA}}[n+1]
\to H^{1}_{\on{LBA}}[n+1]$ is surjective, this is the image of a cohomology 
class in $H^{1}_{\on{QLBA}}[n+1]$. Substracting this cohomology class from 
$[(\tilde m_{\leq n+1},\tilde \Delta_{\leq n+1},\Phi_{\leq n+1})])$, we
obtain a preimage of $[(m_{\leq n+1},\Delta_{\leq n+1})]$ by $\on{red}_{n+1}$. 
 
\section{Structure of the prop $\on{QLBA}$}

In order to establish Theorem \ref{thm:isom}, we study the structure of $\on{QLBA}$. 

\subsection{Products of ideals in props}

If $\cC$ is a prop and $I_{1},...,I_{n}$ are ideals of $\cC$, 
then the product $I_{1}...I_{n}$ 
is the smallest ideal containing the morphisms $f_{1}*...*f_{n}$, where
$f_{i}$ is morphism in $I_{i}$ and $*$ is $\circ$ or $\otimes$. 
One defines in this way the powers $I^{n}$
of an ideal. 

\subsection{Structure of the prop $\on{LBA}$}

Define $\on{LA}$ (resp., $\on{LCA}$) as the prop generated by $\mu\in
\on{LA}({\bf id}^{\otimes 2},{\bf id})$ subject to the antisymmetry and 
Jacobi relation (resp., $\delta\in\on{LCA}({\bf id},{\bf id}^{\otimes 2})$
subject to antisymmetry and the co-Jacobi relation). We have prop morphisms 
$\on{LA},\on{LCA}\to \on{LBA}$. The structure of $\on{LBA}$ is given 
by 

\begin{prop} \label{prop:PBW}
For $(X_{i})_{i\in I}$, $(Y_{j})_{j\in J}$ finite families of objects in 
$\on{Ob(Sch)}$, we have an isomorphism
\begin{equation} \label{str:LBA}
\on{LBA}(\otimes_{i}X_{i},\otimes_{j}Y_{j}) \simeq
\oplus_{(Z_{ij})\in \on{Irr(Sch)}^{I\times J}} 
[\otimes_{i}\on{LCA}(X_{i},\otimes_{j}Z_{ij})] 
\otimes [\otimes_{j}\on{LA}(\otimes_{i}Z_{ij},Y_{j})], 
\end{equation}
whose inverse is the direct sum of the maps $(\otimes_{i}c_{i})\otimes 
(\otimes_{j} a_{j})\mapsto (\otimes_{j}a_{j}) \circ \beta_{I,J} \circ
(\otimes_{i}c_{i})$, where $\beta_{I,J} : \otimes_{i}(\otimes_{j}Z_{ij})
\to \otimes_{j}(\otimes_{i}Z_{ij})$ is the braiding morphism. 
\end{prop}

This is proved in \cite{Enr:univ,Pos}; see also Appendix \ref{app:A}.

\subsection{A filtration on $\on{QLBA}$}

 Let $\langle \varphi\rangle$ be the prop ideal of $\on{QLBA}$
 generated by $\varphi$ and by $\langle \varphi\rangle^{n}$ is $n$th power. 
 For $X,Y\in \on{Ob(Sch)}$, we have a decreasing filtration 
 $\on{QLBA}(X,Y) \supset \langle\varphi\rangle(X,Y)
 \supset \langle\varphi\rangle^{2}(X,Y)\supset...$. 
 As $\varphi$ is homogeneous for the $\ZZ^{2}$-grading, so are the 
 $\langle\vp\rangle^{n}$, i.e., $\langle\vp\rangle^{n}(X,Y)
 =\oplus_{\alpha\in\ZZ^{2}}\langle\vp\rangle^{n}(X,Y)[\alpha]$. 
 
\begin{lemma} 
This filtration is complete, i.e., $\cap_{n\geq 0} \langle\vp\rangle^{n}(X,Y)=0$. 
 \end{lemma}
 
{\em Proof.} Observe that $\langle \vp\rangle^{n}(X,Y)$ is supported in 
 $n(-1,2)+\NN(1,0)+\NN(0,1)+\NN(-1,2)\subset (2n+\NN)(0,1)+\ZZ(1,0)$. 
 Then $\cap_{n\geq 0} \langle\vp\rangle^{n}(X,Y)$ is supported
 in $\cap_{n\geq 0}(2n+\NN)(0,1)+\ZZ(1,0)$, which is empty. 
 So this intersection is zero. \hfill \qed\medskip

 The composition of $\on{QLBA}$ restricts to a map 
 $\langle\vp\rangle^{m}(G,H) \otimes \langle \vp\rangle^{n}(F,G)
 \to \langle \vp\rangle^{n+m}(F,H)$, and the tensor product restricts to 
 $\langle\vp\rangle^{n}(F,G) \otimes \langle\vp\rangle^{n'}(F',G')
 \to \langle\vp\rangle^{n+n'}(F\otimes F',G\otimes G')$, so 
 $\on{QLBA}\supset \langle\vp\rangle\supset...$ is a prop filtration. 
 The associated graded prop is defined by 
 $\on{gr}\on{QLBA}(F,G):= \oplus_{n\geq 0}\on{gr}_{n}\on{QLBA}(F,G)$, 
 where $\on{gr}_{n}\on{QLBA}(F,G) = \langle\vp\rangle^{n}(F,G)
 / \langle\vp\rangle^{n+1}(F,G)$. 

\subsection{The graded prop $\on{LBA}_{\alpha}$} \label{sect:alpha}
Define $P$ to be the prop with the same generators $\tilde\mu,
\tilde\delta,\tilde\varphi$ as $\on{QLBA}$ and the same relations, 
except for the third which is replaced by $\on{Alt}_{3} \circ
(\tilde\delta\otimes\on{id}_{{\bf id}}) \circ\tilde\delta=0$. 

We now construct a prop isomorphic to $P$. The following general construction 
goes back to \cite{EH}. For $C\in \on{Ob(Sch)}$, 
we have a prop $\on{LBA}_{C}$ defined by $\on{LBA}_{C}(F,G):= 
\oplus_{n\geq 0}\on{LBA}(F\otimes S^{n}(C),G)$ (the composition is 
induced by the coproduct $S\to S^{\otimes 2}$). For 
$D\in \on{Ob(Sch)}$, we set 
$\on{LBA}_{C,D}(F,G):= \oplus_{n\geq 0}\on{LBA}(F\otimes S^{n}(C)
\otimes D,G)$; for $\alpha \in \on{LBA}(C,D)$,
we have a map $\on{LBA}(F\otimes S^{n}(C)\otimes D,G)
\to \on{LBA}(F\otimes S^{n+1}(C),G)$, $x\mapsto 
x \circ [\on{id}_{F\otimes S^{n}(C)}\otimes \alpha] 
\circ [\on{id}_{F}\otimes \Delta_{n,1}]$, where $\Delta_{n,1} : 
S^{n+1}(C)\to S^{n}(C)\otimes C$ is the component $n+1 \to (n,1)$ of the 
coproduct $S(C)\to S(C)^{\otimes 2}$. We then have commutative diagrams 
$$\begin{matrix}
\begin{matrix} & \on{LBA}_{C,D}(F,G)\otimes \on{LBA}_{C}(G,H)
\\ & \oplus \on{LBA}_{C}(F,G)\otimes \on{LBA}_{C,D}(G,H) \end{matrix}
& \to & \on{LBA}_{C,D}(F,H)\\
\downarrow & & \downarrow \\
\on{LBA}_{C}(F,G)\otimes \on{LBA}_{C}(G,H) & \to & 
\on{LBA}_{C}(F,H)
\end{matrix}$$
and 
$$\begin{matrix}
\begin{matrix} & \on{LBA}_{C,D}(F,G)\otimes \on{LBA}_{C}(F',G')
\\ & \oplus \on{LBA}_{C}(F,G)\otimes \on{LBA}_{C,D}(F',G')\end{matrix}
& \to & \on{LBA}_{C,D}(F\otimes F',G\otimes G')\\
\downarrow & & \downarrow \\
\on{LBA}_{C}(F,G)\otimes \on{LBA}_{C}(F',G') &\to & 
\on{LBA}_{C}(F\otimes F',G\otimes G')
\end{matrix}$$
induced by the composition and tensor product, which implies that 
if 
$$\on{LBA}_{\alpha}(F,G):= \on{Coker}[\on{LBA}_{C,D}(F,G)
\to \on{LBA}_{C}(F,G)],$$ 
then we have a prop morphism $\on{LBA}_{C}
\to \on{LBA}_{\alpha}$. 

In what follows, we will set $C:= \wedge^{3}$, $D:= \wedge^{4}$, 
$\alpha:= pr_{4} \circ 
\on{Alt}_{4}\circ (\delta\otimes \on{id}_{{\bf id}}) \circ inj_{3}
\in \on{LBA}(\wedge^{3},\wedge^{4})$, where $inj_{3} : \wedge^{3}\to 
{\bf id}^{\otimes3}$ and $pr_{4} : {\bf id}^{\otimes 4}\to \wedge^{4}$
are the canonical injection and projection. 

\begin{lemma} \label{str:LBAalpha}
We have a prop isomorphism $\on{LBA}_{\alpha}\simeq P$. 
\end{lemma}

{\em Proof.} Let $\tilde P$ be the prop with generators 
$\tilde\mu,\tilde\delta,\tilde\varphi$ and only relations: 
Lie bialgebra relations between $\tilde\mu,\tilde\delta$ and 
$\tilde\vp={1\over 6}\on{Alt}_{3}(\tilde\vp)$. We have a 
morphism $\tilde P \to \on{LBA}_{\wedge^{3}}$, defined by 
$\tilde\mu \mapsto \mu\in  \on{LBA}({\bf id}^{\otimes 2}\otimes 
S^{0}(\wedge^{3}),{\bf id})
\subset \on{LBA}_{\wedge^{3}}({\bf id}^{\otimes 2},{\bf id})$; 
$\tilde\delta\mapsto \delta\in \on{LBA}({\bf id}\otimes S^{0}(\wedge^{3}),
{\bf id}^{\otimes 2}) \subset \on{LBA}_{\wedge^{3}}({\bf id},
{\bf id}^{\otimes 2})$; $\tilde\vp\mapsto inj_{3} \in 
\on{LBA}({\bf 1}\otimes S^{1}(\wedge^{3}),{\bf id}^{\otimes 3})\subset 
\on{LBA}_{\wedge^{3}}({\bf 1},{\bf id}^{\otimes 3})$, as $inj_{3}={1\over 6}
\on{Alt}_{3}\circ inj_{3}$. We also have a morphism 
$\on{LBA}_{\wedge^{3}}\to \tilde P$, defined by 
$\on{LBA}_{\wedge^{3}}(F,G) \supset \on{LBA}(F\otimes S^{n}(\wedge^{3}),G)
\supset f \mapsto \on{can}(f) \circ (\on{id}_{F}\otimes S^{n}(\tilde \vp))
\in \tilde P(F,G)$, where $\on{can}:\on{LBA}\to \tilde P$ is the 
prop morphism defined by $\mu,\delta\mapsto\tilde\mu,\tilde\delta$. 
One proves that these are inverse isomorphisms, which induce an isomorphism 
$\on{LBA}_{\alpha}\simeq P$. \hfill \qed\medskip 


\subsection{A graded prop morphism $\on{LBA}_{\alpha}\to
\on{gr}\on{QLBA}$}

\begin{lemma}
There is a unique prop morphism $\on{LBA}_{\alpha}\simeq P \to
\on{gr}\on{QLBA}$, defined by $P({\bf id}^{\otimes 2},{\bf id})
\ni \tilde\mu \mapsto \mu \in \on{LBA}({\bf id}^{\otimes 2},{\bf id})
= \on{gr}_{0}\on{QLBA}({\bf id}^{\otimes 2},{\bf id})$, 
$P({\bf id},{\bf id}^{\otimes 2})
\ni \tilde\delta \mapsto \delta \in \on{LBA}({\bf id},{\bf id}^{\otimes 2})
= \on{gr}_{0}\on{QLBA}({\bf id},{\bf id}^{\otimes 2})$
(we have $\on{QLBA}/\langle\vp\rangle = \on{LBA}$, so 
$\on{gr}_{0}\on{QLBA} = \on{LBA}$), 
$P({\bf 1},{\bf id}^{\otimes 3})\ni \tilde\vp\mapsto [\vp]
\in \on{gr}_{1}\on{QLBA}({\bf 1},{\bf id}^{\otimes 3})$. 
\end{lemma}

{\em Proof.} The images in $\gr^{0}\on{QLBA}$ of 
the Jacobi relation for $\mu$, of the cocycle relation 
between $\mu,\delta$, and of the quasi-co-Jacobi relation between $\mu,\delta,\varphi$
(which hold in $\langle\vp\rangle^{0} = \on{QLBA}$)
are respectively, the Jacobi relation for $[\mu]$, the cocycle relation between 
$[\mu],[\delta]$ and the co-Jacobi relation for $[\delta]$. The images in 
$\on{gr}_{1}\on{QLBA}$ of the relations $\vp = {1\over 6}\on{Alt}_{3}(\vp)$, 
$\on{Alt}_{4}((\delta\otimes\on{id}_{{\bf id}}^{\otimes 2})(\vp))=0$
(which hold in $\langle\vp\rangle$) are the similar relations, with $\delta,\vp$
replaced by $[\delta],[\vp]$. It follows that we have a prop morphism $P\to 
\on{gr}\on{QLBA}$, $\tilde\mu,\tilde\delta,\tilde\vp\mapsto 
[\mu],[\delta],[\vp]$. 
\hfill \qed\medskip 

\subsection{}

\begin{thm} \label{thm:QLBA}
The morphism $\on{LBA}_{\alpha}\to\on{gr}\on{QLBA}$ is a prop isomorphism. 
\end{thm}

{\em Proof.} We say that a prop morphism $\cC\to\cD$ is surjective (resp., injective) 
if the maps $\cC(F,G)\to \cD(F,G)$ are. 

As $\on{QLBA}$ is generated by $\mu,\delta,\vp$, the prop $\on{gr}\on{QLBA}$
is generated by their classes $[\mu],[\delta],[\vp]$, and since the 
generators of $P\simeq \on{LBA}_{\alpha}$ map to these
elements, the morphism $\on{LBA}_{\alpha}\to \on{gr}\on{QLBA}$ is surjective. 

We now prove the injectivity of 
$\on{LBA}_{\alpha}\to \on{gr}\on{QLBA}$. 
For this, we construct a filtered prop morphism 
$\on{QLBA}\to L({\bf LCA}_{\wedge^2})$; composing the 
associated graded morphism 
$\on{gr}\on{QLBA}\to L({\bf LCA}_{\wedge^2})$ with  
$\on{LBA}_{\alpha}\to \on{gr}\on{QLBA}$, we obtain 
a morphism $\on{LBA}_\alpha\to L({\bf LCA}_{\wedge^2})$. 
This morphism factors as $\on{LBA}_\alpha \to 
L({\bf LCA}_\alpha)\to L({\bf LCA}_{\wedge^2})$. 
The injectivity of $\on{LBA}_\alpha \to 
L({\bf LCA}_\alpha)$ is a consequence of a general 
argument (already used in the proof of the structure result
for the prop $\on{LBA}$, see Appendix \ref{app:A}), 
while the injectivity of the second morphism follows from that of 
a morphism $\on{LCA}_\alpha\to \on{LCA}_{\wedge^2}$, which is a 
consequence of Lemma \ref{lemma:3:4}. This establishes the injectivity of 
$\on{LBA}_\alpha \to L({\bf LCA}_{\wedge^2})$ and therefore of 
$\on{LBA}_{\alpha}\to \on{gr}\on{QLBA}$. Let us now proceed
with the details of the proof. 

We first define the auxiliary props mentioned above. 
$\on{LCA}_{\wedge^2}$ is the prop with generators
$\delta_{\on{LCA}}:{\bf id}\to {\bf id}^{\otimes 2}$, 
$r:{\bf 1}\to {\bf id}^{\otimes 2}$, and relations: 
antisymmetry and co-Jacobi for $\delta_{\on{LCA}}$, and
antisymmetry for $r$. Similarly, $\on{LCA}_\alpha$ is the 
prop with generators $\tilde{\tilde{\delta}}: {\bf id} \to 
{\bf id}^{\otimes 2}$ and $\tilde{\tilde\varphi}: 
{\bf 1}\to {\bf id}^{\otimes 3}$, 
and relations: antisymmetry and co-Jacobi for $\tilde{\tilde\delta}$, 
antisymmetry for $\tilde{\tilde\varphi}$, and $\on{Alt}_4 \circ 
(\tilde{\tilde\delta}\otimes \on{id}_{\bf id}^{\otimes 2}) 
\circ \tilde{\tilde\varphi}=0$. One checks that there are 
unique ${\bf Sch}$-props ${\bf LCA}_{\wedge^2}$, ${\bf LCA}_\alpha$ 
associated to these props (for example, ${\bf LCA}_{\wedge^2}({\bf F},
{\bf G}) = \hat\oplus_{i,j}\on{LCA}_{\wedge^2}(F_i,G_j)$ for 
${\bf F} = \hat\oplus F_i$, ${\bf G} = \hat\oplus G_i$). 
We denote by $L\in \on{Ob}({\bf Sch})$ the ``free Lie algebra'' 
Schur functor, i.e., if $V$ is a vector space, then $L(V)$ is 
the free Lie algebra generated by $V$; so $L = L_1\oplus L_2 \oplus...$, 
where $L_1={\bf id}$, $L_2=\wedge^2$, etc. 

We now define the prop morphism $\on{QLBA}\to L({\bf LCA}_{\wedge^2})$. 
The universal version of the Lie algebra bracket on $L(V)$ is
an element $\mu_{free}\in{\bf Sch}(L^{\otimes 2},L)$.  
The prop morphism $\on{QLBA}\to L({\bf LCA}_{\wedge^2})$ is then 
defined by $\mu\mapsto \mu_{free}$ (we identity $\mu_{free}$ with its
image under $L({\bf Sch})\to L({\bf LCA}_{\wedge^{2}})$), 
$\delta\mapsto \delta_{free}+ \on{ad}(r)$, where we identify 
$\delta_{free}$ with its image under $L({\bf LCA})\to L({\bf LCA}_{\wedge^{2}})$
(see Appendix \ref{app:A}), and $\on{ad}(r)\in 
{\bf LCA}_{\wedge^{2}}(L,L^{\otimes 2})$ is 
$(\mu_{free}\otimes \on{id}_{L} + (\on{id}_{L}\otimes 
\mu_{free})\circ (\beta_{L,L}\otimes \on{id}_{L}))
\circ (\on{id}_{L}\otimes \tilde r)$, $\vp\mapsto 
{1\over 2}\on{Alt}_{3} \circ [(\delta_{free}\otimes \on{id}_{L})\circ \tilde r 
- (\on{id}_{L}^{\otimes 2}\otimes \mu_{free}) \circ 
(\on{id}_{L}\otimes \beta_{L,L}\otimes \on{id}_{L})\circ 
(\tilde r\otimes \tilde r)]$ (here $\tilde r = inj_{1}^{\otimes 2}\circ r$, where
$inj_{1}:{\bf id}\to L$ is the canonical injection). 
This morphism is the propic version of the following 
construction: if $(\C,\delta_{\C})$ is a Lie coalgebra and $r_{\C}\in\wedge^{2}(\C)$, 
we consider the twist by $\wedge^{2}(inj_{1}^{\C})(r_{\C})$ 
of the Lie bialgebra $(L(\C),\delta_{L(\C)})$, where
$\delta_{L(\C)}:L(\C)\to L(\C)^{\otimes 2}$ is the unique derivation extending 
$\delta_{\C}$ (where $inj_1^\C:\C\to L(\C)$ is the canonical injection); 
this is a quasi-Lie bialgebra. 

The powers of the prop ideal $\langle r\rangle$ define a 
filtration on the prop $\on{LCA}_{\wedge^2}$; the associated graded
prop $\on{gr}\on{LCA}_{\wedge^2}$ is canonically isomorphic to 
$\on{LCA}_{\wedge^2}$. The prop morphism 
$\on{QLBA}\to L({\bf LCA}_{\wedge^2})$ is compatible with the 
filtrations (as it takes $\varphi$ to $\langle r\rangle$), and the
associated graded morphism $\on{gr}\on{QLBA}\to L({\bf LCA}_{\wedge^2})$
is given by 
$[\mu]\mapsto \mu_{free}$, $[\delta]\mapsto \delta_{free}$
and $[\varphi]\mapsto {1\over 2}\on{Alt}_3 \circ (\delta_{free}
\otimes \on{id}_L) \circ inj_1^{\otimes 2} \circ r$. 
The composed morphism 
\begin{equation} \label{comp:mor}
\on{LBA}_\alpha\to \on{gr}\on{QLBA}\to 
L({\bf LCA}_{\wedge^2})
\end{equation} 
is then given by 
$\tilde\mu\mapsto \mu_{free}$, $\tilde\delta\mapsto 
\delta_{free}$ and $\tilde\varphi\mapsto 
{1\over 2}\on{Alt}_3 \circ (\delta_{free}
\otimes \on{id}_L) \circ inj_1^{\otimes 2} \circ r$. 

We now define two prop morphisms $\on{LBA}_\alpha\to 
L({\bf LCA}_\alpha)$ and $\on{LCA}_\alpha\to 
\on{LCA}_{\wedge^2}$, such that the above morphism 
$\on{LBA}_\alpha\to L({\bf LCA}_{\wedge^2})$ coincides
with $\on{LBA}_\alpha\to L({\bf LCA}_\alpha) 
\to L({\bf LCA}_{\wedge^2})$. 

First define  
$\on{LBA}_\alpha\to L({\bf LCA}_\alpha)$. There
is a unique morphism $\on{LBA}\to L(\on{LCA})$, 
taking $\mu,\delta$ to $\mu_{free},\delta_{free}$
(see Appendix \ref{app:A}); this is the propic 
version of the functor $\{$Lie coalgebras$\}\to\{$Lie 
bialgebras$\}$, $(\C,\delta_\C)\mapsto (L(\C),$ 
free Lie bracket, unique cobracket extending $\delta_{\C})$. 
We define $\on{LBA}_\alpha\to L({\bf LCA}_\alpha)$ by 
$\tilde\mu,\tilde\delta\mapsto \mu_{free},\delta_{free}$
(we identify $\mu_{free},\delta_{free}$ with their images in 
$L({\bf LCA}_\alpha)$) and $\tilde\varphi\mapsto 
inj_1^{\otimes 2} \circ \tilde{\tilde\varphi}$. This morphism 
is the propic version of $\{(\C,\delta_{\C},\varphi_\C)|
(\C,\delta_{\C})$ is a Lie coalgebra, $\varphi_\C\in 
\wedge^3(\C), \on{Alt}_4\circ (\delta_{\C}\otimes 
\on{id}_\C^{\otimes 2})(\varphi_\C)=0\}\to \{(\A,\delta_\A,\mu_\A,
\varphi_\A)|(\A,\mu_\A,\delta_\A)$ is a Lie bialgebra, 
$\varphi_\A\in\wedge^3(\A)$, $\on{Alt}_4\circ (\delta_\A
\otimes \on{id}_\A^{\otimes 2}) \circ \varphi_\A=0\}$, extending the 
above functor by $\varphi_\A:= \wedge^3(inj_1^\C)(\varphi_\C)$. 

We then define the morphism $\on{LCA}_\alpha\to 
\on{LCA}_{\wedge^2}$ by $\tilde{\tilde\delta}\mapsto 
\delta_{\on{LCA}}$, $\tilde{\tilde\vp}\mapsto 
{1\over 2}\on{Alt}_{3}(\delta_{\on{LCA}}\otimes \on{id}_{{\bf id}})\circ r$. 
One checks that (\ref{comp:mor}) coincides with 
$\on{LBA}_\alpha\to L({\bf LCA}_\alpha) 
\to L({\bf LCA}_{\wedge^2})$. 

Let us prove that $\on{LBA}_\alpha\to L({\bf LCA}_\alpha)$
is injective. Using the symmetric group actions, this is equivalent to 
proving that for any 
$p,q\geq 0$, the map 
\begin{equation} \label{map:LBAalpha}
\on{LBA}_\alpha(T_p,T_q)\to 
{\bf LCA}_\alpha(L^{\otimes p},L^{\otimes q})
\end{equation} 
is injective. 

\begin{lemma} \label{lemma:alpha}
The map $\oplus_{Z\in \on{Irr(Sch)}} \on{LCA}_\alpha(T_p,Z)
\otimes \on{LA}(Z,T_q)\to \on{LBA}_\alpha(T_p,T_q)$, induced
by the prop morphisms $\on{LCA}_\alpha,\on{LA}\to \on{LBA}_\alpha$
($\tilde{\tilde\delta},\tilde{\tilde\varphi}
\mapsto\tilde\delta,\tilde\varphi$, $\mu\mapsto\tilde\mu$)
and by composition, is an isomorphism of vector spaces. 
\end{lemma}

{\em Proof of Lemma.}
Recall that $C=\wedge^{3}$, $D = \wedge^{4}$, $\alpha\in \on{LBA}(D,C)$. 
One may construct as above a prop $\on{LCA}_{C}$ by 
$\on{LCA}_{C}(F,G):= \oplus_{n\geq 0}\on{LCA}(F\otimes S^{n}(C),G)$; set  
$\on{LCA}_{C,D}(F,G):= \oplus_{n\geq 0}\on{LCA}(F\otimes S^{n}(C)\otimes
D,G)$, then using the fact that $\alpha\in \on{LCA}(D,C)$, 
one constructs a map $\on{LCA}_{C,D}(F,G)\to \on{LCA}_{C}(F,G)$
and one then checks that $\on{LCA}_{\alpha}(F,G)=
\on{Coker}[\on{LCA}_{C,D}(F,G)\to \on{LCA}_{C}(F,G)]$. 
For $F,G\in\on{Ob(Sch)}$, we have a commutative diagram 
$$
\begin{matrix}
\oplus_{Z\in\on{Irr(Sch)}} \on{LCA}_{C,D}(F,Z)\otimes 
\on{LA}(Z,G) & \stackrel{\simeq}{\to} & \on{LBA}_{C,D}(F,G) \\
\downarrow & & \downarrow \\
\oplus_{Z\in\on{Irr(Sch)}} \on{LCA}_{C}(F,Z)\otimes 
\on{LA}(Z,G) & \stackrel{\simeq}{\to} & \on{LBA}_{C}(F,G) 
\end{matrix}
$$
whose vertical cokernel is an isomorphism 
$$ 
\oplus_{Z\in\on{Irr(Sch)}} \on{LCA}_{\alpha}(F,Z)\otimes 
\on{LA}(Z,G) \stackrel{\simeq}{\to} \on{LBA}_{\alpha}(F,G) ;  
$$
this isomorphism coincides with the map described in the statement of the lemma. 
\hfill \qed\medskip 

We now consider the composite map 
\begin{equation} \label{long:compos}
\oplus_{Z\in\on{Irr(Sch)}} \on{LCA}_\alpha(T_p,Z)
\otimes \on{LA}(Z,T_q)\to \on{LBA}_\alpha(T_p,T_q)
\to {\bf LCA}_\alpha(L^{\otimes p},L^{\otimes q})
\to {\bf LCA}_\alpha(T_p,L^{\otimes q}), 
\end{equation}
where the first map is described in Lemma \ref{lemma:alpha}, 
the middle map is (\ref{map:LBAalpha}), and the last map is 
induced by the injection $T_p = {\bf id}^{\otimes p}
\to L^{\otimes p}$.

\begin{lemma} \label{lemma:recent}
The map (\ref{long:compos}) coincides with the composite map 
$\oplus_{Z}\on{LCA}_{\alpha}(T_{p},Z)\otimes \on{LA}(Z,T_{q}) \simeq 
\oplus_{Z} \on{LCA}_{\alpha}(T_{p},Z)\otimes {\bf Sch}(Z,L^{\otimes q})
\to {\bf LCA}_{\alpha}(T_{p},L^{\otimes q})$, where the first map is 
induced by the  isomorphism 
$\on{LA}(Z,T_{q})\simeq {\bf Sch}(Z,L^{\otimes q})$ and the second by 
composition. 
\end{lemma}

{\em Proof.} Note that the isomorphism 
$\on{LA}(Z,T_{q})\simeq {\bf Sch}(Z,L^{\otimes q})$
is proved in Appendix \ref{app:A}. By using symmetric group 
actions, one shows that it suffices to prove the above statement with 
$Z$ replaced by $T_{N}$.  We have composed prop morphisms
$\rho : \on{LCA}_{\alpha}\to 
\on{LBA}_{\alpha}\to L({\bf LCA}_{\alpha})$
and $\sigma : \on{LA}\to \on{LBA}_{\alpha}\to 
L({\bf LCA}_{\alpha})$; actually $\sigma$ factors through $\on{LA}\to 
L({\bf Sch})$. The map (\ref{long:compos}) (with $Z$
replaced by $T_{N}$) then takes $f\otimes g$ to $\sigma(g) \circ \rho(f)
\circ inj_{1}^{\otimes p}$, where $inj_{1}:{\bf id}\to L$ is the canonical 
morphism, $f\in \on{LCA}_{\alpha}(T_{p},T_{N})$, $g\in \on{Sch}(T_{N},T_{q})$, 
$\rho(f)\in {\bf LCA}_{\alpha}(L^{\otimes p},L^{\otimes N})$, 
$\sigma(g)\in {\bf Sch}(L^{\otimes N},L^{\otimes q})$. 

We have $\rho(f)\circ inj_{1}^{\otimes p} = inj_{1}^{\otimes N}\circ f$, 
as this property can be checked for $f=\delta_{\on{LCA}},r$ and is preserved
by composition and tensor products. Moreover, $\sigma(g) \circ inj_{1}^{\otimes N}
\in {\bf LA}(T_{N},L^{\otimes N})$ is the image $\tilde g$ of $g$ under
$\on{LA}(Z,T_{q})\simeq {\bf Sch}(Z,L^{\otimes q})$. It follows that 
(\ref{long:compos}) coincides with $f\otimes g\mapsto \tilde g
\circ f$, which was to be proved. 
\hfill \qed\medskip 

According to Lemma \ref{lemma:app}, the composite map 
$\oplus_{Z}\on{LCA}_{\alpha}(T_{p},Z)\otimes \on{LA}(Z,T_{q}) \simeq 
\oplus_{Z} \on{LCA}_{\alpha}(T_{p},Z)\otimes {\bf Sch}(Z,L^{\otimes q})
\to {\bf LCA}_{\alpha}(T_{p},L^{\otimes q})$ is an isomorphism, 
so Lemma \ref{lemma:recent} implies that the composite map 
(\ref{long:compos}) is an isomorphism. The first map in (\ref{long:compos})
is also an isomorphism by Lemma \ref{lemma:alpha}, so the map 
(\ref{map:LBAalpha}) is injective, which was to be proved. 

Let us now prove that $\on{LCA}_\alpha \to \on{LCA}_{\wedge^2}$
is injective. For this, we outline the structure of these props. 
We have 
\begin{align*}
& \on{LCA}_{\wedge^{2}}(T_p,T_{z}) \\
& = 
\oplus_{k\geq 0}
\big[ \oplus_{z',z_{1},...,z_{k}|z'+z_{1}+...+z_{k}=z}
\on{Ind}_{\SG_{z'}\times \SG_{z_{1}}\times...\times \SG_{z_{k}}}^{\SG_{z}}
[\on{LCA}(T_{p},T_{z'})\otimes 
\otimes_{i=1}^{k} \on{LCA}(\wedge^{2},T_{z_{i}})] \big]_{\SG_{k}}
\end{align*}
and
\begin{align*}
\on{LCA}_{\alpha}(T_{p},T_{z})
 = \oplus_{z\geq 0}
\big[\oplus_{\begin{matrix} {\scriptstyle z',z_{1},...,z_{k} |}
\\ {\scriptstyle z'+z_{1}+...+z_{k}=z} \end{matrix}} & 
\on{Ind}_{\SG_{z'}\times\SG_{z_{1}}...\times\SG_{z_{k}}}^{\SG_{z}}
[\on{LCA}(T_{p},T_{z'}) \otimes \\
&
\otimes_{i=1}^{k}\on{Coker}\{\on{LCA}(\wedge^{4},T_{z_{i}})
\to \on{LCA}(\wedge^{3},T_{z_{i}})\}]\big]_{\SG_{k}}. 
\end{align*}
The injectivity of $\on{LCA}_{\alpha}\to \on{LCA}_{\wedge^{2}}$, 
is therefore equivalent to that of $\on{Coker}(\on{LCA}(\wedge^{4},T_{z}) 
 \to \on{LCA}(\wedge^{3},T_{z}))
 \to \on{LCA}(\wedge^{2},T_{z})$; in other terms, we have a sequence 
 $\on{LCA}(\wedge^{4},T_{z})\stackrel{-\circ\on{Alt}_{4} \circ 
 (\delta\otimes \on{id}_{{\bf id}}^{\otimes 2})}{\to} 
 \on{LCA}(\wedge^{3},T_{z}) 
 \stackrel{-\circ\on{Alt}_{3}\circ (\delta\otimes \on{id}_{{\bf id}})}{\to} 
 \on{LCA}(\wedge^{2},T_{z})$ where the composite map is zero, and we 
 should prove that the homology vanishes. 

To prove this, we will prove that the second homology of the complex
\begin{equation} \label{comp:LA}
...\to \on{LA}(T_{z},\wedge^{4}) 
\stackrel{
\on{Alt}_{3}\circ (\mu\otimes \on{id}_{{\bf id}}^{\otimes 2}))
\circ -}{\to} 
\on{LA}(T_{z},\wedge^{3})
\stackrel{  
\on{Alt}_{2} \circ (\mu\otimes \on{id}_{{\bf id}})) \circ-}{\to}
\on{LA}(T_{z},\wedge^{2}) 
\stackrel{\mu\circ-}{\to} 
\on{LA}(T_{z},{\bf id}) \to 0
\end{equation}
vanishes. We will prove more generally: 

\begin{lemma} \label{lemma:3:4}
If $z\geq 2$, the complex (\ref{comp:LA}) is acyclic; if $z=1$, 
its homology is $1$-dimensional, concentrated in degree $0$. 
\end{lemma}

{\em Proof.} Let $\cL_{z}$ (resp., $\cA_{z}$) be the free Lie (resp., associative)
algebra with generators $x_{1},...,x_{z}$. The spaces are both graded by $\oplus_{i=1}^{z}
\NN\delta_{i}$, where $|x_{i}|=\delta_{i}$. For $V$ a vector space 
graded by $\oplus_{i=1}^{z}\NN\delta_{i}$, and $I\subset [1,z]$, 
we denote by $V_{I}$ the part of $V$ of degree $\sum_{i\in I}\delta_{i}$. 
We have $\on{LA}(T_{z},\wedge^{k})
\simeq (\wedge^{k}(\cL_{z}))_{[1,z]}$. This isomorphism 
takes the complex (\ref{comp:LA}) to 
\begin{equation} \label{cplx:L}
...
\stackrel{\on{Alt}_{3}\circ (\mu_{\cL}\otimes \on{id}^{\otimes 2})}{\to}
(\wedge^{3}(\cL_{z}))_{[1,z]}
\stackrel{\on{Alt}_{2}\circ (\mu_{\cL}\otimes \on{id})}{\to}
(\wedge^{2}(\cL_{z}))_{[1,z]}
\stackrel{\mu_{\cL}}{\to} (\cL_{z})_{[1,z]} \to 0
\end{equation}
where $\mu_{\cL}$ is the Lie bracket of $\cL_{z}$. 

Let $\on{Part}_{k}(I)$ be the set of $k$-partitions of a set $I$, i.e., of the $k$uples
$(I_{1},...,I_{k})$ with $\sqcup_{i=1}^{k}I_{i}=I$. The group 
$\SG_{k}$ acts on $\on{Part}_{k}([1,z])$, and we have a decomposition 
$$(\wedge^{k}(\cL_{z}))_{[1,z]} = 
\oplus_{[(I_{1},...,I_{k})]\in \on{Part}_{k}([1,z])/\SG_{k}}
(\wedge^{k}(\cL_{z}))_{[(I_{1},...,I_{k})]},  
$$ 
where the summand in the r.h.s. is the space of antisymmetric tensors in 
$\oplus_{\sigma\in\SG_{n}} (\cL_{z})_{I_{\sigma(1)}}\otimes...\otimes
(\cL_{z})_{I_{\sigma(k)}}$. 

We have a bijection $\{(I'_{1},[(I_{2},...,I_{k})]) | I'_{1}\subset [2,z], 
[(I_{2},...,I_{k})]\in \on{Part}_{k-1}([2,z]-I'_{1})/\SG_{k-1}\} \to 
\on{Part}_{k}([1,z])$, 
taking $(I'_{1},[(I_{2},...,I_{k})])$ to $[(I'_{1}\sqcup \{1\},I_{2},...,I_{k})]$. 
The inverse bijection takes $[(I_{1},...,I_{k})]$ to 
$(I_{i}-\{1\},[(I_{1},...,I_{i-1},I_{i+1},...,I_{k})])$, where 
$i\in [1,k]$ is the index such that $1\in I_{i}$. 

For $(I'_{1},[(I_{2},...,I_{k})])$ in the first set, we have an isomorphism 
$$
(\wedge^{k}(\cL_{z}))_{[(I'_{1}\sqcup \{1\},I_{2},...,I_{k})]} 
\simeq 
(\cL_{z})_{I'_{1}\sqcup \{1\}} \otimes (\wedge^{k-1}(\cL_{z}))_{[(I_{2},...,I_{k})]}
$$
(whose inverse is given by $\on{Alt}_{k}$, or, up to a factor, by the sum of all cyclic
permutations if $k$ is odd, and their alternated sum if $k$ is even), 
which gives rise to an isomorphism 
$$
(\wedge^{k}(\cL_{z}))_{[1,z]} \simeq 
\oplus_{(I'_{1},[(I_{2},...,I_{k})])}
(\cL_{z})_{I'_{1}\sqcup \{1\}}\otimes (\wedge^{k-1}(\cL_{z}))_{[(I_{2},...,I_{k})]}
\subset (\cL_{z}\otimes \wedge^{k-1}(\cL_{z}))_{[1,z]}. 
$$
We have a complex 
\begin{equation} \label{big:cplx}
...\to (\cL_{z}\otimes \wedge^{2}(\cL_{z}))_{[1,z]}
\to (\cL_{z}\otimes \cL_{z})_{[1,z]}\to (\cL_{z})_{[1,z]} \to 0, 
\end{equation}
where the differential $(\cL_{z}\otimes \wedge^{k}(\cL_{z}))_{[1,z]}\to 
 (\cL_{z}\otimes \wedge^{k-1}(\cL_{z}))_{[1,z]}$ is induced by 
 $x_{0}\otimes (x_{1}\wedge...\wedge x_{k})
\mapsto \sum_{i=1}^{k} (-1)^{i+1}[x_{0},x_{i}]\otimes 
(x_{1}\wedge...\check x_{i}...\wedge x_{k}) + 
\sum_{1\leq i<j\leq k} (-1)^{j+1}x_{0}\otimes (x_{1}\wedge ... \wedge [x_{i},x_{j}]
\wedge ... \check x_{j}...\wedge x_{k})$. If $I,J\subset [1,z]$ are disjoint, 
we have $[(\cL_{z})_{I},(\cL_{z})_{J}] \subset (\cL_{z})_{I\cup J}$, which implies 
that if 
$$
C_{k}:= \oplus_{(I'_{1},[(I_{2},...,I_{k})])}
(\cL_{z})_{I'_{1}\sqcup \{1\}} \otimes 
(\wedge^{k-1}(\cL_{z}))_{[(I_{2},...,I_{k})]}, 
$$
then 
$$
...\to C_{2}\to C_{1} \to 0
$$
is a subcomplex of (\ref{big:cplx}), isomorphic to (\ref{cplx:L}). 

For $I'\subset [2,z]$, we have an isomorphism $(\cA_{z})_{I'}\to 
(\cL_{z})_{I'\cup\{1\}}$, given by $x_{i_{1}}...x_{i_{s}}\mapsto 
[[[x_{1},x_{i_{1}}],x_{i_{2}}],...,x_{i_{s}}]$; the inverse isomorphism 
is the restriction of the map $(\cA_{z})_{I'\cup\{1\}}\to (\cA_{z})_{I'}$
taking a monomial not starting with $x_{1}$ to $0$, and a monomial 
starting with $x_{1}$ to the same monomial with the $x_{1}$ removed
(see \cite{Bbk}). 

The compatibility of these isomorphisms with the Lie bracket can be 
described as follows: for $I',I\subset[2,z]$ disjoint, 
we have a commutative diagram 
$$
\begin{matrix}
(\cA_{z})_{I'} \otimes (\cL_{z})_{I} &\to & (\cA_{z})\\
\downarrow & & \downarrow\\
(\cL_{z})_{I'\cup\{1\}}\otimes (\cL_{z})_{I} & \to & 
(\cL_{z})_{I\cup I'\cup \{1\}}
\end{matrix}
$$
where the upper horizontal map is induced by the product in $\cA_{z}$
($\cL_{z}$ being viewed as a subspace of $\cA_{z}$) and the lower horizontal map 
is induced by the Lie bracket of $\cL_{z}$.

We have a complex 
\begin{equation} \label{cplx:A}
...\to (\cA_{z}\otimes \wedge^{2}(\cL_{z}))_{[2,z]} \to 
(\cA_{z} \otimes \cL_{z})_{[2,z]}\to (\cA_{z})_{[2,z]}\to 0, 
\end{equation}
where the map $(\cA_{z}\otimes \wedge^{k}(\cL_{z}))_{[2,z]}
\to (\cA_{z}\otimes \wedge^{k-1}(\cL_{z}))_{[2,z]}$ is induced by 
$x_{0}\otimes (x_{1}\wedge...\wedge x_{k})\mapsto 
\sum_{i=1}^{k}(-1)^{i+1} x_{0}x_{i}\otimes x_{1}\wedge...\check x_{i}
...\wedge x_{k}+\sum_{1\leq i<j\leq k} (-1)^{j+1}x_{0}\otimes
x_{1}\wedge...\wedge [x_{i},x_{j}]\wedge...\check x_{j}...\wedge x_{k}$.
The isomorphisms $(\cL_{z})_{I'_{1}\cup\{1\}}\simeq (\cA_{z})_{I'_{1}}$
induce isomorphisms $C_{k}\simeq \oplus_{(I'_{1},[(I_{2},...,I_{k})])}
(\cA_{z})_{I'_{1}}\otimes (\wedge^{k-1}(\cL_{z}))_{[(I_{2},...,I_{k})]}
= (\cA_{z}\otimes \wedge^{k-1}(\cL_{z}))_{[2,z]}$, which are compatible with the
differentials. Hence the complex (\ref{cplx:A}) is isomorphic to 
$...\to C_{2}\to C_{1}\to C_{0}\to 0$. 

The complex (\ref{cplx:A}) is the degree $\delta_{2}+...+\delta_{k}$ part of the 
complex 
\begin{equation} \label{full:A}
...\to \cA_{z} \otimes \wedge^{2}(\cL_{z}) \to \cA_{z}\otimes 
\cL_{z} \to \cA_{z}\to 0, 
\end{equation}
where the differentials are defined by the same formulas. 

Define a complete increasing filtration on (\ref{full:A}) by 
$\on{Fil}^{n}[\cA_{z}\otimes \wedge^{k}(\cL_{z})] = 
(\cA_{z})_{\leq n-k} \otimes \wedge^{k}(\cL_{z})$, 
where $(\cA_{z})_{\leq n}$ is the part of degree $\leq n$
of $\cA_{z}\simeq U(\cL_{z})$ (i.e., the span of products of 
$\leq n$ elements of $\cL_{z}$). The associated graded complex is the 
sum over $n\geq 0$ of complexes 
$\wedge^{n}(\cL_{z})\to ...\to S^{n-1}(\cL_{z})\otimes 
\cL_{z}\to S^{n}(\cL_{z})\to 0$, 
which add up to the Koszul complex
$$
...\to S(\cL_{z})\otimes \wedge^{2}(\cL_{z}) 
\to S(\cL_{z})\otimes \cL_{z}\to S(\cL_{z})\to 0, 
$$
where the differential $S(\cL_{z})\otimes \wedge^{k}(\cL_{z})
\to S(\cL_{z})\otimes \wedge^{k-1}(\cL_{z})$ is 
$f\otimes (x_{1}\wedge...\wedge x_{k})\mapsto \sum_{i=1}^{k}
(-1)^{i+1}fx_{i}\otimes (x_{1}\wedge...\check x_{i}...\wedge x_{k})$. 

Now if $V$ is a vector space, the Koszul complex 
$$
C(V):= [...\to S(V)\otimes \wedge^{2}(V) \to S(V)\otimes V\to S(V)\to 0]
$$
is a sum of complexes, graded by $\NN$
(this degree corresponds to $p+q$ in $S^{p}(V)\otimes \wedge^{q}(V)$). 
It is well-known that the homology of this complex is concentrated in 
homological degree $0$ and in degree $0$, where it is equal to 
$\kk$. Recall a proof. One checks this directly when $V$ is one-dimensional; 
we have isomorphisms $C(V\oplus W)\simeq C(V)\otimes C(W)$
of $\NN$-graded complexes, which 
implies the statement when $V$ is finite dimensional. It follows that the 
Koszul complex in ${\bf Sch}$ $...\to S\otimes \wedge^{2}\to S \otimes {\bf id}
\to S\to 0$ has its homology concentrated in homological degree $0$ and degree $0$, 
where it equals ${\bf 1}$. This implies the statement in general. 

It follows that the homology of (\ref{full:A}) is concentrated in degree 
$0$, where it is equal to $\kk$; a non-trivial homology class is that of 
$1\in \cA_{z}$. It follows that the degree $\delta_{2}+...+\delta_{z}$
part of this complex is acyclic if $z\geq 2$, i.e., (\ref{cplx:A}) is acyclic if $z\geq 2$. 
The computation of the homology of (\ref{full:A}) is straightforward when $z=1$.  
\hfill \qed\medskip 

This ends the proof of Theorem \ref{thm:QLBA}. \hfill \qed\medskip

\section{Comparison of cohomology groups} \label{sec:comparison}

We now prove Theorem \ref{thm:isom}. The morphism of complexes $C_{\on{QLBA}}
\to C_{\on{LBA}}$ is surjective, so we have an exact sequence $0\to 
\on{Ker}(C_{\on{QLBA}}
\to C_{\on{LBA}})\to C_{\on{QLBA}}\to C_{\on{LBA}}\to 0$, inducing a long 
exact sequence in cohomology. The isomorphisms $H^{i}_{\on{QLBA}}\simeq 
H^{i}_{\on{LBA}}$ will then follow from the vanishing of the relative cohomology, 
i.e., the cohomology of the complex $\on{Ker}(C_{\on{QLBA}}\to C_{\on{LBA}})$. 

Note that the complex $C_{\on{QLBA}}$ has a complete descending filtration
$F^{i}(C^{p,q}_{\on{QLBA}}) := \langle\varphi\rangle^{i}(\wedge^{p},
\wedge^{q})$. The associated graded complex is the Schouten complex 
$C_{\on{LBA}_{\alpha}}$ of $\on{LBA}_{\alpha}$, equipped with the differential 
$[\tilde\mu\oplus\tilde\delta,-]$; unlike $C_{\on{QLBA}}$, this is the total complex 
of a bicomplex, as the differentials have now degrees $(1,0)$ and $(0,1)$. 
The relative complex $\on{Ker}(C_{\on{QLBA}}\to C_{\on{LBA}})$ coincides with the 
first step $F^{1}(C_{\on{QLBA}})$ of the complex $C_{\on{QLBA}}$; its associated
graded is the positive degree (in $\tilde\varphi$) part of the Schouten complex  
$C_{\on{LBA}_{\alpha}}$. To prove that the relative complex is acyclic, it
then suffices to prove that the positive degree part in $\tilde\varphi$ of 
the complex $C_{\on{LBA}_{\alpha}}$ is acyclic. 

Explicitly, recall that $C_{\on{LBA}_{\alpha}} = \oplus_{p,q}
\on{LBA}_{\alpha}(\wedge^{p},
\wedge^{q})$, and denoting by $\on{LBA}_{\alpha}^{(i)}(X,Y)$ the degree $i$
(in $\tilde\varphi$) of $\on{LBA}_{\alpha}(X,Y)$, the bicomplex
$(C_{\on{LBA}_{\alpha}},[\tilde\mu\oplus\tilde\delta,-])$ splits up as 
$\oplus_{i\geq 0} (C^{(i)}_{\on{LBA}_{\alpha}},[\tilde\mu\oplus\tilde\delta,-])$,
where $(C^{(i)}_{\on{LBA}_{\alpha}})^{p,q} = \on{LBA}_{\alpha}^{(i)}
(\wedge^{p},\wedge^{q})$, and we wish to prove that for $i>0$, the total cohomology of
$(C^{(i)}_{\on{LBA}_{\alpha}},[\tilde\mu\oplus\tilde\delta,-])$ is zero. 
For this, we will prove that the lines of this complex are acyclic. 

We will prove more generally: 

\begin{thm} \label{thm:vanishing}
Let $C,D$ be homogeneous Schur functors of positive degrees, let $\kappa\in 
\on{LCA}(C,D)$. Let $\on{LBA}_{\kappa}(X,Y):= \on{Coker}[\on{LBA}(D\otimes
X,Y)\to \on{LBA}(C\otimes X,Y)]$. Then for any $q\geq 0$, the complex
$(\on{LBA}_{\kappa}(\wedge^{p},\wedge^{q}),[\mu,-])_{p\geq 0}$ is acyclic. 
\end{thm}

{\em Proof.} Let us make this complex explicit. For $Z\in \on{Irr(Sch)}$, define
$\mu_{Z}\in\on{LA}({\bf id}\otimes Z,Z)$ and $\tilde\mu_{Z}\in 
\on{LA}(Z\otimes{\bf id},Z)$ as follows: $\mu_{T^{p}}\in 
\on{LA}({\bf id} \otimes T^{p},T_{p})$ is the universal version of
$x\otimes x_{1}\otimes...\otimes x_{p} \mapsto \sum_{i=1}^{p}
x_{1}\otimes ...\otimes [x,x_{i}]\otimes ... \otimes x_{p}$; as it is 
$\SG_{p}$-equivariant, it decomposes under $\on{LA}({\bf id}\otimes T_{p},T_{p})
\simeq \oplus_{Z,W||Z|=|W|=p} \on{LA}({\bf id}\otimes Z,W)\otimes 
\on{Vect}(\pi_{Z},\pi_{W})$ as $\oplus_{Z}\mu_{Z}\otimes \on{id}_{\pi_{Z}}$. 
We then set $\tilde\mu_{Z}:= -\mu_{Z}\circ \beta_{Z,{\bf id}}$, where
$\beta_{Z,{\bf id}} : Z\otimes{\bf id}\to {\bf id}\otimes Z$ is the braiding 
morphism. 

Then $[\mu,-] : \on{LBA}(C\otimes \wedge^{p},\wedge^{q})\to \on{LBA}(C\otimes
\wedge^{p+1},\wedge^{q})$ is the composed map $\on{LBA}(C\otimes\wedge^{p},
\wedge^{q})\to \on{LBA}(C\otimes\wedge^{p}\otimes {\bf id},\wedge^{q})
\to \on{LBA}(C\otimes \wedge^{p+1},\wedge^{q})$, where the first map 
is $x\mapsto x\circ (\on{id}_{C}\otimes \tilde\mu_{\wedge^{p}}) - 
\tilde\mu_{\wedge^{q}} \circ (x\otimes \on{id}_{{\bf id}})$ and the second map is 
$y\mapsto y\circ \on{Alt}_{p+1}$. We have a similar differential, with $C$ replaced
by $D$, and $\kappa$ induces a commutative diagram 
$$
\begin{matrix}
\on{LBA}(D\otimes\wedge^{p},\wedge^{q}) &\to &\on{LBA}(D\otimes
\wedge^{p}\otimes{\bf id},\wedge^{q}) &\to & \on{LBA}(D\otimes\wedge^{p+1},
\wedge^{q}) \\
\downarrow & &\downarrow & & \downarrow\\
\on{LBA}(C\otimes\wedge^{p},\wedge^{q}) &\to &\on{LBA}(C\otimes
\wedge^{p}\otimes{\bf id},\wedge^{q}) &\to & \on{LBA}(C\otimes\wedge^{p+1},
\wedge^{q}) 
\end{matrix}$$
The cokernel of this diagram is $\on{LBA}_{\kappa}(\wedge^{p},\wedge^{q})
\to \on{LBA}_{\kappa}(\wedge^{p}\otimes{\bf id},\wedge^{q})\to 
\on{LBA}_{\kappa}(\wedge^{p+1},\wedge^{q})$ and the composed map 
is the differential of our complex. 

Recall that for $X_{i},Y\in \on{Ob(Sch)}$, $i=1,...,n$, we have an isomorphism 
$\on{LBA}(X_{1}\otimes...\otimes X_{n},Y)\simeq \oplus_{Z_{1},...,Z_{n}
\in \on{Irr(Sch)}} \on{LCA}(X_{1},Z_{1})\otimes ...\otimes 
\on{LCA}(X_{n},Z_{n})
\otimes \on{LA}(Z_{1}\otimes ...\otimes Z_{n},Y) = \oplus_{Z_{1},...,Z_{n}}
\on{LBA}(X_{1}\otimes ...\otimes X_{n},Y)_{Z_{1},...,Z_{n}}$. 
The inverse isomorphism 
is the direct sum of the maps $c_{1}\otimes ...\otimes c_{n}\otimes a
\mapsto a\circ (c_{1}
\otimes ...\otimes c_{n})$. If $X_{i}$ is homogeneous of positive degree, 
$\on{LCA}(X_{i},{\bf 1})=0$, so the above sum may be restricted by the condition 
$|Z_{i}|>0$.  

We now define a complex $0\to C^{0}\stackrel{d^{0,1}}{\to} 
C^{1}\to...$ as follows. The analogue of the above complex 
$[\mu,-] : \on{LBA}(C\otimes \wedge^{p},\wedge^{q}) \to 
\on{LBA}(C\otimes \wedge^{p+1},\wedge^{q})$
(with $C$ replaced by $Z$) admits a subcomplex, namely 
$C^{p}_{Z,q}:= \oplus_{Z'\in\on{Irr(Sch)}}
\on{LBA}(Z\otimes \wedge^{p},\wedge^{q})_{Z,Z'}$; 
$$
d^{p,p+1}_{Z,q} : C^{p}_{Z,q}\to C^{p+1}_{Z,q}
$$
is then the restriction of the differential $[\mu,-]$. We then have an isomorphism 
between the complexes $(\on{LBA}(C\otimes\wedge^{p},\wedge^{q}),
[\mu,-])_{p\geq 0}$
and $\oplus_{Z\in\on{Irr(Sch)},|Z|>0} \on{LCA}(C,Z) \otimes (C^{p}_{Z,q},
d^{p,p+1}_{Z,q})_{p\geq 0}$. We have a similar isomorphism replacing $C$ by $D$, 
and these isomorphisms are compatible with the morphisms of complexes 
induced by $\kappa$. 
Taking cokernels, we get an isomorphism of complexes 
$$
(\on{LBA}_{\kappa}(\wedge^{p},\wedge^{q}),[\mu,-])_{p\geq 0}
\simeq \oplus_{Z\in\on{Irr(Sch)},|Z|>0} 
\on{Coker}[\on{LCA}(D,Z)\to \on{LCA}(C,Z)]
\otimes (C_{Z,q}^{p},d_{Z,q}^{p,p+1})_{p\geq 0}. 
$$

We now prove the acyclicity of $(C_{Z,q}^{p},d_{Z,q}^{p,p+1})_{p\geq 0}$, 
for any $q\geq 0$ and any $Z\in \on{Irr(Sch)}$, $|Z|>0$. To lighten notation, 
we will denote it $(C^{p},d^{p,p+1})_{p\geq 0}$. We reexpress 
this complex as follows. View $C^{p}$ as the antisymmetric part (under the action of 
$\SG_{p}$) of 
$\tilde C^{p} := \oplus_{Z_{1},...,Z_{p}\in \on{Irr(Sch)}}
\on{LBA}(Z\otimes {\bf id}^{\otimes p},\wedge^{q})_{Z,Z_{1},...,Z_{p}}
\subset \on{LBA}(Z\otimes {\bf id}^{\otimes p},\wedge^{q})$ (we may restrict this
sum by the conditions $|Z_{i}|>0$).
Define 
$$
\tilde d^{p,p+1} : \on{LBA}(Z\otimes {\bf id}^{\otimes p},\wedge^{q})\to 
\on{LBA}(Z\otimes {\bf id}^{\otimes p+1},\wedge^{q})
$$ 
by $\tilde d^{p,p+1}(x):= x \circ (\on{id}_{Z}\otimes \mu\otimes 
\on{id}_{{\bf id}}^{\otimes p-1}) \circ (\sum_{1\leq i<j\leq p+1} 
(-1)^{i+j}\beta_{ij})
+ \mu_{\wedge^{q}} \circ (\on{id}_{{\bf id}}\otimes x) \circ 
(\sum_{1\leq i\leq p+1}(-1)^{i+1}\beta_{i})$, where $\beta_{ij}$ 
is the automorphism of 
$Z\otimes {\bf id}^{\otimes p+1}$, universal version of 
$z\otimes x_{1}\otimes...\otimes x_{p+1}\mapsto 
z\otimes x_{i}\otimes x_{j}\otimes x_{1}\otimes...\check x_{i}...\check x_{j}...
\otimes x_{p+1}$, and $\beta_{i}:Z\otimes{\bf id}^{\otimes p+1}\to 
{\bf id}\otimes Z \otimes {\bf id}^{\otimes p}$ is the universal version of 
$z\otimes x_{1}\otimes...\otimes x_{p+1}\mapsto x_{i}\otimes z\otimes
x_{1}\otimes ... \check x_{i}...\otimes x_{p+1}$. Then 
$\tilde d^{p,p+1}$ restricts to $d^{p,p+1} : C^{p}\to C^{p+1}$. 

We now introduce a filtration on $C^{p}$. Let $(\tilde C^{p})^{\leq p'}
\subset \tilde C^{p}$ be the sum of all terms where $\on{card}\{i|Z_{i}
={\bf id}\} \leq p'$. This subspace is invariant under the action of $\SG_{p}$, so its
totally antisymmetric part is a subspace $(C^{p})^{\leq p'}\subset 
C^{p}$. 

\begin{lemma}
$d^{p,p+1}((C^{p})^{\leq p'}) \subset (C^{p+1})^{\leq p'+1}$. 
\end{lemma}

{\em Proof.} To prove this, we will prove that 
$\tilde d^{p,p+1}((\tilde C^{p})^{\leq p'}) \subset 
(\tilde C^{p+1})^{\leq p'+1}$. If $x\in \on{LBA}(Z\otimes 
{\bf id}^{\otimes p},\wedge^{q})_{Z,Z_{1},...,Z_{p}}$, then 
$\mu_{\wedge^{q}}\circ (\on{id}_{{\bf id}}\otimes x) \circ \beta_{i}$
is clearly in $\on{LBA}(Z\otimes {\bf id}^{\otimes p+1},\wedge^{q})_{Z,Z_{1},...,
Z_{i-1},{\bf id},Z_{i},...,Z_{p}}$. Here $\on{card}\{i|Z_{i}={\bf id}\}$ has
been increased by 1. Moreover, for any $W\in\on{Irr(Sch)}$, the image of
$\on{LCA}({\bf id},W)\to \on{LBA}({\bf id}^{\otimes 2},W)$, $c\mapsto 
c\circ \mu$ lies in $\oplus_{\stackrel{W_{1},W_{2}\in \on{Irr(Sch)},|W_{i}|>0,} 
{\stackrel{|W_{1}|+|W_{2}|=|W|+1}{}}}\on{LBA}
({\bf id}^{\otimes 2},W)_{W_{1},W_{2}}$. So $x\circ 
(\on{id}_{Z}\otimes\mu\otimes\on{id}_{{\bf id}}^{\otimes p-1}) \circ 
\beta_{ij}$ lies in 
\begin{equation} \label{temp:exp}
\oplus_{W_{1},W_{2}\in\on{Irr(Sch)}||W_{1}|+|W_{2}| = 
|Z_{1}|+1}\on{LBA}(Z\otimes {\bf id}^{\otimes p+1},\wedge^{q}
)_{Z,Z_{2},...,Z_{i},W_{1},Z_{i+1},...,Z_{j-1},W_{2},Z_{j},...,Z_{p}}.
\end{equation}
When $(W_{1},W_{2})\in \on{Irr(Sch)}$ are such that 
$|W_{i}|>0,|W_{1}|+|W_{2}|=|Z|+1$,  $\{i||W_{i}|=1\}$ is $\leq 1$ if
$|Z_{1}|>1$ and $=2$ if $|Z_{1}|=1$. So in the summands of (\ref{temp:exp}), 
$\on{card}\{i||Z_{i}|=1\}$ is increased by at most $1$. 
\hfill \qed\medskip 

It follows that the differential $d^{p,p+1}$ is compatible with the filtration 
$(C^{p})^{\leq 0} \subset (C^{p})^{\leq 1} \subset ...
\subset (C^{p})^{\leq p} = C^{p}$. To prove that it is 
acyclic, we will prove that the associated graded complex is acyclic. 
For this, we first determine this associated graded complex. 

For $p'+p''=p$, let $\tilde C^{p',p''}:= \oplus_{Z''_{1},...,Z''_{p''}\in\on{Irr(Sch)},
|Z''_{i}|\neq 0,1} \on{LBA}(Z\otimes {\bf id}^{\otimes p'} \otimes {\bf id}^{p''},
\wedge^{q})_{Z,{\bf id},...,{\bf id},Z''_{1},...,Z''_{p''}}$. Let $C^{p',p''}$ be the 
antisymmetric part of this space w.r.t. the action of $\SG_{p'}\times \SG_{p''}$. 

\begin{lemma} 
$(C^{p})^{\leq p'}/(C^{p})^{\leq p'-1} = C^{p',p''}$, where $p''=p-p'$. 
\end{lemma}

{\em Proof.} $(\tilde C^{p})^{\leq p'}/(\tilde C^{p})^{\leq p'-1} = 
\oplus_{Z_{1},...,Z_{p}\in \on{Irr(Sch)}, |Z_{i}|>0,
\on{card}\{i|Z_{i}={\bf id}\}=p'}\on{LBA}(Z\otimes {\bf id}^{\otimes p},
\wedge^{q})_{Z,Z_{1},...,Z_{p}}$. $(C^{p})^{\leq p'}/(C^{p})^{\leq p'-1}$
is the $\SG_{p}$-antiinvariant part of this space, which identifies with the 
$\SG_{p'}\times\SG_{p''}$-antiinvariant part of $\tilde C^{p',p''}$, i.e., 
$C^{p',p''}$. The isomorphism\footnote{For $M$ a module over 
$\prod_{i}\SG_{p_{i}}$, we denote by $M^{(\prod_{i}\SG_{p_{i}})^{-}}$ the 
component of $M$ of type $\otimes_{i}\eps_{i}$, where $\eps_{i}$ is the 
signature character of $\SG_{p_{i}}$.} 
$$[\oplus_{\on{card}\{i|Z_{i}={\bf id}\}}
\on{LBA}(Z\otimes {\bf id}^{\otimes p},\wedge^{q})]^{\SG_{p}^{-}} \to 
[\oplus_{|Z''_{i}|>1}\on{LBA}(Z\otimes{\bf id}^{\otimes p},\wedge^{q}
)_{Z,{\bf id},...,{\bf id},Z''_{1},...,Z''_{p''}}]^{(\SG_{p'}\times\SG_{p''})^{-}}$$
is given by projection on the relevant components, and the inverse isomorphism 
is given by the action of $(1/p!)\sum_{\sigma\in\SG_{p}}\epsilon(\sigma)\sigma$ (or 
$(p'!p''!/p!)\sum_{\sigma\in\SG_{p',p''}}\epsilon(\sigma)\sigma$, where 
$\SG_{p',p''}$ is the set of $p',p''$-shuffle permuations). 
\hfill \qed \medskip 

Define
$$
\tilde d^{p',p'+1|p''} : \on{LBA}(Z\otimes {\bf id}^{\otimes p'} \otimes 
{\bf id}^{\otimes p''},\wedge^{q}) \to \on{LBA}(Z\otimes {\bf id}^{\otimes p'+1}
\otimes {\bf id}^{\otimes p''},\wedge^{q})
$$
by $x\mapsto x \circ (\on{id}_{Z}\otimes \mu\otimes 
\on{id}_{{\bf id}}^{\otimes p-1}) \circ (\sum_{1\leq i<j\leq p'+1}
(-1)^{i+j+1}\beta_{ij}) + x \circ (\mu_{Z}\otimes 
\on{id}_{{\bf id}}^{\otimes p})\circ (\sum_{1\leq i\leq p'+1}(-1)^{i+1}\beta_{i})$. 

\begin{lemma}
The map $\tilde d^{p',p'+1|p''}$
restricts to maps $\tilde C^{p',p''}\to \tilde C^{p'+1,p''}$ and 
$d^{p',p'+1|p''} : C^{p',p''}\to C^{p'+1,p''}$, 
and the map $(C^{p})^{\leq p'}/(C^{p})^{\leq p'-1}\to 
(C^{p+1})^{\leq p'+1}/(C^{p+1})^{\leq p'}$ induced by $d^{p,p+1}$ 
coincides with $d^{p',p'+1|p''}$, where $p''=p-p'$. 
\end{lemma}

For each $p''$, $(C^{p',p''},d^{p',p'+1|p''})_{p'\geq 0}$ is therefore a complex
(this can be checked directly); it is embedded in the similar complex, where the
restrictions $|Z''_{i}|\neq 0,1$ are dropped, which is the universal version of the 
complex computing $H^{p'}(\A,Z(\A)^{*})\otimes \wedge^{p''}(\A)^{*}
\otimes\wedge^{q}(\A)$, where $\A$ is a Lie bialgebra. 

{\em Proof.} If $x\in \on{LBA}(Z\otimes {\bf id}^{\otimes p'} \otimes 
{\bf id}^{\otimes p''},\wedge^{q})_{Z,{\bf id},...,{\bf id},Z''_{1},...,Z''_{p''}}$, 
then one checks that both 
$x\circ (\on{id}_{Z}\otimes \mu \otimes \on{id}_{{\bf id}}^{\otimes p-1})
\circ \beta_{ij}$ and $x\circ (\mu_{Z}\otimes \on{id}_{{\bf id}}^{\otimes p}) 
\circ \beta_{i}$ lie in $\on{LBA}(Z\otimes {\bf id}^{\otimes p'+1}\otimes 
{\bf id}^{\otimes p''},\wedge^{q})_{Z,{\bf id},...,{\bf id},Z''_{1},...,Z''_{p''}}$, 
which implies that $\tilde d^{p',p'+1|p''}$ induces a map $\tilde C^{p',p''}
\to \tilde C^{p'+1,p''}$. The map $\tilde d^{p',p'+1|p''}$ maps the 
$\SG_{p'}\times\SG_{p''}$-antisymmetric part of $\on{LBA}(Z\otimes 
{\bf id}^{\otimes p'}
\otimes {\bf id}^{\otimes p''},\wedge^{q})$ to its analogue with $p'$ increased by $1$, 
so it restricts to a map $d^{p',p'+1|p''} : C^{p',p''}\to C^{p'+1,p''}$. 

Let us now show that the map $C^{p',p-p'}\to C^{p'+1,p-p'}$ induced by 
$d^{p,p+1} : (C^{p})^{\leq p'}\to (C^{p+1})^{\leq p'+1}$ is $d^{p',p'+1|p-p'}$. 

Let $Z_{1},...,Z_{p}\in \on{Irr(Sch)}$ be such that $Z_{i}={\bf id}$ for $i\leq p'$
and $|Z_{i}|>1$ if $i\geq p'+1$. Let  $y\in 
\on{LBA}(Z\otimes{\bf id}^{\otimes p},\wedge^{q})_{Z,Z_{1},...,Z_{p}}$ 
be of the form $a\circ (\on{id}_{Z}\otimes c_{1}\otimes...
\otimes c_{p})$, where $c_{i}\in \on{LCA}({\bf id},Z_{i})$ and 
$a\in \on{LA}(Z\otimes Z_{1}\otimes...\otimes Z_{p},\wedge^{q})$. 
Let $x:= y \circ (\on{id}_{Z}\otimes (\sum_{\sigma\in\SG_{p'}\times
\SG_{p''}}\eps(\sigma)\sigma))$, and $\tilde x:= y \circ  (\on{id}_{Z}
\otimes (\sum_{\sigma\in \SG_{p}}\eps(\sigma)\sigma))$. Then 
$x\in C^{p',p''}$, $\tilde x\in (\tilde C^{p})^{\leq p'}$, and 
$x$ corresponds to the class of $\tilde x$ under $C^{p',p''}\simeq 
(\tilde C^{p})^{\leq p'}/(\tilde C^{p})^{\leq p'-1}$. 

Let us compute $\tilde d^{p,p+1}(\tilde x)$. We have  
\begin{align} \label{contr:ij}
 & \tilde x \circ (\on{id}_{Z}\otimes\mu\otimes 
 \on{id}_{{\bf id}}^{\otimes p-2})
\circ (\sum_{i<j}(-1)^{i+j}\beta_{ij}) \\ \nonumber 
 & = 
\sum_{1\leq i<j\leq p} (-1)^{i+j}\sum_{\sigma\in\SG_{p}}
\eps(\sigma) a \circ (\on{id}_{Z}\otimes\beta_{\sigma}) \circ 
[\on{id}_{Z}\otimes (c_{\sigma(1)}\circ\mu)\otimes c_{\sigma(2)}\otimes...
\otimes c_{\sigma(p)}] \circ \beta_{ij},
\end{align} 
where $\beta_{\sigma} : Z_{\sigma(1)}\otimes...\otimes Z_{\sigma(p)}
\to Z_{1}\otimes...\otimes Z_{p}$ is the braiding map. 

We now use the fact that if $c\in\on{LCA}({\bf id},Z)$, then 
\begin{equation} \label{c:mu}
c \circ \mu = \mu_{Z} \circ (\on{id}_{{\bf id}}\otimes c) + \tilde\mu_{Z}
\circ (c\otimes \on{id}_{{\bf id}}) + \kappa(c), 
\end{equation}
where $\kappa(c)\in \on{LBA}({\bf id}^{\otimes 2},Z)$ is such that: 

$\bullet$ $\kappa(c)\in \oplus_{|W_{1}|,|W_{2}|>1}\on{LBA}(
{\bf id}^{\otimes 2},Z)_{W_{1},W_{2}}$ if $|Z|>1$, 

$\bullet$ $\kappa(c) = -c\circ\mu$ if $Z={\bf id}$. 

(\ref{c:mu}) is proved as follows: it is obvious when $Z={\bf id}$; 
we first prove it when $Z=T_{p}$ ($p>0$) and $c = (\delta\otimes
\on{id}_{\bf id}^{\otimes p-2})\circ...\circ\delta$ (iterating the use 
of the cocycle identity); as this element
generates the $\SG_{p}$-module $\on{LCA}({\bf id},T_{p})$, this
implies the identity when $Z=T_{p}$. The case of $Z\in \on{Irr(Sch)}$, 
$|Z|=p$ is derived from there by taking isotypic components under 
the action of $\SG_{p}$. 

When $|Z_{\sigma(1)}|>1$, the contribution of $\kappa(c_{\sigma(1)})$ to 
(\ref{contr:ij})
belongs to $(\tilde C^{p+1})^{\leq p'}$. The class of (\ref{contr:ij})
in $(\tilde C^{p+1})^{\leq p'+1}/ (\tilde C^{p+1})^{\leq p'}$ 
is then the same as that of 
\begin{align} \label{contr:ij:2}
& \sum_{1\leq i<j\leq p+1} (-1)^{i+j} \sum_{\sigma\in\SG_{p}}
a \circ (\on{id}_{Z}\circ\beta_{\sigma})\circ 
(\on{id}_{Z} \otimes (\mu_{Z_{\sigma(1)}} \circ [\on{id}_{{\bf id}}
\otimes c_{\sigma(1)}]) \otimes c_{\sigma(2)}\otimes...\otimes c_{\sigma(p)})
\circ \beta_{ij}\\
& 
+ \sum_{1\leq i<j\leq p+1} (-1)^{i+j} \sum_{\sigma\in\SG_{p}}
a \circ  (\on{id}_{Z}\circ\beta_{\sigma})\circ
(\on{id}_{Z} \otimes (\tilde\mu_{Z_{\sigma(1)}} \circ 
[ c_{\sigma(1)} \otimes \on{id}_{{\bf id}}]) 
\otimes c_{\sigma(2)}\otimes...\otimes c_{\sigma(p)})
\circ \beta_{ij} \nonumber  \\
& \nonumber 
+\sum_{1\leq i<j\leq p+1} (-1)^{i+j+1}\sum_{\sigma\in\SG_{p}|\sigma(1)\in[1,p']}
\eps(\sigma) a \circ  (\on{id}_{Z}\otimes\beta_{\sigma}) \circ 
[\on{id}_{Z}\otimes (c_{\sigma(1)}\circ\mu)\otimes c_{\sigma(2)}\otimes...
\otimes c_{\sigma(p)}] \circ \beta_{ij}. 
\end{align}
The first line may be rewritten as follows. Let 
$\alpha_{j}\in\SG_{p}$ be the cycle $\alpha_{j}(1)=2$,...,
$\alpha_{j}(j-2)=j-1$, $\alpha(j-1)=1$, $\alpha_{j}(j)=j$,..., $\alpha_{j}(p)=p$. 
In terms of $\tau:= \sigma\circ\alpha_{j}$, this line expresses as
$$
\sum_{j\in [1,p+1]}\sum_{i<j} \sum_{\tau\in\SG_{p}}
(-1)^{i}\eps(\tau)
a \circ  (\on{id}_{Z}\circ\beta_{\tau})\circ
(\on{id}_{Z}\otimes c_{\tau(1)}\otimes... \otimes
[\mu_{Z_{\tau(j-1)}} \circ (\on{id}_{{\bf id}}\otimes c_{\tau(j-1)})] 
\otimes ... \otimes c_{\tau(p)}) \circ \gamma_{ij}, 
$$
where $\gamma_{ij}\in \on{Aut}(Z\otimes {\bf id}^{\otimes p})$ 
is the categorical version of $z\otimes x_{1}\otimes ...
\otimes x_{p+1}
\mapsto z\otimes x_{1}\otimes... \check x_{i}...\otimes x_{j-1}\otimes x_{i}
\otimes x_{j+1}\otimes...\otimes x_{p+1}$. In the same way, one shows that the 
second line has the same expression, with the condition $i<j$ replaced by $i>j$ and 
$\gamma_{ij}$ the categorical version of $z\otimes x_{1}\otimes... \otimes 
x_{p+1}\mapsto z\otimes x_{1}\otimes... \otimes x_{j-1}\otimes x_{i}
\otimes x_{j+1}\otimes...\check x_{i}...\otimes x_{p+1}$. 

Adding up these lines, and using the identity 
$$
\mu_{W}\circ (\on{id}_{{\bf id}}\otimes a) = 
\sum_{\alpha=1}^{k} a \circ (\on{id}_{W_{1}}\otimes ... 
\otimes \mu_{W_{\alpha}}\otimes ... \otimes \on{id}_{W_{k}}) \circ 
\beta_{\alpha}, $$
in $\on{LA}({\bf id}\otimes W_{1}\otimes...\otimes W_{k},W)$, 
where $a\in \on{LA}(W_{1}\otimes...\otimes W_{k},W)$ and 
$\beta_{\alpha}$ is the braiding ${\bf id}\otimes W_{1}\otimes ... \otimes W_{k}
\to W_{1}\otimes ...\otimes W_{\alpha-1}\otimes{\bf id}\otimes 
W_{\alpha}\otimes... \otimes W_{k}$, we express the contribution of (\ref{contr:ij})
as (last line of (\ref{contr:ij:2})) $ + \mu_{\wedge^{q}} \circ (\on{id}_{\bf id}
\otimes \tilde x) \circ (\sum_{i=1}^{p+1}(-1)^{i}\beta_{i}) 
+ \tilde x \circ (\mu_{Z}\otimes \on{id}_{{\bf id}}^{\otimes p})\circ
(\sum_{i=1}^{p+1} (-1)^{i+1}\beta_{i})$. 

The class of $\tilde d^{p,p+1}(\tilde x)$ in $(\tilde C^{p+1})^{\leq p'+1}/
(\tilde C^{p+1})^{\leq p}$ is therefore the same as that of 
(last line of (\ref{contr:ij:2})) $+\tilde x \circ (\mu_{Z}\otimes
\on{id}_{\bf id}^{\otimes p}) \circ (\sum_{i=1}^{p+1}
(-1)^{p+1}\beta_{i})$. To evaluate its image in $C^{p'+1,p''}$, 
we apply the projection of 
$\oplus_{Z_{1},...,Z_{p+1}}
\on{LBA}(Z\otimes {\bf id}^{\otimes p+1},\wedge^{q})_{Z,Z_{1},...,Z_{p+1}}$
on the sum of components with $Z_{1}=...=Z_{p'+1}={\bf id}$, 
$|Z_{p'+1}|,...,|Z_{p+1}|>1$ along the other components. 

We have $\tilde x \circ (\mu_{Z}\otimes \on{id}_{\bf id}^{\otimes p})
\circ \beta_{i} = \sum_{\sigma\in \SG_{p}} \eps(\sigma) a \circ (\on{id}_{Z}
\otimes \beta_{\sigma})
\circ (\mu_{Z}\otimes c_{\sigma(1)}\otimes...\otimes c_{\sigma(p)}) \circ
\beta_{i}$, and the summand corresponding to $\sigma$ belongs to 
$\on{LBA}(Z\otimes{\bf id}^{\otimes p},\wedge^{q})_{Z,
Z_{\sigma(1)},...,Z_{\sigma(i-1)},{\bf id},Z_{\sigma(i)},...,Z_{\sigma(p)}}$; 
the projection is the identity on the 
components such that $i\in [1,p'+1]$ and 
$\sigma\in\SG_{p'}\times\SG_{p''}$ and zero on the
other ones. So the projection of 
$\tilde x \circ (\mu_{Z}\otimes \on{id}_{\bf id}^{\otimes p})
\circ (\sum_{i=1}^{p+1}(-1)^{i+1}\beta_{i})$ is $x \circ 
(\mu_{Z}\otimes \on{id}_{\bf id}^{\otimes p})
\circ (\sum_{i=1}^{p'+1}(-1)^{i+1}\beta_{i})$.  

Let us compute the projection of the last line of (\ref{contr:ij:2}). The term in this line
corresponding to $i,j,\sigma$ belongs to $\on{LBA}(Z\otimes{\bf id}^{\otimes p},
\wedge^{q})_{Z,Z_{\sigma(2)},...,Z_{\sigma(i)},{\bf id},Z_{\sigma(i+1)},...,
Z_{\sigma(j-1)},{\bf id},Z_{\sigma(j)},...,Z_{\sigma(p)}}$. The projection is then 
the identity on the terms such that $i,j\in[1,p'+1]$ and $\sigma\in\SG_{p'}\times
\SG_{p''}$ and zero on the other ones. The projection of the last line of (\ref{contr:ij:2})
is therefore $x\circ 
(\on{id}_{Z}\otimes\mu\otimes\on{id}_{\bf id}^{\otimes p-1})\circ 
(\sum_{1\leq i<j\leq p'+1}(-1)^{i+j+1} \beta_{ij})$. 

The projection of $\tilde d^{p,p+1}(\tilde x)$ is then the sum of these 
projections, i.e., $\tilde d^{p',p'+1|p''}(x)$. \hfill \qed\medskip

The associated graded of the complex $(C^{p},d^{p,p+1})_{p\geq 0}$ 
is therefore $\oplus_{p''\geq 0} (C^{p',p''},d^{p',p'+1|p''})_{p'\geq 0}$. 

We now prove that for each $p''\geq 0$, the complex 
$(C^{p',p''},d^{p',p'+1|p''})_{p'\geq 0}$ is acyclic. 
 
For $\ul Z'' = (Z''_{1},...,Z''_{p''})\in \on{Irr(Sch)}$, let
$$
\tilde d^{p',p'+1}_{\ul Z''} : 
\on{LA}(Z\otimes {\bf id}^{\otimes p'} \otimes
(\otimes_{i}Z''_{i}),\wedge^{q}) \to 
\on{LA}(Z\otimes {\bf id}^{\otimes p'+1} \otimes
(\otimes_{i}Z''_{i}),\wedge^{q}) 
$$
be defined by the same formula as $\tilde d^{p',p'+1|p''}$, 
replacing $\on{id}_{{\bf id}}^{\otimes p-1}$, 
$\on{id}_{{\bf id}}^{\otimes p}$ by 
$\on{id}_{{\bf id}}^{\otimes p'-1} \otimes \on{id}_{\otimes_{i}Z''_{i}}$, 
$\on{id}_{{\bf id}}^{\otimes p'} \otimes \on{id}_{\otimes_{i}Z''_{i}}$. 
Let $C^{p'}_{\ul Z''}$ be the antisymmetric part of  
$\on{LA}(Z\otimes {\bf id}^{\otimes p'} \otimes
(\otimes_{i}Z''_{i}),\wedge^{q})$ (w.r.t. the action of $\SG_{p'}$). 
Then $\tilde d_{\ul Z''}^{p',p'+1}$ restricts to a differential 
$d_{\ul Z''}^{p',p'+1} : C^{p'}_{\ul Z''}\to C^{p'+1}_{\ul Z''}$; moreover,
we have an isomorphism between $(C^{p',p''},d^{p',p'+1|p''})_{p'\geq 0}$
and the antisymmetric part (w.r.t. the action of $\SG_{p''}$) of  
$$
\oplus_{Z''_{1},...,Z''_{p''}\in \on{Irr(Sch)}, |Z''_{i}|\neq 0,1}
\otimes_{i=1}^{p''}\on{LA}({\bf id},Z''_{i})
\otimes (C^{p'}_{\ul Z''},d_{\ul Z''}^{p',p'+1})_{p'\geq 0}. 
$$
Since the differential of this complex is $\SG_{p''}$-equivariant, it  
suffices to prove that each complex $(C^{p'}_{\ul Z''},d^{p',p'+1}_{\ul Z''})_{p'\geq 0}$
is acyclic. 

Let $z:=|Z|$, $N:= \sum_{i}|Z''_{i}|$; let 
$$
\tilde d^{p',p'+1}_{z,N,q} : \on{LA}({\bf id}^{\otimes z} \otimes 
{\bf id}^{\otimes p'} \otimes {\bf id}^{\otimes N},{\bf id}^{\otimes q})
\to \on{LA}({\bf id}^{\otimes z} \otimes 
{\bf id}^{\otimes p'+1} \otimes {\bf id}^{\otimes N},{\bf id}^{\otimes q})
$$
be defined by the same formula as $\tilde d^{p',p'+1}_{\ul Z''}$, replacing 
$\otimes_{i}Z''_{i}$ by ${\bf id}^{\otimes N}$ and $\mu_{Z}$ by 
$\mu_{{\bf id}^{\otimes z}}$. Let $C^{p'}_{z,N,q}$ be the antisymmetric part of 
$\on{LA}({\bf id}^{\otimes z}\otimes {\bf id}^{\otimes p'} \otimes 
{\bf id}^{\otimes N},{\bf id}^{\otimes q})$ (w.r.t. the action of $\SG_{p'}$).
Then $\tilde d^{p',p'+1}_{z,N,q}$ restricts to a differential $d^{p',p'+1}_{z,N,q} : 
C^{p'}_{z,N,q}\to C^{p'+1}_{z,N,q}$. The complex 
$(C^{p'}_{z,N,q},d^{p',p'+1}_{z,N,q})_{p'\geq 0}$ is equipped with 
a natural action of $\SG_{z}\times \prod_{i}\SG_{|Z''_{i}|}\times\SG_{q}$, and 
$(C^{p'}_{\ul Z''},d_{\ul Z''}^{p',p'+1})$ is an isotypic component of 
this action. It suffices therefore to prove that 
$(C^{p'}_{z,N,q},d^{p',p'+1}_{z,N,q})_{p'\geq 0}$ is acyclic. 
 
In what follows, we denote by $\cL(u_{1},..,u_{s})$ (resp., $\cA(u_{1},...,u_{s})$)
the free Lie (resp., associative) algebra generated by $u_{1},...,u_{p}$. 
These spaces are graded by $\oplus_{i\in[1,p]}\NN\delta_{i}$ and for 
$S\subset [1,p]$, we denote by $\cL(u_{1},..,u_{s})_{S}$, 
$\cA(u_{1},..,u_{s})_{S}$ the subspaces of degree $\oplus_{i\in S}\delta_{i}$. 
In the case of two sets of generating variables $(u_{1},...u_{s})$ and $(v_{1},...,v_{t})$, 
the spaces are graded by $\oplus_{i\in[1,s]\sqcup [1,t]}\NN\delta_{i}$ and we use
the same notation for homogeneous subspaces.

\begin{lemma} \label{lemma:tensor}
We have an isomorphism of complexes
\begin{equation} \label{iso:C}
C^{\bullet}_{z,N,q} \simeq \oplus_{
\stackrel{\sqcup_{\alpha=1}^{q}I_{\alpha}=[1,z], }
{\stackrel{\sqcup_{\alpha=1}^{q}J_{\alpha}=z+[1,N]}{}}} 
\otimes_{\alpha=1}^{q} C^{\bullet}_{|I_{\alpha}|,|J_{\alpha}|,1}.
\end{equation} 
\end{lemma}

{\em Proof.} Identify $C^{p'}_{z,N,q}$ with  
$[\cL(a_{1},...,a_{z+N},x_{1},...,x_{p'})^{\otimes q}]^{\SG_{p'}^{-}}_{[1,z+N]
\sqcup[1,p']}$, which is the part of the $q$th tensor power of $\cL(a_{1},...,x_{p'})$, 
multilinear in $a_{1},...,a_{z+N},x_{1},...,x_{p'}$, and antisymmetric in 
$x_{1},...,x_{p'}$ ($a_{1},...,a_{z}$ correspond to the $z$ factors of 
${\bf id}^{\otimes z}$, $a_{z+1},...,a_{z+N}$ to the $N$ factors of 
${\bf id}^{\otimes N}$, and $x_{1},...,x_{p'}$ to the $p'$ factors of 
${\bf id}^{\otimes p'}$). The differential $d^{p',p'+1}_{z,N,q}$ then expresses as 
\begin{align} \label{d:p:p+1}
\nonumber F(a_{1},...,a_{z+N},x_{1},...,x_{p'}) \mapsto 
& \sum_{1\leq i<j\leq p'+1} (-1)^{i+j+1}
F(a_{1},...,a_{z+N},[x_{i},x_{j}],x_{1},...\check x_{i}...
\check x_{j}...,x_{p'+1})
\\
 & +\sum_{1\leq i\leq p'+1} \sum_{z'=1}^{z}(-1)^{i+1}
F(a_{1},...,[x_{i},a_{z'}],...,a_{z+N},x_{1},...\check x_{i}...,x_{p'+1}). 
\end{align}
On the other hand, we have an isomorphism $\on{LA}({\bf id}^{\otimes N},
{\bf id}^{\otimes q}) \simeq \oplus_{I_{1}\sqcup...\sqcup I_{q}=[1,N]}
\otimes_{\alpha=1}^{q}\on{LA}({\bf id}^{\otimes |I_{\alpha}|},{\bf id})$, 
with inverse given by the sum of maps $a_{1}\otimes...\otimes a_{q}
\mapsto (a_{1}\otimes...\otimes a_{q})\circ \beta_{I_{1},...,I_{q}}$, 
where $\beta_{I_{1},...,I_{q}} : {\bf id}^{\otimes N} \to 
\otimes_{\alpha}{\bf id}^{\otimes|I_{\alpha}|}$ is the braiding induced
by the maps $[1,N]\to \sqcup_{\alpha} [1,|I_{\alpha}|]$, taking $I_{\alpha}$
to $[1,|I_{\alpha}|]$ by preserving the order. Analyzing the action of $\SG_{N}$
on the set of $q$-compositions of $[1,N]$, we derive an isomorphism 
$\on{LA}(\wedge^N,{\bf id}^{\otimes q})\simeq \oplus_{N_{1}+...+N_{q}=N}
\otimes_{\alpha=1}^{q}\on{LA}(\wedge^{N_{\alpha}},{\bf id})$, 
with inverse given by the direct sum of the maps $a_{1}\otimes...\otimes a_{q}
\mapsto (a_{1}\otimes...\otimes a_{q})\circ \beta_{N_{1},...,N_{q}}$, 
where $\beta_{N_{1},...,N_{q}} : \wedge^{N} \to 
\wedge^{N_{1}}\otimes...\otimes \wedge^{N_{q}}$ is the composite Schur morphism 
$\wedge^{n}\hookrightarrow {\bf id}^{\otimes N} \simeq 
\otimes_{\alpha} {\bf id}^{\otimes N_{\alpha}} \twoheadrightarrow 
\otimes_{\alpha}\wedge^{N_{\alpha}}$. One proves similarly that we have an 
isomorphism $\on{LA}({\bf id}^{\otimes z} \otimes \wedge^{p'} \otimes 
{\bf id}^{\otimes N},{\bf id}^{\otimes q}) \simeq 
\oplus_{\sqcup_{\alpha} I_{\alpha}=[1,z],
\sqcup_{\alpha}J_{\alpha}=z+[1,N],\sum_{\alpha}p_{\alpha}=p'}
\otimes_{\alpha=1}^{q}\on{LA}({\bf id}^{\otimes|I_{\alpha}|}
\otimes\wedge^{p_{\alpha}} \otimes {\bf id}^{\otimes|J_{\alpha}|},{\bf id})$, 
with inverse induced by the direct sum of the maps 
$\otimes_{\alpha} a_{\alpha} \mapsto (\otimes_{\alpha}a_{\alpha}) \circ 
\beta_{(I_{\alpha}),(J_{\alpha}),(p_{\alpha})}$, where 
$\beta_{(I_{\alpha}),(J_{\alpha}),(p_{\alpha})} : {\bf id}^{\otimes z}
\otimes \wedge^{p'} \otimes {\bf id}^{\otimes N} \to 
\otimes_{\alpha} ({\bf id}^{\otimes|I_{\alpha}|} \otimes \wedge^{p_{\alpha}}
\otimes {\bf id}^{\otimes |J_{\alpha}|})$ is constructed from the above Schur 
morphisms. It follows that we have the isomorphism (\ref{iso:C}) between 
graded vector spaces. Let us show that it is compatible with differentials. 

For $\sqcup_{\alpha}I_{\alpha}=[1,z]$, $\sqcup_{\alpha}J_{\alpha}=z+[1,N]$, 
$\sum_{\alpha}p_{\alpha}=p'$, the map $\otimes_{\alpha=1}^{q}
C^{p_{\alpha}}_{|I_{\alpha}|,|J_{\alpha}|,1}\to C^{p'}_{z,N,q}$
identifies with the map 
$$
\otimes_{\alpha} [\cL(a_{1},...,x_{p'})^{\SG_{p_{\alpha}}^{-}}_{I_{\alpha}
\sqcup J_{\alpha}\sqcup (p_{1}+...+p_{\alpha-1}+[1,p_{\alpha}])}]
\to [\cL(a_{1},...,x_{p'})^{\otimes q}]^{\SG_{p'}^{-}}_{[1,z+N]\sqcup [1,p']}, 
$$
$\otimes_{\alpha}F_{\alpha} \mapsto 
\sum_{\sigma\in\SG_{p_{1},...,p_{\alpha}}}\eps(\sigma)\sigma * 
(\otimes_{\alpha}F_{\alpha})$, where $\SG_{p_{1},..,p_{\alpha}}$ is the set of 
shuffle permutations of $\SG_{p'}$ (preserving the order of the elements of 
$[1,p_{1}]$, $p_{1}+[1,p_{2}]$, etc.) and $*$ is the permutation action on 
$x_{1},...,x_{p'}$. The projection $C^{p'}_{z,N,q}\to \otimes_{\alpha}
C^{p_{\alpha}}_{|I_{\alpha}|,|J_{\alpha}|,1}$ to the component indexed by 
$((I_{\alpha}),(J_{\alpha}))$ along the other components
can then be described as follows: the composite map  
\begin{align*} [\cL(a_{1},...,x_{p'})^{\otimes q}]_{[1,z+N]\sqcup [1,p']}
& \simeq \oplus_{\stackrel{\sqcup I_{\alpha}=[1,z],}
{\stackrel{\sqcup J_{\alpha}=z+[1,N],}
{\stackrel{\sqcup P_{\alpha}=[1,p']}{}}}}
\otimes_{\alpha}\cL(a_{1},...,x_{p'})_{I_{\alpha}\sqcup J_{\alpha}\sqcup 
P_{\alpha}} \\
 & \to \otimes_{\alpha} \cL(a_{1},...,x_{p'})_{I_{\alpha}\sqcup 
J_{\alpha}\sqcup (p_{1}+...+p_{\alpha-1}+[1,p_{\alpha}])}
\end{align*} (where the 
second map is the projection along all other components) restricts to 
$$[\cL(a_{1},...,x_{p'})^{\otimes q}]_{[1,z+N]\sqcup [1,p']}^{\SG_{p'}^{-}}
\to \otimes_{\alpha}[\cL(a_{1},...,x_{p'})^{\SG_{p_{\alpha}}^{-}}_{I_{\alpha}
\sqcup J_{\alpha}\sqcup (p_{1}+...+p_{\alpha-1}+[1,p_{\alpha}])}],$$ 
which identifies with the projection $C^{p'}_{z,N,q}\to \otimes_{\alpha}
C^{p_{\alpha}}_{|I_{\alpha}|,|J_{\alpha}|,1}$.  

Extending formula (\ref{d:p:p+1}) defining $d^{p',p'+1}_{z,N,q}$, 
we define a map 
$$\tilde d^{p',p'+1}_{z,N,q} : 
[\cL(a_{1},...,x_{p'})^{\otimes q}]_{[1,z+N]\sqcup [1,p']} \to 
[\cL(a_{1},...,a_{z+N},x_{1},...,x_{p'+1})^{\otimes q}]_{[1,z+N]\sqcup 
[1,p'+1]}.$$

It follows that the map $d^{p',p'+1}_{z,N,q} : C^{p'}_{z,N,q}\to C^{p'+1}_{z,N,q}$
may be identified with the composite map 
\begin{align*}
& \oplus_{\stackrel{\stackrel{\sqcup_{\alpha}I_{\alpha}=[1,z]}{
\sqcup_{\alpha}J_{\alpha}=z+[1,N]}}{\sum_{\alpha}p_{\alpha}=p'}}
\otimes_{\alpha=1}^{q} [\cL(a_{1},...,x_{p'})_{I_{\alpha},J_{\alpha},
p_{1}+...+p_{\alpha-1}+[1,p_{\alpha}]}^{\SG^{-}_{p_{\alpha}}}]
\\
 & \stackrel{\oplus (\sum_{\sigma\in \SG_{p_{1},...,p_{n}}}
\eps(\sigma)\sigma)}{\to} 
[\cL(a_{1},...,x_{p'})^{\otimes q}]_{[1,z+N]\sqcup[1,p']} \stackrel{
\tilde d^{p',p'+1}_{z,N,q}}{\to}
[\cL(a_{1},...,x_{p'+1})^{\otimes q}]_{[1,z+N]\sqcup[1,p'+1]} \\
 & \simeq 
 \oplus_{\stackrel{\sqcup_{\alpha}\tilde I_{\alpha}=[1,z]}
 {\stackrel{\sqcup_{\alpha}\tilde J_{\alpha}=z+[1,N]}{
\sqcup_{\alpha}\tilde P_{\alpha}=[1,p'+1]}}} 
 \otimes_{\alpha=1}^{q} [\cL(a_{1},...,x_{p'})_{\tilde I_{\alpha},
 \tilde J_{\alpha},
\tilde P_{\alpha}}]
\\
  & \twoheadrightarrow
 \oplus_{\stackrel{\sqcup_{\alpha}\tilde I_{\alpha}=[1,z]}
 {\stackrel{\sqcup_{\alpha}\tilde J_{\alpha}=z+[1,N]}{
 \sum_{\alpha}\tilde p_{\alpha}=p'+1}}} 
 \otimes_{\alpha=1}^{q} [\cL(a_{1},...,x_{p'})_{\tilde I_{\alpha},
 \tilde J_{\alpha},
\tilde p_{1}+...+\tilde p_{\alpha-1}+[1,\tilde p_{\alpha}]}]
\end{align*}

We have a decomposition $\tilde d^{p',p'+1}_{z,N,q} = \sum_{1\leq i<j\leq p'}
\tilde d^{ij} + \sum_{i=1}^{p'+1}\sum_{z'=1}^{z} \tilde d^{iz'}$. 
Then: 

$\bullet$ $\tilde d^{ij}$ takes the summand indexed by 
$((I_{\alpha})_{\alpha},(J_{\alpha})_{\alpha},
(P_{\alpha})_{\alpha})$ to the summand indexed by 
$((I_{\alpha})_{\alpha},(J_{\alpha})_{\alpha},
(\tilde P^{ij}_{\alpha})_{\alpha})$, where $(\tilde P^{ij}_{\alpha})_{\alpha}$
is the partition of $[1,p'+1]$ given by  
$\tilde P^{ij}_{\alpha} = [(P_{\alpha}\cap [2,i])-1]
\cup (P_{\alpha}\cap [i+1,j-1])\cup [(P_{\alpha}\cap [j,p])+1]$
if $1\notin P_{\alpha}$, and the union of the same set with $\{i,j\}$ if 
$1\in P_{\alpha}$ (all these unions are disjoint); 

$\bullet$ $\tilde d^{iz'}$ takes the summand indexed by 
$((I_{\alpha})_{\alpha},(J_{\alpha})_{\alpha},
(P_{\alpha})_{\alpha})$ to the summand indexed by 
$((I_{\alpha})_{\alpha},(J_{\alpha})_{\alpha},
(\tilde P_{\alpha}^{i,z'})_{\alpha})$, where $(\tilde P_{\alpha}^{i,z'})_{\alpha}$
is the partition of $[1,p'+1]$
given by $\tilde P^{i,z'}_{\alpha} = (P_{\alpha}\cap [1,i-1])\cup [(P_{\alpha}
\cap [i,p'])+1]$ if $z'\notin I_{\alpha}$, and the union of the same set with $\{i\}$
if $z'\in I_{\alpha}$ (all these unions are disjoint). 

As the partitions $(I_{\alpha})_{\alpha}$ and $(J_{\alpha})_{\alpha}$ of 
$[1,z]$ and $z+[1,N]$ are not modified,  (\ref{iso:C}) is a decomposition of 
complexes. If $(\tilde P_{\alpha})_{\alpha}$ is one of the partitions
$(\tilde P^{ij}_{\alpha})_{\alpha}$ ot $(\tilde P^{iz'}_{\alpha})_{\alpha}$, then 
the sequence $(|\tilde P_{\alpha}|)_{\alpha=1,...,q}$ has the form $(p_{\alpha}
+\delta_{\alpha\beta})_{\alpha}$, where $\beta\in [1,q]$ and 
$p_{\alpha}=|P_{\alpha}|$. 

Fix $\beta\in [1,q]$ and set $p^{\beta}_{\alpha}:= p_{\alpha}
+\delta_{\alpha\beta}$. The partition $(\tilde P_{\alpha})_{\alpha}$ coincides with 
$( p^{\beta}_{1}+...+ p^{\beta}_{\alpha-1}
+[1, p^{\beta}_{\alpha}])_{\alpha}$ 
if: 

(a) $P_{\alpha}=1+p_{1}+...+p_{\alpha-1}+[1,p_{\alpha}]$ if 
$\alpha<\beta$, $P_\beta = [1+p_{1}+...+p_{\beta-1}+[1,p_{\beta}-1]]\sqcup \{1\}$, 
$P_{\alpha}=p_{1}+...+p_{\alpha-1}+[1,p_{\alpha}]$ if $\alpha>\beta$
and $p_{1}+...+p_{\beta-1}+1\leq i<j\leq p_{1}+...+p_{\beta}+1$; 
in that case, $(\tilde P^{ij}_{\alpha})_{\alpha}$
is given by $\tilde P^{ij}_{\alpha}=p_{1}+...+p_{\alpha-1}+[1,p_{\alpha}]$
for $\alpha<\beta$, $\tilde P^{ij}_{\beta}=p_{1}+...+p_{\beta-1}+[1,p_{\beta}+1]$, 
and $\tilde P^{ij}_{\alpha} = 1+p_{1}+...+p_{\alpha-1}+[1,p_{\alpha}]$ if $\alpha>\beta$. 
In particular, $i<j$ belong to $\tilde P^{ij}_{\beta}$; 

(b) $P_{\alpha} = p_{1}+...+p_{\alpha-1}+[1,p_{\alpha}]$
for any $\alpha$, $p_{1}+...+p_{\beta-1}+1\leq i\leq p_{1}+...+p_{\beta}+1$
and $z'\in I_{\beta}$; in that case, $(\tilde P^{iz'}_{\alpha})_{\alpha}$ 
is given by $\tilde P^{iz'}_{\alpha}=p_{1}+...+p_{\alpha-1}+[1,p_{\alpha}]$
for $\alpha<\beta$, $\tilde P^{iz'}_{\beta}=p_{1}+...+p_{\beta-1}+[1,p_{\beta}+1]$
and $\tilde P^{iz'}_{\beta}=1+p_{1}+...+p_{\beta-1}+[1,p_{\beta}]$
for $\alpha>\beta$. In particular, $i\in \tilde P^{iz}_\beta$ and $z'\in I_{\beta}$. 

Let now $\otimes_{\alpha}F_{\alpha}(a_{1},...,x_{p'})$ belong to 
$\otimes_{\alpha} \cL(a_{1},...,x_{p'})_{I_{\alpha}\sqcup J_{\alpha}\sqcup 
(p_{1}+...+p_{\alpha-1}+[1,p_{\alpha}])}^{\SG^{-}_{p_{\alpha}}}$. 
The image of this element in $\cL(a_{1},...,x_{p'})_{[1,z+N]\sqcup [1,p']}^{\SG_{p'}^{-}}$
is $(\sum_{\sigma\in\SG_{p_{1},...,p_{q}}}\eps(\sigma)\sigma)*(\otimes_{\alpha}
F_{\alpha})$. Let us apply $d^{p',p'+1}$ to this element, and let us project the 
result to $\oplus_{\beta=1}^{q} \otimes_{\alpha}
\cL(a_{1},...,x_{p'+1})^{\SG_{p_{\alpha}^{\beta}}}_{I_{\alpha}\sqcup J_{\alpha}
\sqcup (p^{\beta}_{1}+...+p^{\beta}_{\alpha-1}+[1,p_{\alpha}^{\beta}])}$. 

According to what we have seen, the nontrivial contributions to the summand indexed by 
$\beta$ are: 

$\bullet$ for $i<j$ in $p_{1}+...+p_{\beta-1}+[1,p_{\beta}+1]$, 
the projection of $\tilde d^{ij}(\eps(\sigma)\sigma*(\otimes_{\alpha}
F_{\alpha}))$, where $\sigma$ is the shuffle permutation taking the  
$p_{1}+...+p_{\alpha-1}+[1,p_{\alpha}]$ to $P_{\alpha}$ described in (a) above; 

$\bullet$ for $i\in p_{1}+...+p_{\beta-1}+[1,p_{\beta}+1]$ and $z'\in I_{\beta}$, 
the projection of $\tilde d^{iz'}(\otimes_{\alpha}F_{\alpha})$, 
where $\tilde d^{iz'}$
is the summand of $\tilde d^{p',p'+1}_{z,N,q}$ corresponding to $(i,z')$. 

Let $d^{p_{\beta},p_{\beta}+1}_{\beta} : 
\cL(a_{1},...,x_{p'})_{I_{\beta},J_{\beta},p_{1}+...+p_{\beta-1}
+[1,p_{\beta}]}^{\SG^{-}_{p_{\beta}}} \to \cL(a_{1},...,x_{p'+1})_{I_{\beta},
J_{\beta},p_{1}+...+p_{\beta-1}+[1,p_{\beta}+1]}^{
\SG^{-}_{p_{\beta}+1}}$ be the differential 
of the complex $C^{\bullet}_{|I_{\beta}|,|J_{\beta}|,1}$ and let 
$d_{\beta}^{ij}$, $d_{\beta}^{iz'}$ be its components. 
We have $\tilde d^{ij}(\sigma*(\otimes_{\alpha}F_{\alpha})) = 
F_{1}\otimes...\otimes d_{\beta}^{ij}(F_{\beta})\otimes...\otimes 
F_{q}$ (to prove this equality, note that the $x_{1}$ 
present in the $\beta$th factor of $\sigma*
(\otimes_{\alpha}F_{\alpha})$ gets replaced by $[x_{i},x_{j}]$ in both 
sides; the signs coincide since the ``usual'' indices of variables $x_{i},x_{j}$ are shifts of 
$i,j$ by the same quantity, and this does not alter $(-1)^{i+j+1}$), while 
$\eps(\sigma)=(-1)^{p_{1}+...+p_{\beta-1}}$; on the other hand, 
$\tilde d^{iz'}(\otimes_{\alpha}F_{\alpha})=(-1)^{p_{1}+...+p_{\beta-1}}
F_{1}\otimes...\otimes
d_{\beta}^{iz'}(F_{\beta})\otimes...\otimes F_{q}$ (here the sign is due to the fact 
that the index of $x_{i}$ is, in the usual ordering, $i-(p_{1}+...+p_{\beta-1})$). 
It follows that the contribution to the summand indexed by $\beta$ is 
$(-1)^{p_{1}+...+p_{\beta-1}} F_{1}\otimes...\otimes d_{\beta}^{p_{\beta},
p_{\beta}+1}(F_{\beta})\otimes...\otimes F_{q}$. So the projection of 
$d^{p',p'+1}((\sum_{\sigma\in\SG_{p_{1},...,p_{q}}}\eps(\sigma)\sigma)*
(\otimes_{\alpha}F_{\alpha}))$ is $$
[\sum_{\beta}(-1)^{p_{1}+...+p_{\beta-1}} \on{id}\otimes... \otimes 
d^{p_{\beta},p_{\beta}+1} \otimes ...\otimes \on{id}](\otimes_{\alpha}
F_{\alpha}), 
$$
as was to be proved. 
\hfill \qed\medskip 

As $z\neq 0$, for each partition $(I_{1},...,I_{q})$ of $[1,z]$, there exists $i$
such that $|I_{i}|\neq 0$. So renaming $|I_{i}|,|J_{i}|$ by $z,N$, 
it suffices to prove that if $z\neq 0$, then $C^{\bullet}_{z,N,1}$ is acyclic. 

Recall that $C^{p'}_{z,N,1}\simeq \cL(a_{1},...,x_{p'})_{[1,z+N]\sqcup[1,p']
}^{\SG^{-}_{p'}}$ and $d_{z,N,1}^{p',p'+1}$ is given by (\ref{d:p:p+1}). On the other
hand, the map $a\mapsto \on{ad}(a)(a_{1})$ gives rise to an isomorphism 
$$
\cA^{p'}_{z,N,1} := 
\cA(a_{2},...,x_{p'})_{[2,z+N]\sqcup [1,p']}^{\SG_{p'}^{-}} \simeq 
C_{z,N,1}^{p'},
$$ 
where $\cA(u_{1},...,u_{n})$ is the free associative algebra
generated by $u_{1},..,u_{n}$ and $\on{ad}:\cA(u_{1},...,u_{s})\to 
\on{End}(\cL(u_{1},...,u_{s}))$ is the algebra morphism derived from the
adjoint action of $\cL(u_{1},...,u_{s})$ on itself.
The differential $d^{p',p'+1}_{z,N,1} : 
\cA^{p'}_{z,N,1} \to \cA^{p'+1}_{z,N,1}$ is given by 
\begin{align} \label{d:assoc}
& Q(a_{2},...,a_{z+N},x_{1},...,x_{p'}) \mapsto  
\sum_{1\leq i<j\leq p'+1} (-1)^{i+j+1} Q(a_{2},...,a_{z+N},[x_{i},x_{j}],x_{1},
...\check x_{i}...\check x_{j}...,x_{p'+1})
\\ \nonumber 
& +\sum_{i=1}^{p'+1} (-1)^{i+1}
\Big( Q(a_{2},...,a_{z+N},x_{1},...\check x_{i}...,x_{p'+1})x_{i}
+ \sum_{z'=2}^{z} Q(a_{2},...,[x_{i},a_{z'}],...,a_{z+N},
x_{1},...\check x_{i}...,x_{p'+1})\Big), 
\end{align} 
as $\on{ad}(a)([x,a_{1}])  = \on{ad}(ax)(a_{1})$, for any $x\in 
\cL(a_{1},...,x_{p'})$ and $a\in \cA(a_{1},...,x_{p'})$. 

We have an isomorphism 
$$
\cA^{p'}_{z,N,1} \simeq \oplus_{\sigma\in \on{Perm}(\{2,...,z+N\})} 
\cA_{\sigma}^{p'}, 
$$
where $\cA_{\sigma}^{p'}:= (\cA(x_{1},...,x_{p'})^{\otimes z+N}
)_{[1,p']}^{\SG^{-}_{p'}}$, whose inverse is the direct sum of the maps induced by 
$$
\otimes_{\alpha=1}^{z+N}Q_{\alpha}(x_{1},...,x_{p'}) \mapsto 
Q_{1}(x_{1},...,x_{p'})a_{\sigma(2)} Q_{2}(x_{1},...,x_{p'})a_{\sigma(3)}...
a_{\sigma(z+N)}Q_{z+N}(x_{1},...,x_{p'}). 
$$
The explicit formula (\ref{d:assoc}) shows that if $Q(a_{1},..,x_{p'})$ is a 
multilinear monomial, then the image of $Q$ by the extension of 
$d^{p',p'+1}_{z,N,1}$ given by the same formula is a linear 
combination of monomials, where the $a_{i}$ appear in the same order as in 
$Q(a_{1},...,x_{p'})$. It follows that for each $\sigma\in\on{Perm}(\{2,...,z+N\})$, 
$\cA_{\sigma}^{\bullet}$ is a subcomplex of $\cA^{\bullet}_{z,N,1}$, and that 
we have a direct sum decomposition of the complex $\cA^{\bullet}_{z,N,1}$
\begin{equation} \label{decomp:A}
\cA^{\bullet}_{z,N,1} \simeq 
\oplus_{\sigma}\cA^{\bullet}_{\sigma}. 
\end{equation}

The acyclicity of $\cA^{\bullet}_{z,N,1}$ is then a consequence of that of each
subcomplex $\cA^{\bullet}_{\sigma}$, which we now prove. Let us fix 
$\sigma\in\on{Perm}(\{2,...,z+N'\})$. There is a unique linear map 
$$\tilde d_{\sigma}^{p',p'+1} : (\cA(x_{1},...,x_{p'})^{\otimes z+N})_{[1,p']}
\to (\cA(x_{1},...,x_{p'+1})^{\otimes z+N})_{[1,p'+1]},$$ given by 
\begin{align*}
& \tilde Q(x_{1},...,x_{p'}) \mapsto \sum_{1\leq i<j\leq p'+1} (-1)^{i+j+1}
\tilde Q([x_{i},x_{j}],x_{1},...\check x_{i}...\check x_{j}...,x_{p'+1})
\\
 & + \sum_{i=1}^{p'+1} (-1)^{i+1}
\Big(
\tilde Q(x_{1},...\check x_{i}...,x_{p'+1})
[\sum_{
\stackrel{\alpha\in [\sigma^{-1}([2,z])-1]}{\stackrel{\sqcup\{z+N\}}{}}
} x_{i}^{(\alpha)}]
- [\sum_{\alpha\in \sigma^{-1}([2,z])} x_{i}^{(\alpha)}]
\tilde Q(x_{1},...\check x_{i}...,x_{p'+1}) \Big). 
\end{align*}
where $f^{(\alpha)} = 1^{\otimes \alpha-1}\otimes f \otimes 
1^{\otimes z+N-\alpha}$. If we set $\eps_{1}=0$, $\eps_{z+N+1}=1$, 
and 
\begin{equation}\label{def:epsilon}
\eps_{\alpha}=1\Leftrightarrow \sigma(\alpha)\in [2,z], \quad 
\eps_{\alpha}=0\Leftrightarrow \sigma(\alpha)\in z+[1,N],  
\end{equation}
for $\alpha\in [2,z+N]$, then this map is 
\begin{align} \label{def:d:sigma}
& \tilde Q(x_{1},...,x_{p'}) \mapsto \sum_{1\leq i<i\leq p'+1}
(-1)^{i+j+1} \tilde Q([x_{i},x_{j}],x_{1}...\check x_{i}...\check x_{j}...,x_{p'+1})
\\ & \nonumber 
+\sum_{\alpha\in[1,z+N]}\sum_{i=1}^{p'+1}
(-1)^{i+1} \big[\eps_{\alpha+1} \tilde Q(x_{1},...\check x_{i}...,x_{p'+1})
x_{i}^{(\alpha)} - \eps_{\alpha}x_{i}^{(\alpha)} \tilde Q(x_{1},...\check x_{i}
...,x_{p'+1})\big]. 
\end{align}

The map $\tilde d_{\sigma}^{p',p'+1}$
then restricts to a linear map between the subspaces of totally antisymmetric 
tensor (under the actions of $\SG_{p'}$ on the left side and $\SG_{p'+1}$ on 
the right side), which coincides with $d^{p',p'+1}_{\sigma}$. 

For $\eps,\eps'\in \{0,1\}$, define the ``elementary'' complexes 
$\cE_{\eps,\eps'}^{\bullet}$ as follows. We set $\cE_{\eps,\eps'}^{p'}:= 
\cA(x_{1},...,x_{p'})_{[1,p']}^{\SG_{p'}^{-}}$, and define
$d_{\eps,\eps'}^{p',p'+1} : \cE_{\eps,\eps'}^{p'}\to \cE_{\eps,\eps'}^{p'+1}$
by 
\begin{align} \label{def:diff}
& (d_{\eps,\eps'}^{p',p'+1}E)(x_{1},...,x_{p'+1}):= 
\sum_{1\leq i<j\leq p'+1} (-1)^{i+j+1} E([x_{i},x_{j}],x_{1},...\check x_{i}...
\check x_{j}...,x_{p'+1}) \\ \nonumber 
 & + \eps\sum_{i=1}^{p'+1} 
 (-1)^{i}x_{i}E(x_{1},...\check x_{i}...,x_{p'+1})
 + \eps' \sum_{i=1}^{p'+1}(-1)^{i+1} E(x_{1},...\check x_{i}...,x_{p'+1}) x_{i}. 
\end{align}

\begin{lemma}
For $\eps,\eps'\in \{0,1\}$, $\cE_{\eps,\eps'}^{\bullet}:= 
(\cE^{p'}_{\eps,\eps'},d^{p',p'+1}_{\eps,\eps'})_{p'\geq 0}$ is a complex. 
\end{lemma}

{\em Proof.} Note first that for any $p'\geq 0$, $\cE_{\eps,\eps'}^{p'}$ is 
$1$-dimensional, spanned by $e_{p'}(x_{1},...,x_{p'}):= 
\sum_{\sigma\in \SG_{p'}} \eps(\sigma) x_{\sigma(1)}...x_{\sigma(p')}$. 

If $\G$ is a Lie algebra, let $U(\G)_{\eps,\eps'}$ be the universal 
enveloping algebra of $\G$, equipped with the trivial $\G$-module structure 
if $(\eps,\eps')=(0,0)$, the left (resp., right) regular $\G$-module structure 
if $(\eps,\eps')=(1,0)$ (resp., $(0,1)$), and the adjoint $\G$-module structure if 
$(\eps,\eps')=(1,1)$. Let $(C_{\eps,\eps'}^{p'}(\G),
d_{\eps,\eps'}^{p',p'+1}(\G))$ be the cochain complex computing 
the cohomology of $\G$ in these modules. We have 
$C_{\eps,\eps'}^{p'}(\G) = \on{Hom}(\wedge^{p'}(\G),U(\G))$. 
There is a unique linear map $\cE_{\eps,\eps'}^{p'}\to C_{\eps,\eps'}^{p'}(\G)$, 
taking $e_{p'}$ to the composite map $\wedge^{p'}(\G) \to \G^{\otimes p'}
\to U(\G)$, where the last map is the product map, and 
one checks that the diagram 
$\begin{matrix} \cE^{p'}_{\eps,\eps'} & \stackrel{d^{p',p'+1}_{\eps,\eps'}}{\to}
& \cE^{p'+1}_{\eps,\eps'}\\
\downarrow & & \downarrow\\
C^{p'}_{\eps,\eps'}(\G)  & \stackrel{d^{p',p'+1}_{\eps,\eps'}(\G)}{\to}& 
C^{p'+1}_{\eps,\eps'}(\G)\end{matrix}$ commutes. Since 
$C^{\bullet}_{\eps,\eps'}(\G)$ is a complex, and there exists a Lie algebra
$\G$ such that the morphisms $\cE^{p'}_{\eps,\eps'}\to C^{p'}_{\eps,\eps'}(\G)$
are injective (for example, $\G$ is a free Lie algebra with countably many generators), 
$\cE^{\bullet}_{\eps,\eps'}$ is also a complex. 
\hfill \qed\medskip 

\begin{lemma} \label{A:E}
We have a isomorphism of complexes $\cA^{\bullet}_{\sigma} \simeq 
\cE^{\bullet}_{0,\eps_{2}}\otimes \cE^{\bullet}_{\eps_{2},\eps_{3}} \otimes...
\otimes \cE^{\bullet}_{\eps_{z+N},1}$, where $(\eps_{2},...,\eps_{z+N})$ is as in 
(\ref{def:epsilon}). 
\end{lemma}

{\em Proof.} The proof is parallel to that of Lemma \ref{lemma:tensor}. 
Let us set 
$$
\cA^{p'}:= \cA^{p'}_{\sigma} = 
[\cA(x_{1},...,x_{p'})^{\otimes z+N}]_{[1,p']}^{\SG_{p'}^{-}}, \quad 
\tilde \cA^{p'}:= [\cA(x_{1},...,x_{p'})^{\otimes z+N}]_{[1,p']};  
$$
if $p_{1}+...+p_{z+N}=p'$, set 
$$
\cA_{p_{1},...,p_{z+N}}:= \otimes_{\alpha=1}^{z+N}
\cA(x_{1},...,x_{p'})_{p_{1}+...+p_{\alpha-1}+[1,p_{\alpha}]}^{
\SG^{-}_{p_{\alpha}}}
$$
and if $\sqcup_{\alpha=1}^{z+N}P_{\alpha}=[1,p']$, set 
$$
\tilde\cA_{P_{1},...,P_{z+N}}:= 
\otimes_{\alpha=1}^{z+N} \cA(x_{1},...,x_{p'})_{P_{\alpha}}. 
$$
We have a decomposition $$\tilde\cA^{p'} = \oplus_{\sqcup_{\alpha}P_{\alpha}
=[1,p']}\tilde\cA_{P_{1},...,P_{z+N}}. 
$$
We will define the support of an element $x$ of $\tilde\cA^{p'}$ as the set of 
partitions $(P_{1},...,P_{z+N})$ of $[1,p']$ such that the component
$x_{(P_{1},...,P_{z+N})}$ is nonzero. 
We also have natural morphisms $\cA_{p_{1},...,p_{z+N}}\to \cA^{p'}$, 
given by $x\mapsto (\sum_{\sigma\in\SG_{p_{1},...,p_{z+N}}}\eps(\sigma)
\sigma)*x$, where $*$ is the permutation action of $\SG_{p'}$ on $x_{1},...,x_{p'}$. 
The direct sum of these morphisms gives rise to an isomorphism
$$
\oplus_{p_{1}+...+p_{z+N}=p'}\cA_{p_{1},...,p_{z+N}} \simeq \cA^{p'}.  
$$
As the l.h.s. identifies with $\oplus_{p_{1}+...+p_{z+N}=p'}
\cE_{0,\eps_{1}}^{p_{1}}\otimes...\otimes\cE_{\eps_{z+N},1}^{p_{z+N}}$, 
we obtain the identification $\cA_{\sigma}^{\bullet}\simeq 
\otimes_{\alpha=1}^{z+N} \cE^{\bullet}_{\eps_{\alpha},\eps_{\alpha+1}}$
at the level of graded vector spaces. We now show that this identification is compatible 
with the differentials. 

The composite map 
$$
\oplus_{\sum_{\alpha}p_{\alpha}=p'} \cA_{p_{1},...,p_{z+N}}
\simeq \cA^{p'} \stackrel{\on{can}}{\hookrightarrow} \tilde\cA^{p'} 
\stackrel{\pi}{\twoheadrightarrow}
\oplus_{\sum_{\alpha}p_{\alpha}=p'} \tilde\cA_{[1,p_{1}],p_{1}+[1,p_{2}],...,
p_{1}+...+p_{z+N-1}+[1,p_{z+N}]}, 
$$
where the last map is the projection along the components indexed by the 
other (non-consecutive) partitions, is the canonical inclusion map. 
It follows that the map $(\otimes_{\alpha}\cE^{\bullet}_{\eps_{\alpha},
\eps_{\alpha+1}})^{p'} \to (\otimes_{\alpha}\cE^{\bullet}_{\eps_{\alpha},
\eps_{\alpha+1}})^{p'+1}$ may be identified with the composite map 
$$
\oplus_{\sum_{\alpha}p_{\alpha}=p'}\cA_{p_{1},...,p_{z+N}}
\simeq \cA^{p'} \hookrightarrow \tilde\cA^{p'}
\stackrel{\tilde d_{\sigma}^{p',p'+1}}{\to}
\tilde\cA^{p'+1}\twoheadrightarrow 
\oplus_{\sum_{\alpha}\tilde p_{\alpha}=p'+1}
\tilde\cA_{[1,\tilde p_{1}],...,\tilde p_{1}+...+\tilde p_{z+N-1}+[1,
\tilde p_{z+N}]}. 
$$
Let now $Q_{\alpha}\in \cE^{p_{\alpha}}_{\eps_{\alpha},\eps_{\alpha+1}}
\simeq \cA(x_{1},...,x_{p'})_{p_{1}+...+p_{\alpha-1}+[1,p_{\alpha}]}^{
\SG^{-}_{p_{\alpha}}}$ and $Q:= \otimes_{\alpha}Q_{\alpha}\in 
\cA_{p_{1},...,p_{z+N}}$. The image of this element in $\tilde\cA^{p'}$ is 
$(\sum_{\sigma\in\SG_{p_{1},...,p_{z+N}}}\eps(\sigma)\sigma)*Q$. 
The summand $\eps(\sigma)\sigma*Q$ belongs to $\tilde\cA_{P_{1}^{\sigma},...,
P_{z+N}^{\sigma}}$, where $P_{\alpha}^{\sigma}:= 
\sigma(p_{1}+...+p_{\alpha-1}+[1,p_{\alpha}])$.

Decompose $\tilde d^{p',p'+1}_{\sigma} : \tilde\cA^{p'}\to \tilde\cA^{p'+1}$
as a sum $\tilde d^{p',p'+1}_{\sigma} = \sum_{1\leq i<j\leq p'+1}
\tilde d^{ij} + \sum_{i=1}^{p'+1}\sum_{\alpha\in[1,z+N]} \tilde d^{i\alpha}$. 
If $\sqcup_{\alpha}P_{\alpha}=[1,p']$, then $\tilde d^{ij}(\tilde\cA_{P_{1},...,P_{\alpha}})
\subset \tilde\cA_{P_{1}^{ij},...,P_{z+N}^{ij}}$ and $\tilde d^{i\alpha}(
\tilde\cA_{P_{1},...,P_{z+N}}) \subset \cA_{P_{1}^{i\alpha},...,P_{z+N}^{i\alpha}}$, 
where

$\bullet$ $(P_{1}^{ij},...,P_{z+N}^{ij})$ is given by 
$P^{\alpha}_{ij}:= [(P_{\alpha}\cap [2,i])-1]\sqcup 
(P_{\alpha}\cap [i+1,j-1]) \sqcup [(P_{\alpha}\cap [j,p'])+1]$
if $1\notin P_{\alpha}$, and the union of the same set with $\{i,j\}$
if $1\in P_{\alpha}$; 

$\bullet$ $(P_{1}^{i\alpha},...,P_{z+N}^{i\alpha})$ is given by 
$P_{\gamma}^{i\alpha} = (P_{\gamma}\cap [1,i-1])\sqcup 
[(P_{\gamma}\cap [i,p'])+1]$ if $\gamma\neq\alpha$, and the union of the 
same set with $\{i\}$ if $\gamma=\alpha$. 

Note that the sequences $(|P_{1}^{ij}|,...,|P_{z+N}^{ij}|)$ and $(|P_{1}^{i\alpha}|,...,
|P_{z+N}^{i\alpha}|)$ are necessarily of the form $(p_{1}^{\beta},...,
p_{z+N}^{\beta}):= (p_{1}+\delta_{1\beta},...,
p_{z+N}+\delta_{z+N,\beta})$, where $\beta\in [1,z+N]$ is the index
such that $1\in P_{\alpha}$ in the first case, and $\alpha$ in the second case. Then: 

(a) for any $i,j$ ($1\leq i<j\leq p'+1$) and any $\beta\in [1,z+N]$, 
$(P_{1}^{ij},...,P_{z+N}^{ij})$ coincides with $([1,p_{1}^{\beta}],...,
p_{1}^{\beta}+...+p_{z+N-1}^{\beta}+[1,p_{z+N}^{\beta}])$
if $P_{\alpha}=1+p_{1}+...+p_{\alpha-1}+[1,p_{\alpha}]$
for $\alpha<\beta$, $P_{\beta}= (1+p_{1}+...+p_{\beta-1}+[1,p_{\beta}-1])
\sqcup\{1\}$, and $P_{\alpha}=p_{1}+...+p_{\alpha-1}+[1,p_{\alpha}]$
for $\beta>\alpha$, and $i,j\in p_{1}+...+p_{\beta-1}+[1,p_{\beta}+1]$; 

(b) for $i\in [1,p'+1]$ and $\alpha\in[1,z+N]$, 
$(P_{1}^{i\alpha},...,P_{z+N}^{i\alpha})$ coincides with 
$([1,p_{1}^{\alpha}],...,
p_{1}^{\alpha}+...+p_{z+N-1}^{\alpha}+[1,p_{z+N}^{\alpha}])$ if 
$P_{\alpha}=p_{1}+...+p_{\alpha-1}+[1,p_{\alpha}]$ for any $\alpha$
and $i\in p_{1}+...+p_{\alpha-1}+[1,p_{\alpha+1}]$.  

If $i,j$ are such that $1\leq i<j\leq p'+1$, then the condition 
on $\sigma\in \SG_{p_{1},...,p_{z+N}}$ for the support 
$\tilde d^{ij}(\eps(\sigma)\sigma*Q)$ to consist in a consecutive partition of 
$[1,p'+1]$ is therefore: there exists $\beta\in [1,z+N]$ such that $i,j
\in p_{1}+...+p_{\beta-1}+[1,p_{\beta}+1]$, and $\sigma$
is the shuffle permutation taking $[1,p_{1}],p_{1}+[1,p_{2}],...,p_{1}+...+
p_{z+N-1}+[1,p_{z+N}]$ to the partition described in (a) above. 

If $i\in [1,p'+1]$ and $\alpha\in [1,z+N]$, 
then the condition on $\sigma\in\SG_{p_{1},...,p_{z+N}}$
for the support of $\tilde d^{i\alpha}(\eps(\sigma)\sigma*Q)$ to consist in a 
consecutive partition of $[1,p'+1]$ is therefore: $\sigma = \on{id}$ and $i\in 
p_{1}+...+p_{\alpha-1}+[1,p_{\alpha}+1]$. 

In the first case, we have $\eps(\sigma)=(-1)^{p_{1}+...+p_{\beta-1}}$
and  
$\pi\circ \tilde d^{ij}(Q) = Q_{1}\otimes...\otimes \tilde d^{ij}_{\beta}
(Q_{\beta}) \otimes...\otimes Q_{z+N}$; in the second case,  
$\pi\circ\tilde d^{i\beta}(Q)$ is $ (-1)^{p_{1}+...+p_{\beta-1}}
Q_{1}\otimes...\otimes \tilde d_{\eps_{\beta},\eps_{\beta+1}}^{i}
(Q_{\beta}) \otimes...\otimes Q_{z+N}$. Here 
$\tilde d_{\eps_{\beta},\eps_{\beta+1}}^{p_{\beta},p_{\beta+1}} : 
\cA(x_{1},...,x_{p'})_{p_{1}+...+p_{\beta-1}+[1,p_{\beta}]}
\to \cA(x_{1},...,x_{p'+1})_{p_{1}+...+p_{\beta-1}+[1,p_{\beta+1}]}$
is decomposed
as $\tilde d_{\eps_{\beta},\eps_{\beta+1}}^{p_{\beta},p_{\beta+1}} 
= \sum_{p_{1}+...+p_{\beta-1}+1
\leq i<j\leq p_{1}+...+p_{\beta}+1} \tilde d_{\beta}^{ij} + 
\sum_{i\in p_{1}+...+p_{\beta-1}+[1,p_{\beta}+1]}
 \tilde d^{i}_{\eps_{\beta},\eps_{\beta+1}}$. 

Then 
$$
\pi\circ \tilde d_{\sigma}^{p',p'+1}\circ \on{can}(Q) = 
\sum_{\beta=1}^{z+N} (-1)^{p_{1}+...+p_{\beta-1}}
Q_{1}\otimes...\otimes \tilde d_{\eps_{\beta},\eps_{\beta+1}}^{
p_{\beta},p_{\beta}+1}(Q_{\beta})
\otimes...\otimes Q_{z+N}, 
$$
which proves our claim. 
\hfill \qed\medskip 

\begin{prop} \label{E:acyclic}
The complexes $\cE^{\bullet}_{0,1}$ and $\cE^{\bullet}_{1,0}$ 
are acyclic; moreover, for $\eps\in\{0,1\}$, 
$H^{p'}(\cE^{\bullet}_{\eps,\eps})$ is zero for any 
$p'\neq 0$ and $\kk$ for $p'=0$. 
\end{prop}

{\em Proof.} If $u_{1},...,u_{n}$ are free variables, let 
$\kk = \cA_{\leq 0}(u_{1},...,u_{n}) \subset ...
\subset \cA_{\leq i}(u_{1},...,u_{n}) \subset ...\subset \cA(u_{1},...,u_{n})$ 
be the increasing PBW filtration of $\cA(u_{1},...,u_{n})$, induced by its identification 
with $U(\cL(u_{1},...,u_{n}))$. The symmetrization isomorphism 
$\cA(u_{1},...,u_{n}) \simeq S(\cL(u_{1},...,u_{n}))$ identifies 
$\cA_{\leq i}(u_{1},...,u_{n})$ with $\oplus_{i'\leq i} 
S^{i'}(\cL(u_{1},...,u_{n}))$. 
The graded space associated to this filtration is the free Poisson algebra
$\cP(u_{1},...,u_{n}) = S(\cL(u_{1},...,u_{n}))$; its degree $i$ part is
$\cP[i](u_{1},...,u_{n}) = S^{i}(\cL(u_{1},...,u_{n}))$. 

Define a filtration on $\cE^{\bullet}_{\eps,\eps'}$ by 
$F_{u}(\cE^{p'}_{\eps,\eps'}):= 
\cA_{\leq u}(x_{1},...,x_{p'})_{[1,p']}^{\SG_{p'}^{-}}$ for $u\geq 0$. 
If $E(x_1,...,x_{p'})\in \cA_{\leq u}(x_1,...,x_{p'})_{[1,p']}$, then: 
$E([x_i,x_j],x_1,...\check x_i...\check x_j...,x_{p'+1})\in 
\cA_{\leq u}(x_1,...,x_{p'+1})_{[1,p'+1]}$;  
$$x_i E(x_1,...\check x_i...,x_{p'+1}), 
E(x_1,...\check x_i...,x_{p'+1})x_i
\in \cA_{\leq u+1}(x_1,...,x_{p'+1})_{[1,p'+1]},$$ while 
$[x_i,E(x_1,...\check x_i...,x_{p'+1})] 
\in \cA_{\leq u}(x_1,...,x_{p'+1})_{[1,p'+1]}$. 
It follows that for $\eps\in \{0,1\}$, we have 
$$d_{\eps,\eps}^{p',p'+1}(F_{u}(\cE^{p'}_{\eps,\eps})) \subset 
F_{u}(\cE^{p'+1}_{\eps,\eps}),$$ 
while for $\eps\neq\eps'$ in $\{0,1\}$, 
$$d_{\eps,\eps'}^{p',p'+1}(F_{u}(\cE^{p'}_{\eps,\eps'})) \subset 
F_{u+1}(\cE^{p'+1}_{\eps,\eps'}). 
$$ 
The associated graded complex is $\cP_{\eps,\eps'}^{\bullet}$, 
where 
\begin{equation} \label{def:P}
\cP_{\eps,\eps'}^{p'} 
= \cP(x_{1},...,x_{p'})_{[1,p']}^{\SG_{p'}^{-}}
= \oplus_{u\geq 0} 
\cP[u](x_{1},...,x_{p'})_{[1,p']}^{\SG_{p'}^{-}},
\end{equation}
with differential 
$$
\on{gr}d_{\eps,\eps'}^{p',p'+1} : \cP^{p'}_{\eps,\eps'}\to 
\cP^{p'+1}_{\eps,\eps'}$$
given by 
\begin{align*}
(\on{gr}d^{p',p'+1}_{\eps,\eps}P)(x_{1},...,x_{p'+1}):= 
& \sum_{1\leq i<j\leq p'+1} (-1)^{i+j+1} P(\{x_{i},x_{j}\},x_{1},...\check x_{i}...
\check x_{j}...,x_{p'+1}) \\
 & + \eps\sum_{i=1}^{p'+1} (-1)^{i}
\{x_{i},P(x_{1},...\check x_{i}...,x_{p'+1})\}
\end{align*}
for $\eps\in \{0,1\}$, 
$$
(\on{gr}d^{p',p'+1}_{0,1}P)(x_{1},...,x_{p'+1}):= 
\sum_{i=1}^{p'+1} (-1)^{i+1}x_{i}P(x_{1},...\check x_{i}...,x_{p'+1}), 
$$
and $\on{gr}d^{p',p'+1}_{1,0} = -\on{gr}d^{p',p'+1}_{0,1}$
(when $\eps'=\eps$, the commutators give rise to brackets in the 
associated graded differential, while if $\eps\neq \eps'$, 
the only part of the differential with nontrivial contribution to the 
associated graded differential is the second line of (\ref{def:diff})). 
The differentials $\on{gr}d_{\eps,\eps}^{p',p'+1}$ have degree $0$, and 
the differentials $\on{gr}d_{\eps,\eps'}^{p',p'+1}$ have degree $1$
(if $\eps'\neq \eps$) with respect to the $\NN$-grading on 
$\cP^{\bullet}_{\eps,\eps'}$ induced by (\ref{def:P}). 
We therefore have direct sum decompositions 
$$\cP_{\eps,\eps}^{\bullet} = \oplus_{u\in\ZZ} \cP_{\eps,\eps}^{\bullet}[u],
\quad  
\cP_{\eps,\eps'}^{\bullet} = \oplus_{u\in\NN} \cP_{\eps,\eps'}^{\bullet}\{u\} 
\on{\ (if\ }\eps'\neq \eps), 
$$
where for any $\eps,\eps'$, we set 
$\cP_{\eps,\eps'}^{p'}[u]:= \cP[u](x_{1},...,x_{p'})_{[1,p']}^{\SG_{p'}^{-}}$
and $\cP_{\eps,\eps'}^{p'}\{u\} = \cP_{\eps,\eps'}^{p'}[u+p']$. 

\begin{lemma}
For $n,u\geq 0$, $\cP_{\eps,\eps'}^{n}[u]$ have the following values: 

$\bullet$ if $n=2m$, $\cP_{\eps,\eps'}^{2m}[m]$ is $1$-dimensional, 
spanned by 
$$p_{2m}(x_{1},...,x_{2m}):= \sum_{\sigma\in \SG_{2,...,2}} \eps(\sigma)
\{x_{\sigma(1)},x_{\sigma(2)}\}...\{x_{\sigma(2m-1)},x_{\sigma(2m)}\}$$
and $\cP_{\eps,\eps'}^{2m}[u]=0$ for $u\neq m$; 

$\bullet$ if $n=2m+1$, $\cP_{\eps,\eps'}^{2m+1}[m+1]$ is $1$-dimensional, 
spanned by $$p_{2m+1}(x_{1},...,x_{2m+1}):= 
\sum_{\sigma\in\SG_{1,2,...,2}} \eps(\sigma) x_{\sigma(1)}
\{x_{\sigma(2)},x_{\sigma(3)}\}...\{x_{\sigma(2m)},x_{\sigma(2m+1)}\},$$ 
and $\cP_{\eps,\eps'}^{2m+1}[u]=0$ for $u\neq m+1$. 
\end{lemma}

{\em Proof of Lemma.} As the category of $\SG_{n}$-modules is semisimple, 
the $\SG_{n}$-modules $\cA(x_{1},...,x_{n})_{[1,n]}$ and 
$\cP(x_{1},...,x_{n})_{[1,n]}$ are equivalent. It follows that 
$\cP(x_{1},...,x_{n})_{[1,n]}^{\SG_{n}^{-}}$ is $1$-dimensional. 
Since this space is equal to $\oplus_{u\geq 0} \cP[u](x_{1},...,x_{n}
)_{[1,n]}^{\SG_{n}^{-}}$, it follows that exactly one of these summands is 
$1$-dimensional, and the others are zero. It then remains to prove that 
$p_{n}\in \cP^{n}_{\eps,\eps'}[[(n+1)/2]]$ and $p_{n}\neq 0$, where $[x]$
is the integral part of $x$.

If $n=2m$, we have $p_{2m}(x_{1},...,x_{2m}) = 2^{-m}\sum_{\sigma\in
\SG_{2m}}\eps(\sigma)\{x_{\sigma(1)},x_{\sigma(2)}\}...\{x_{\sigma(2m-1)},
x_{\sigma(2m)}\}$, so $p_{2m}$ is $\SG_{n}$-antiinvariant; so 
$p_{2m}\in \cP^{2m}_{\eps,\eps'}[m]$; and if $\Gamma$ is the set of 
$\sigma\in\SG_{2m}$, such that $\sigma(1)<\sigma(3)<...<\sigma(2m-1)$
and $\sigma(2i+1)<\sigma(2i+2)$ for $i=0,...,m-1$ (this identifies with the set
of partitions of $[1,2m]$ in subsets of cardinality 2, modulo permutation of the subsets), 
we have 
$p_{2m}(x_{1},...,x_{2m}) = m! \sum_{\sigma\in\Gamma}\eps(\sigma)
\{x_{\sigma(1)},x_{\sigma(2)}\}...\{x_{\sigma(2m-1)},x_{\sigma(2m)}\}$, 
and as the summands in this expression are linearly independant, $p_{2m}\neq 0$. 

If $n=2m+1$, we have similarly $$p_{2m+1}(x_{1},...,x_{2m+1}) = 2^{-m}
\sum_{\sigma\in\SG_{2m+1}}\eps(\sigma) x_{\sigma(1)}\{x_{\sigma(2)},
x_{\sigma(3)}\}...\{x_{\sigma(2m)},x_{\sigma(2m+1)}\},$$ which implies that 
$p_{2m+1}$ is $\SG_{n}$-antiinvariant, and $$p_{2m+1}(x_{1},...,x_{2m+1})
= m!\sum_{\sigma\in\Gamma}\eps(\sigma)x_{\sigma(1)}\{x_{\sigma(2)},
x_{\sigma(3)}\}...\{x_{\sigma(2m)},x_{\sigma(2m+1)}\},$$ where $\Gamma$
is the set of permutations $\sigma\in\SG_{n}$ such that $\sigma(2)<\sigma(4)<...
<\sigma(2m)$ and $\sigma(2i)<\sigma(2i+1)$ for $i=1,...,m$, which implies that 
$p_{2m+1}$ is nonzero, as the summands in this expression are linearly independent. 
\hfill \qed\medskip 

{\em End of proof of Proposition \ref{E:acyclic}.}
For $u\in \ZZ$, the complex $\cP_{0,1}^{\bullet}\{u\}$ is $0\to \cP_{0,1}^{0}[u]\to 
\cP_{0,1}^{1}[u+1]\to....$. For $u>0$, the groups of this complex are all zero, so 
$\cP_{0,1}^{\bullet}\{u\}$ is acyclic. 
For $u\leq 0$, this complex is $0\to...\to 0\to \cP_{0,1}^{2m}[m]\to 
\cP_{0,1}^{2m+1}[m+1]\to 0\to...$, where $m=-u$. The nontrivial map in this 
complex is $p_{2m}\mapsto \on{gr}d_{0,1}^{2m,2m+1}(p_{2m}) = p_{2m+1}$, 
which is an isomorphism, so $\cP_{0,1}^{\bullet}\{u\}$ is acyclic. It follows
that $\cP_{0,1}^{\bullet}$ is acyclic. As the differential of $\cP_{1,0}^{\bullet}$
is the negative of that of $\cP_{0,1}^{\bullet}$, $\cP_{1,0}^{\bullet}$ is acyclic as
well. 

Let $\eps\in\{0,1\}$ and $u\in\NN$. The complex $\cP_{\eps,\eps}^{\bullet}[u]$
is $0\to \cP_{\eps,\eps}^{0}[u]\to \cP_{\eps,\eps}^{1}[u]\to...$; 
if $i=u$, this complex is $0\to \kk\to 0\to 0 \to ...$, 
whose cohomology is $1$-dimensional, concentrated in degree $0$; if $u>0$, this complex
is $0\to...\to 0\to \cP^{2u-1}_{\eps,\eps}[u]\to \cP^{2u}_{\eps,\eps}[u]
\to 0\to 0...$; the nontrivial map in this complex is
$p_{2u-1}\mapsto \on{gr}d^{2u-1,2u}_{\eps,\eps}(p_{2u-1}) = up_{2u}$ if
$\eps=0$ and $-up_{2u}$ if $\eps=1$. As this is an isomorphism in both cases, 
$\cP^{\bullet}_{\eps,\eps}[u]$ is acyclic for $u>0$. It follows that the cohomology
of $\cP^{\bullet}_{\eps,\eps}$ is $1$-dimensional, concentrated in degree $0$. 

This implies that $\cE^{\bullet}_{\eps,\eps'}$ is acyclic for $\eps\neq\eps'$, 
and that the cohomology of $\cE^{\bullet}_{\eps,\eps}$ is concentrated in 
degree 0. As $d^{0,1}_{\eps,\eps}=0$, we have in degree 0, 
$H^{0}(\cE^{\bullet}_{\eps,\eps}) = \cE^{0}_{\eps,\eps}\simeq \kk$. 
\hfill \qed\medskip 

\begin{remark}
If $\G$ is a Lie algebra, we have natural maps 
\begin{equation} \label{maps}
H^{\bullet}(\cE^{\bullet}_{\eps,\eps'}) \to H^{\bullet}(\G,
U(\G)_{\eps,\eps'}).
\end{equation} When $(\eps,\eps')=(0,1)$ and $\G$ is finite dimensional, 
then $H^n(\G,U(\G)_{0,1}) = \kk$ if $n=\on{dim}\G$, and $=0$ otherwise. 
Indeed, if $C^n(\G) := \wedge^n(\G) \otimes U(\G)$, then the differential 
$d_\G^{n,n+1} : C^n(\G) \to C^{n+1}(\G)$ is given by 
$\omega\otimes x\mapsto \delta(\omega)\otimes x + 
\sum_{\alpha=1}^{\on{dim}\G}(\omega\wedge e^\alpha)\otimes (e_\alpha x)$, 
where $(e^\alpha)_\alpha,(e_\alpha)_\alpha$ are dual bases of $\G^*$
and $\G$ and $\delta:\wedge^n(\G^*)\to\wedge^{n+1}(\G^*)$ is induced by the 
Lie coalgebra structure of $\G^*$. For $i\in\ZZ$, set 
$F_i(C^n(\G)):= \wedge^n(\G^*) \otimes
U(\G)_{\leq n+i}$ (where the last term is the subspace of elements of 
degree $\leq n+i$ for the PBW filtration). Then $d^{n,n+1}(F_i(C^n(\G)))
\subset F_i(C^{n+1}(\G))$, so $...\subset F_i(C^n(\G))\subset ...\subset 
C^n(\G)$ is a complete filtration of $C^n(\G)$. The associated graded
complex is $\tilde C^n(\G):= \wedge^n(\G^*)\otimes S(\G)$, with 
differential $\tilde d^{n,n+1} : \tilde C^n(\G)\to \tilde C^{n+1}(\G)$, 
$\omega \otimes x\mapsto \sum_{\alpha=1}^{\on{dim}\G} (\omega\wedge e^\alpha)
\otimes (e_\alpha x)$. This complex only depends on the vector space
structure of $\G$; if we denote it by $\tilde C^\bullet(\G)$, then we have an 
isomorphism $\tilde C^\bullet(\G_1\oplus\G_2) \simeq \tilde C^\bullet(\G_1)
\otimes \tilde C^\bullet(\G_2)$, so  $C^\bullet(\G)\simeq 
C^\bullet(\kk)^{\otimes \on{dim}\G}$. As the cohomology of 
$\tilde C^\bullet(\kk)$ is 1-dimensional, concentrated in degree 1, 
the cohomology of $\tilde C^\bullet(\G)$ is 1-dimensional, concentrated in 
degree $\on{dim}\G$. It follows that $C^\bullet(\G)$ is acyclic in 
every degree $\neq \on{dim}\G$, and its cohomology in degree 
$\on{dim}\G$ has dimension $\leq 1$. If $\omega\in \wedge^{\on{dim}\G}(\G)$
is nonzero, then $\omega\otimes 1\in C^n(\G)$ is a nontrivial cocycle, 
so the cohomology of $C^\bullet(\G)$ coincides with that of 
$\tilde C^\bullet(\G)$. 
As $U(\G)_{1,0}\simeq U(\G)_{0,1}$ (using the antipode), we have 
$H^{\bullet}(\G,U(\G)_{1,0}) \simeq H^{\bullet}(\G,U(\G)_{0,1})$.


When $\eps\neq\eps'$, (\ref{maps}) is the zero map. If $\eps=\eps'$, 
then the map $\kk = H^{0}(\cE^{\bullet}_{\eps,\eps}) \to 
H^{0}(\G,U(\G)_{\eps,\eps})$ takes $1$ to the class of 
$1\in U(\G)_{\eps,\eps}$ (which is invariant, both under 
the trivial and the adjoint actions of $\G$ on $U(\G)$). \hfill \qed\medskip  
\end{remark}

{\em End of proof of Theorem \ref{thm:vanishing}.} 
One of the pairs $(0,\eps_{2}),(\eps_{2},\eps_{3}),...,(\eps_{z+N},1)$
necessarily coincides with $(0,1)$; call it $(\eps_{i},\eps_{i+1})$. 
According to Proposition \ref{E:acyclic}, the corresponding complex
$\cE^{\bullet}_{\eps_{i},\eps_{i+1}}$ is then acyclic. 
Lemma \ref{A:E} and the K\"unneth formula then imply 
that $\cA^{\bullet}_{\sigma}$ is acyclic. 
This being valid for any $\sigma$, the decomposition (\ref{decomp:A})
then implies that $\cA^{\bullet}_{z,N,1}$ is acyclic, as claimed. 
\hfill \qed\medskip 

\section{Compatibility of quantization functors with twists}

We first prove the compatibility of quantization functors of quasi-Lie 
bialgebras with twists of quasi-Lie bialgebras; we derive from there the compatibility
of quantization functors of Lie bialgebras with twists of Lie bialgebra (a result which 
was obtained in \cite{EH} in the case of Etingof-Kazhdan quantization functors). 

\subsection{Twists of quasi-Lie bialgebras}

Let $\on{QLBA}_{f}$ be the prop with the same generators as $\on{QLBA}$
with the additional $f\in \on{QLBA}_{f}(\wedge^{2},{\bf id})$ and the same relations. 
This prop is $\NN$-graded, if we extend the degree $\on{deg}_{\mu}+\on{deg}_{\delta}$
in $\on{QLBA}$ by $|f|=1$. 

We then have $\on{QLBA}_{f}(X,Y) = \on{QLBA}(S(\wedge^{2}) \otimes X,Y)$. 
The filtration of $\on{QLBA}_{f}$ induced by the degree is such that 
$$\on{QLBA}_{f}^{\geq n}(X,Y) = \oplus_{k\geq 0}
\on{QLBA}^{\geq n-k}(S^{k}(\wedge^{2})\otimes X,Y).$$ 
It follows that $\on{QLBA}_{f}(X,Y) \subset \on{QLBA}_{f}^{\geq v_{f}(|X|,|Y|)}
(X,Y)$, where $v_{f}(|X|,|Y|) = \on{inf}\{v(|X|+2k,|Y|)+k, k\geq 0\}$ and 
$v(|X|,|Y|)={1\over 3}||X|-|Y||$. 
As $v_{f}(|X|,|Y|)\geq v(|X|,|Y|)$, $\on{QLBA}_{f}$ gives rise to a 
topological prop ${\bf QLBA}_{f}$. 

We have two prop morphisms $\kappa_{i} : \on{QLBA}\to \on{QLBA}_{f}$, 
defined by 
$$
\kappa_{1} : \mu,\delta,\varphi\mapsto\mu,\delta,\varphi, 
$$
$$
\kappa_{2}:\mu\mapsto\mu, \delta\mapsto \delta + \on{Alt}_{2}\circ 
(\on{id}_{\bf id}\otimes\mu)\circ (f\otimes \on{id}_{\bf id}),
$$
$$
\varphi\mapsto 
\varphi + {1\over 2}\on{Alt}_{3}\circ [(\delta\otimes\on{id}_{\bf id})
\circ f + (\on{id}_{\bf id}\otimes 
\mu\otimes\on{id}_{\bf id})\circ (f\otimes f)];  
$$
this is the universal version of the operation of twisting of a quasi-Lie 
bialgebra structure. The prop morphisms $\kappa_{i}$ extend to topological props. 

Let $(m,\Delta,\Phi,\eps,\eta)$ be a QSB in ${\bf QLBA}$ (i.e., a quantization 
functor of quasi-Lie bialgebras). 
For $i=1,2$, set $(m_{i},\Delta_{i},\Phi_{i},\eta_{i},\eta_{i})
:= \kappa_{i}(m,\Delta,\Phi,\eps,\eta)$. 
Then $(m_{i},\Delta_{i},\Phi_{i},\eps_{i},\eta_{i})$ are
QSQB's in ${\bf QLBA}_{f}$. 

\begin{prop} \label{prop:twists}
The QSQB's $(m_{i},\Delta_{i},\Phi_{i},\eps_{i},\eta_{i})$ are related by 
equivalence and twist, i.e., there exists 
$F\in {\bf QLBA}_{f}({\bf 1},S^{\otimes 2})$ and 
$i\in {\bf QLBA}_{f}(S,S)$ such that $(\eps\otimes\on{id}_{S}) \circ F = 
(\on{id}_{S}\otimes \eps)\circ F = \eta$, $F = inj_{0}^{\otimes 2}$
and $i=\on{id}_{S}$ mod $\on{QLBA}_{f}^{\geq 1}$, such that 
$(m_{2},\Delta_{2},\Phi_{2},\eps_{2},\eta_{2}) = i *\{ F\star
(m_{1},\Delta_{1},\Phi_{1},\eps_{1},\eta_{1})\}$. 
\end{prop}

This implies that the quantization functors of quasi-Lie bialgebras take
a pair of quasi-Lie bialgebras related by a classical twist to a pair of 
quasi-Hopf QUE algebras related by a quantum twist. 

{\em Proof.} As $(\mu,\delta,\varphi)$ defines a quasi-Lie bialgebra structure on
${\bf id}$ on $\on{QLBA}_{f}$, we have $[\mu\oplus\delta\oplus\varphi,
\mu\oplus\delta\oplus\varphi]=0$ (in $C_{\on{QLBA}_{f}}$), 
so the bracket with $\mu\oplus\delta\oplus\varphi$ defines a complex 
structure on $C_{\on{QLBA}_{f}}$. This gives rise to cohomology groups
$H^{i}_{\on{QLBA}_{f}}$, which are graded by the grading of $\on{QLBA}_{f}$. 

\begin{lemma} \label{lemma:coh}
The canonical maps $H^{i}_{\on{QLBA}_{f}}\to H^{i}_{\on{LBA}}$ are 
isomorphisms. 
\end{lemma}

{\em Proof of Lemma.} As before, we will show that the relative complex
is acyclic. This relative complex is filtered by the powers of the ideal 
$\langle\varphi\rangle_{f}$ of $\on{QLBA}_{f}$. The associated
graded prop of $\on{QLBA}_{f}$ w.r.t. this filtration is $\on{LBA}_{\alpha,f}$
defined by $\on{LBA}_{\alpha,f}(X,Y) = \on{LBA}_{\alpha}(S(\wedge^{2})
\otimes X,Y) = \on{Coker}[\on{LBA}(S(\wedge^{2})\otimes S(\wedge^{3})
\otimes \wedge^{4}\otimes
X,Y)\to \on{LBA}(S(\wedge^{2})\otimes S(\wedge^{3})\otimes X,Y)]$, 
induced by the morphism $S(\wedge^{2})\otimes S(\wedge^{3}) \to 
S(\wedge^{2}) \otimes S(\wedge^{3}) \otimes \wedge^{4}$, 
tensor product with $\on{id}_{S(\wedge^{2})}$ of 
$S(\wedge^{3}) \to S(\wedge^{3})^{\otimes 2} \to \wedge^{3}\otimes 
S(\wedge^{3}) \stackrel{[\on{Alt}_{4} \circ (\delta\otimes
\on{id}_{{\bf id}}^{\otimes 2})]\otimes \on{id}_{S(\wedge^{3})}}{\to}
\wedge^{4}\otimes S(\wedge^{3})$ (the first morphism is the coproduct in $S$, 
the second morphism is the projection $S(\wedge^{3})\to\wedge^{3}$). 
The associated graded complex of the relative complex is the positive degree
part of $C_{\on{LBA}_{\alpha,f}}$, w.r.t. the degree on $\on{LBA}_{\alpha,f}$
induced by the filtration, i.e., for $k\geq 0$, the degree $k$ part of 
$\on{LBA}_{\alpha,f}(X,Y)$ is 
$\on{Coker}[\on{LBA}(S(\wedge^{2})\otimes S^{k-1}(\wedge^{3})
\otimes \wedge^{4}\otimes X,Y)\to \on{LBA}(S(\wedge^{2})\otimes 
S^{k}(\wedge^{3})\otimes X,Y)]$; by convention $S^{-1}(\wedge^{3})
=0$. 

For $k>0$, the degree $k$ part of the associated graded complex 
has the form $C^{p,q}_{\on{QLBA}_{f}}\{k\} = \on{Coker}[\on{LBA}
(S(\wedge^{2})\otimes S^{k-1}(\wedge^{3})\otimes \wedge^{4}
\otimes \wedge^{p+1},\wedge^{q+1})
\to \on{LBA}(S(\wedge^{2})\otimes S^{k}(\wedge^{3})\otimes 
\wedge^{p+1},\wedge^{q+1})]$ ($p,q\geq -1$), equipped with the differential 
$[\mu\oplus\delta,-]$. For $k>0$, both 
$S(\wedge^{2})\otimes S^{k-1}(\wedge^{3})\otimes \wedge^{4}$ and 
$S(\wedge^{2})\otimes S^{k}(\wedge^{3})$ are sums of Schur functors
of positive degrees, so Theorem \ref{thm:vanishing} implies 
that the lines of $C^{p,q}_{\on{QLBA}_{f}}\{k\}$ are acyclic, so for each $k>0$,
this complex is acyclic. So the relative complex 
$\on{Ker}(C_{\on{QLBA}_{f}}\to C_{\on{LBA}})$ is acyclic. 
 \hfill \qed\medskip 

{\em End of proof of Proposition \ref{prop:twists}.}
Recall that we have a prop morphism 
$\pi : \on{QLBA}\to \on{LBA}$ defined by 
$\pi: \mu,\delta,\varphi\mapsto \mu,\delta,0$. 
We also have a prop morphism $\pi_{f}:\on{QLBA}_{f}\to \on{LBA}$, 
defined by $\pi_{f} : \mu,\delta,\varphi,f\mapsto \mu,\delta,0,0$. 
These morphisms extend to topological props and satisfy 
$\pi_{f}\circ\kappa_{i}=\kappa$ for $i=1,2$. 

Then $\pi_{f}(m_{i},\Delta_{i},\Phi_{i},\eps_{i},\eta_{i}) 
= \pi_{f}\circ\kappa_{i}
(m,\Delta,\Phi,\eps,\eta) = \pi(m,\Delta,\Phi,\eps,\eta)$ for $i=1,2$. 
The classical limits of $(m_{i},\Delta_{i},\Phi_{i},\eps_{i},\eta_{i})$
are $(\mu,\delta,Z)$ for $i=1,2$ (as the additional terms $\on{Alt}_{2}
\circ (\on{id}_{\bf id}\otimes\mu) \circ (f\otimes \on{id}_{\bf id})$, 
${1\over 2}\on{Alt}_{3} \circ [(\delta\otimes\on{id}_{\bf id})
\circ f + (\on{id}_{\bf id}\otimes\mu\otimes\on{id}_{\bf id})\circ
(f\otimes f)]$ have degree $>1$). 

As in the proof of Theorem \ref{thm:main} (Subsection 
\ref{pf:main}), Lemma \ref{lemma:coh} 
implies that two QSQB's in ${\bf QLBA}_{f}$ with the same image 
in ${\bf LBA}$ and the same classical limit 
are related by twist and equivalence. It follows that 
$(m_{i},\Delta_{i},\Phi_{i},\eps_{i},\eta_{i})$ ($i=1,2$) 
are related by twist and equivalence. \hfill \qed\medskip 

\subsection{Twists of Lie bialgebras} In \cite{EH}, we introduced the prop 
$\on{LBA}_{f}$ of Lie bialgebras with a twist; its generators are $\mu,\delta,f
\in \on{LBA}_{f}({\bf 1},\wedge^{2})$
and its relations as the same as the relations of $\on{LBA}$, together with 
$\on{Alt}_{3}\circ [(\delta\otimes\on{id}_{\bf id})\circ f + 
(\on{id}_{\bf id}\otimes\mu\otimes\on{id}_{\bf id})\circ (f\otimes f)]=0$. 
We defined prop morphisms $\kappa_{i}^{0} : \on{LBA}\to \on{LBA}_{f}$ by 
$\kappa_{1} : \mu,\delta\mapsto\mu,\delta$ and $\kappa_{2} : 
\mu,\delta\mapsto \mu,\delta + \on{Alt}_{2}\circ (\on{id}_{\bf id}
\otimes\mu)\circ (f\otimes\on{id}_{\bf id})$. 

Let $(m_{0},\Delta_{0},\eps_{0},\eta_{0})$ be a QSB in ${\bf LBA}$, 
quantizing $(\mu,\delta)$ (i.e., a quantization functor of Lie bialgebras).
Set $(m_{i}^{0},\Delta_{i}^{0},\eps_{i}^{0},\eta_{i}^{0}):= 
\kappa_{i}^{0}(m_{0},\Delta_{0},\eps_{0},\eta_{0})$. These are QSB's 
in ${\bf LBA}_{f}$. 

\begin{prop} \label{prop:5:2}
The QSB's $(m_{i}^{0},\Delta_{i}^{0},\eps_{i}^{0},\eta_{i}^{0})$
are related by a and equivalence and a bialgebra twist. More explicitly, 
there exists $F_{0}\in {\bf LBA}_{f}({\bf 1},S^{\otimes 2})$ and $i_{0}
\in {\bf LBA}_{f}(S,S)$, whose reduction mod ${\bf LBA}_{f}^{\geq 1}$
is $inj_{0}^{\otimes 2}$ and $\on{id}_{S}$, such that $(\eps_{0}\otimes
\on{id}_{S})\circ F_{0} = (\on{id}_{S}\otimes \eps_{0})\circ F_{0}
= \eta_{0}$, $[\eta_{0}\otimes F_{0}] * [(\on{id}_{S}\otimes
\Delta_{0}) \circ F_{0}] = [F_{0}\otimes \eta_{0}] * [(\Delta_{0}
\otimes \on{id}_{S}) \circ F_{0}]$, and 
$(m_{2}^{0},\Delta_{2}^{0},\eps_{2}^{0},\eta_{2}^{0}) = i_{0} * 
(m_{1}^{0}, F_{0} * \Delta_{1}^{0} * F_{0}^{-1},\eps_{1}^{0},\eta_{1}^{0})$. 
\end{prop}

{\em Proof.} There are uniquely defined prop morphisms 
$\tilde\pi_{f} : \on{QLBA}_{f}\to \on{LBA}_{f}$ and 
$\pi_{f}^{0} : \on{LBA}_{f}\to \on{LBA}$, such that
$\tilde\pi_{f} : \mu,\delta,\varphi,f\mapsto \mu,\delta,0,f$
and $\pi_{f}^{0} : \mu,\delta,f\mapsto \mu,\delta,0$. 
We then have $\on{id}_{\on{LBA}} = \pi_{f}^{0}\circ \kappa_{i}^{0}$
for $i=1,2$, $\kappa_{i}^{0}\circ \pi = \tilde\pi_{f}\circ \kappa_{i}$, 
$\pi_{f} = \pi_{f}^{0}\circ \tilde\pi_{f}$. All the information on 
compositions can be summarized in the commutative diagram 
$$
\begin{matrix}
\on{QLBA} & \stackrel{\kappa_{i}}{\to} & \on{QLBA}_{f} & 
\searrow\scriptstyle{\pi_{f}} \\
\downarrow\scriptstyle{\pi} & & \downarrow
\scriptstyle{\tilde\pi_{f}}
& \on{LBA} \\
\on{LBA} & \stackrel{\kappa_{i}^{0}}{\to} & \on{LBA}_{f} & \nearrow
\scriptstyle{\pi_{f}^{0}} 
\end{matrix}
$$
together with the relations $\pi = \pi_{f}\circ \kappa_{1,2}$, $\on{id}_{\on{LBA}}
= \pi_{f}^{0}\circ \kappa_{1,2}^{0}$. These morphisms induce morphisms 
between completed props.

According to Theorem \ref{thm:main}, we may lift 
$(m_{0},\Delta_{0},\eps_{0},\eta_{0})$ to a QSQB
$(m,\Delta,\Phi,\eps,\eta)$ in ${\bf QLBA}$ with classical limit
$(\mu,\delta,\varphi)$; this means that 
$\pi(m,\Delta,\Phi,\eps,\eta)= (m_{0},\Delta_{0},\eta_{0}^{\otimes 3},
\eps_{0},\eta_{0})$. Set 
$$(m_{i},\Delta_{i},\Phi_{i},\eps_{i},\eta_{i}):= 
\kappa_{i}(m,\Delta,\Phi,\eps,\eta),$$ then 
$\pi(m_{i},\Delta_{i},\Phi_{i},\eps_{i},\eta_{i}) = 
(m_{i}^{0},\Delta_{i}^{0},\eta_{i}^{\otimes 3},\eps_{i},\eta_{i})$. 

Let then $i,F$ be as in Proposition \ref{prop:twists}. Let $i_{0}:= \tilde\pi_{f}(i)$, 
$F_{0}:= \tilde\pi_{f}(F)$. As $F$ is a twist relating $\Phi_{1}$ and $i^{-1}*\Phi_{2}$, 
and as the images of $\Phi_{1},i^{-1}(\Phi_{2})$ 
under $\tilde\pi_{f}$ are $\eta_{1}^{\otimes 3}$, 
$F_{0}$ is a twist relating $\eta_{1}^{\otimes 3}$ with itself, i.e., 
it satisfies the announced cocycle relation. The image under $\tilde\pi_{f}$
of the statement that the $(m_{i},\Delta_{i},\Phi_{i},\eps_{i},\eta_{i})$
are related by $(i,F)$ is then that the $(m_{i}^{0},\Delta_{i}^{0},\eps_{i}^{0},
\eta_{i}^{0})$ are related by $(i_{0},F_{0})$. \hfill \qed\medskip

\begin{appendix} 
\section{Structure of the prop $\on{LBA}$} \label{app:A}

The following structure theorem of the prop LBA was proved in \cite{Enr:univ,Pos}. 
We reformulate here this proof using the language of props. In \cite{EH}, 
we derived Proposition \ref{prop:PBW} from Theorem \ref{thm:PBW} below. 

\begin{thm} \label{thm:PBW}
If $F,G\in \on{Ob(Sch)}$, then the map $\oplus_{Z\in\on{Irr(Sch)}}
\on{LCA}(F,Z)\otimes \on{LA}(Z,G)\to \on{LBA}(F,G)$ 
induced by composition and the 
prop morphisms $\on{LCA}\to \on{LA}$, $\on{LBA}\to \on{LA}$
is a linear isomorphism. 
\end{thm}

{\em Proof.} It suffices to prove this when $F,G\in \on{Irr(Sch)}$, and then 
(using the action of $\SG_{n},\SG_{m}$) for $F = T_{n}$, $G = T_{m}$. 
Using the cocycle relation, and the isomorphism of the l.h.s. with 
$\oplus_{z\geq 0}
(\on{LCA}(T_{n},T_{z}) \otimes \on{LA}(T_{z},T_{m}))_{\SG_{z}}$, 
one proves that the morphism is surjective. We now prove that it is injective. 
We have $$\on{LBA}(T_{n},T_{m}) = \oplus_{a,b\geq 0|a-b=n-m}
\on{LBA}(T_{n},T_{m})[a,b] = \oplus_{z\geq \on{min}(n,m)}
\on{LBA}(T_{n},T_{m})[z-m,z-n],$$ and the morphism is the direct of 
over $z\geq \on{min}(n,m)$ of the maps
$\oplus_{Z\in \on{Irr(Sch)}||Z|=z}
\on{LCA}(T_{n},Z)\otimes \on{LA}(Z,T_{m})
\to \on{LBA}(T_{n},T_{m})[z-m,z-n]$. It remains to show that each of the 
maps is injective. 

There is a unique morphism $\on{LBA}\to L({\bf LCA})$ (the generators of 
$\on{LBA}$ are $\mu,\delta$, the generator of $\on{LCA}$ is $\delta_{\on{LCA}}$), 
taking $\mu$ to 
$\mu_{free} : L^{\otimes 2}\to L$ and $\delta$ to the unique $\delta_{free} : 
L\to L^{\otimes 2}$, such that ${\bf id} \to L \stackrel{\delta_{free}}{\to} 
L^{\otimes 2}$ is  $\delta_{\on{LCA}}:{\bf id}
\stackrel{\delta_{\on{LCA}}}{\to} 
{\bf id}^{\otimes 2}\to L^{\otimes 2}$ and 
$\delta\circ \mu_{free} = ((\mu_{free}\otimes \on{id}_{L})
\circ (\on{id}_{L}\otimes \beta_{L,L})+\on{id}_{L}\otimes \mu_{free})
\circ (\delta_{free}\otimes \on{id}_{L}) + (\mu_{free}\otimes \on{id}_{L}
+ (\on{id}_{L}\otimes\mu_{free}) \circ (\beta_{L,L}\otimes \on{id}_{L})) 
\circ (\on{id}_{L}\otimes \delta_{free})$. The prop $\on{LCA}$ is $\ZZ$-graded, 
with $\on{deg}\delta_{\on{LCA}}=1$, then the morphism $\on{LBA}\to 
L({\bf LCA})$ is compatible with the morphism $\ZZ^{2}\to\ZZ$, $(1,0)\mapsto 0$, 
$(0,1)\mapsto 1$. 

We then have maps $\on{LBA}(T_{n},T_{m})\to L({\bf LCA})(T_{n},T_{m}) = 
{\bf LCA}(L^{\otimes n},L^{\otimes m}) \to {\bf LCA}(T_{n},L^{\otimes m})$, 
where the last map is induced by ${\bf id}\to L$, which restrict to 
$\on{LBA}(T_{n},T_{m})[z-m,z-n]\to {\bf LCA}(T_{n},L^{\otimes m})[z-n]
= \on{LCA}(T_{n},(L^{\otimes m})_{z}) 
$, where the index $z$ denotes the 
(Schur functor) degree $z$ part. 

\begin{lemma} \label{lemma:app}
If $X$ is any prop and $F\in \on{Ob(Sch)}$, we have an isomorphism 
$X(F,(L^{\otimes m})_{z}) \simeq \oplus_{Z\in\on{Irr(Sch)},|Z|=z}X(F,Z)\otimes 
\on{LA}(Z,T_{m})$. 
\end{lemma}

{\em Proof of Lemma.}
We have isomorphisms $\on{LA}(T_{z},{\bf id}) \simeq$ multilinear part
of the free Lie algebra in $z$ ordered generators $\simeq \on{Sch}(T_{z},L_{z})$. 
So if $|Z|=z$, $\on{LA}(Z,{\bf id}) \simeq \on{Sch}(Z,L_{z})$, which 
may be expressed as $L_{z} = \oplus_{|Z|=z} \on{LA}(Z,{\bf id}) \otimes Z$. 

So 
\begin{align*}
& X(F,(L^{\otimes m})_{z}) = \oplus_{z_{1}+...+z_{m}=z}
X(F,\otimes_{i=1}^{m}
L_{z_{i}}) 
\\ &  = \oplus_{|Z_{1}|+...+|Z_{m}|=z} X(F,\otimes_{i=1}^{m}Z_{i})
\otimes (\otimes_{i=1}^{m}\on{LA}(Z_{i},{\bf id}))
\\ & = \oplus_{|Z_{1}|+...+|Z_{m}|=z,|Z|=z} X(F,Z) \otimes \on{Sch}(Z,
\otimes_{i=1}^{m}Z_{i})
\otimes (\otimes_{i=1}^{m}\on{LA}(Z_{i},{\bf id}))
\\ & = \oplus_{|Z|=z}X(F,Z)\otimes \on{LA}(Z,T_{m}),
\end{align*} where 
the last equality follows from  $\on{LA}(Z,T_{m}) = 
\oplus_{Z_{1},...,Z_{m}\in \on{Irr(Sch)}|\sum_{i}|Z_{i}|=z}
\on{Sch}(Z,\otimes_{i}Z_{i}) \otimes \otimes_{i}\on{LA}(Z_{i},{\bf id})$, 
for $Z\in \on{Ob(Sch)}$ (see \cite{EH}). 
\hfill \qed \medskip 

{\em End of proof of Theorem.} We have constructed a map 
$\on{LBA}(T_{n},T_{m})[z-m,z-n]\to \oplus_{Z\in\on{Irr(Sch)},|Z|=z}
\on{LCA}(T_{n},Z)\otimes \on{LA}(Z,T_{m})$, and one proves that is a section of 
the morphism $\oplus_{Z\in \on{Irr(Sch)}||Z|=z}
\on{LCA}(T_{n},Z)\otimes \on{LA}(Z,T_{m})
\to \on{LBA}(T_{n},T_{m})[z-m,z-n]$, which is therefore 
injective. \hfill \qed\medskip 

\end{appendix}

\end{document}